\documentclass[12pt]{article}
\usepackage{fullpage}
\usepackage{graphicx}
\usepackage{epic}
\usepackage{eepic}
\usepackage{amsmath,amssymb,mathrsfs}
\usepackage[francais]{babel}
\usepackage[T1]{fontenc}
\usepackage[latin1]{inputenc}
\usepackage{soul}
\usepackage{pst-eucl}
\usepackage[all, ps]{xy}
\newtheorem{theom}{Théorème}
\newtheorem{lem}{Lemme}[subsection]
\newtheorem{duf}{Definition}[subsection]
\newtheorem{theo}{Théorème}[subsection]
\newtheorem{prop}{Proposition}[subsection]

\author{Raphaël Beuzart-Plessis}
\title{Expression d'un facteur epsilon de paire par une formule intégrale}
\begin{document}

\maketitle

\bigskip

\textbf{Introduction}

\vspace{4mm}

Soient $F$ un corps local non archimédien de caractéristique nulle, $E$ une extension quadratique de $F$ et $c$ le $F$-automorphisme non trivial de $E$. Fixons deux entiers naturels $m<d$ de parités différentes et deux représentations lisses irréductibles tempérées $\pi$, $\sigma$ de $GL_d(E)$ et $GL_m(E)$ respectivement. On suppose les conditions suivantes vérifiées

\begin{center}
$\pi\simeq \pi^\vee\circ c$, $\sigma\simeq \sigma^\vee\circ c$
\end{center}

\noindent où $\pi^\vee$ (resp. $\sigma^\vee$) désigne la représentation contragrédiente de $\pi$ (resp. de $\sigma$) et où l'action de $c$ sur $E$ est naturellement étendue à $GL_d(E)$ et $GL_m(E)$. De telles représentations seront qualifiées de conjuguées-duales. Jacquet, Piatetskii-Shapiro et Shalika ont défini dans [JPSS] une fonction holomorphe $s\mapsto \epsilon(s,\pi\times \sigma,\psi_E)$. Cette fonction epsilon dépend d'un caractère continu non trivial $\psi_E$ de $E$. On fixe un caractère continu non trivial $\psi$ de $F$ et on pose $\psi_E=\psi\circ Tr$ où $Tr$ est l'application trace de l'extension $E/F$. Le résultat principal de cet article est une formule intégrale exprimant $\epsilon(1/2,\pi\times \sigma,\psi_E)$. Une telle formule est démontrée ici en vue d'une application à la conjecture locale de Gan-Gross-Prasad pour les groupes unitaires ([GGP] conjecture 17.3). On a montré en [B], l'existence d'une formule similaire pour une multiplicité apparaissant dans cette conjecture, qui prédit un lien entre cette multiplicité et certains facteurs epsilon de paires comme ci-dessus. L'existence d'une telle formule n'est donc pas étonnante, en tout cas si l'on croit à la conjecture 17.3 de [GGP]. Remarquons que la formule montrée ici est aussi l'exacte analogue de celle obtenue par Waldspurger ([W3]), dans le cadre de la conjecture de Gan-Gross-Prasad pour les groupes spéciaux orthogonaux (il s'agissait alors de facteurs epsilon de paires pour des représentations autoduales de groupes linéaires). \\

 Enonçons maintenant le résultat principal de cet article. Notons $J_d\in GL_d(E)$ la matrice de coefficients $(J_d)_{i,j}=(-1)^{i}\delta_{i,d-j}$ pour tous $i,j=1,\ldots,d$ ($\delta_{i,j}$ le symbole de Kronecker). Notons $\theta_d$ l'automorphisme de $GL_d(E)$ donné par $g\mapsto J_d {}^t c(g)^{-1} J_d^{-1}$. Puisque $\pi$ est conjuguée-duale, on peut l'étendre en une représentation du groupe non connexe $GL_d^+(E)=GL_d(E)\rtimes \{1,\theta_d\}$. En fait, il y a deux extensions possibles. La théorie des fonctionnelles de Whittaker permet d'en privilégier une. Rappelons que l'on a fixé un caractère non trivial $\psi_E$ de $E$. Soit $U_d$ le sous-groupe de $GL_d$ des matrices unipotentes triangulaires supérieures. On définit un caractère $\xi$ de $U_d(E)$ par la formule

$$\displaystyle \xi(u)=\psi_E\big(\sum_{i=1}^{d-1}u_{i,i+1}\big)$$

\noindent pour tout $u\in U_d(F)$ et où les $u_{i,j}$ désignent les coefficients de la matrice $u$. Appelons fonctionnelle de Whittaker pour $\pi$ toute forme linéaire $\ell$ sur $E_\pi$ (l'espace de la représentation $\pi$) vérifiant $\ell\circ \pi(u)=\xi(u)\ell$ pour tout $u\in U_d(E)$. Puisque $\pi$ est tempérée, l'espace des fonctionnelles de Whittaker est une droite. Soit $\pi^+$ un prolongement de $\pi$ à $GL_d^+(E)$, l'application linéaire $\ell\mapsto \ell\circ \pi^+(\theta_d)$ préserve la droite des fonctionnelles de Whittaker et il existe un unique prolongement de sorte qu'elle y induise l'identité. C'est ce prolongement que l'on privilégie et on note $\tilde{\pi}$ la restriction de ce dernier à $\tilde{GL}_d(E)=GL_d(E)\theta_d$. On définit de la même façon $\tilde{\sigma}$. Soit $V$ un $E$-espace vectoriel de dimension $d$. Posons $G=R_{E/F}GL(V)$ ($R_{E/F}$ désigne la restriction des scalaires à la Weil) et $\tilde{G}$ la variété algébrique (définie sur $F$) des formes sesquilinéaires non dégénérées sur $V$ (linéaires en la deuxième variable). La variété $\tilde{G}$ est munie de deux actions commutantes à droite et à gauche de $G$. On peut identifier, au moyen d'une base, $GL_d(E)$ à $G(F)$ et $\tilde{GL}_d(E)$ à $\tilde{G}(F)$. On fixe de même un $E$-espace vectoriel $W$ de dimension $m$, on pose $H=R_{E/F}GL(W)$, $\tilde{H}$ la variété algébrique des formes sesquilinéaires non dégénérées sur $W$ et on identifie $GL_m(E)$ à $H(F)$ et $\tilde{GL}_m(E)$ à $\tilde{H}(F)$. Fixons aussi une injection de $W$ dans $V$, une décomposition $V=W\oplus Z$ et une forme hermitienne $\zeta$ sur $Z$ somme orthogonale de plans hyperboliques et de la forme hermitienne de dimension $1$, $(x,y)\mapsto 2\nu c(x)y$ pour un certain $\nu\in F^\times$. On identifie alors $H$ au sous-groupe des éléments de $G$ qui agissent trivialement sur $Z$ et on plonge $\tilde{H}$ dans $\tilde{G}$ en envoyant toute forme sesquilinéaire non dégénérée $\tilde{x}$ sur $W$ sur la somme de $\tilde{x}$ et de $\zeta$. On doit, pour énoncer la formule intégrale introduire un ensemble de "sous-tores" tordus de $\tilde{H}$. Soient $W=W'\oplus W''$ une décomposition et $\zeta''$ une forme hermitienne sur $W''$ telle que les groupes unitaires de $\zeta''$ et $\zeta''\oplus \zeta$ soient quasi-déployés. Notons $\tilde{H}'$ la variété des formes sesquilinéaires non dégénérées sur $W'$ qui est munie de deux actions (à droite et à gauche) de $H'=R_{E/F} GL(W')$. On appelle sous-tore tordu maximal de $\tilde{H}'$ tout couple $(T',\tilde{T}')$ où $T'$ est un sous-tore maximal de $H'$ défini sur $F$ et $\tilde{T}'$ est une sous-variété définie sur $F$ de $\tilde{H}'$ qui est le normalisateur de $T'$ et d'un sous-groupe de Borel (pas forcément défini sur $F$) $B'$ contenant $T'$. Soit $(T',\tilde{T}')$ un sous-tore tordu maximal de $\tilde{H}'$, alors il existe un automorphisme $\theta$ de $T'$ tel que $\tilde{x}t'=\theta(t')\tilde{x}$ pour tous $\tilde{x}\in \tilde{T}'$, $t'\in T'$. On note $T'_\theta$ la composante neutre du sous-groupe fixé par $\theta$ de $T'$ et on dit que $(T',\tilde{T}')$ est anisotrope si $T'_\theta$ l'est. Soit $(T',\tilde{T}')$ un sous-tore tordu maximal anisotrope de $\tilde{H}'$ et notons $\tilde{T}$ la sous-variété des éléments de $\tilde{H}$ qui s'écrivent comme somme directe d'un élément de $\tilde{T}'$ et de $\zeta''$. On note $\underline{\mathcal{T}}$ l'ensemble des "sous-tores" tordus $\tilde{T}$ de $\tilde{H}$ obtenus de cette façon. Soit $\tilde{T}\in \underline{\mathcal{T}}$, la décomposition $W=W'\oplus W''$, la forme hermitienne $\zeta''$ et le tore tordu $(T',\tilde{T}')$ sont uniquement déterminés. Identifions $T'$ à un sous-groupe de $H$ en laissant agir $T'$ trivialement sur $W''$. Alors l'action par conjugaison de $T'$ sur $\tilde{T}$ définit un automorphisme $\theta$ de $\tilde{T}$ et on note $\tilde{T}(F)/\theta$ l'ensemble des orbites dans $\tilde{T}(F)$. Munissons $T'_\theta(F)$ de la mesure de Haar de masse totale $1$. On peut munir $\tilde{T}(F)/\theta$ d'une structure de $F$-espace analytique et d'une mesure caractérisées de la façon suivante: pour tout $\tilde{t}\in \tilde{T}/\theta$, l'application

$$T'_\theta(F)\to \tilde{T}(F)/\theta$$
$$t'\mapsto t'\tilde{t}$$

\noindent est un isomorphisme local qui préserve localement les mesures. Notons $Norm_H(\tilde{T})$ le normalisateur dans $H$ de $\tilde{T}$. Il contient naturellement $U(\zeta'')\times T'$ comme sous-groupe algébrique (où $U(\zeta'')$ est le groupe unitaire de $\zeta''$). On pose

$$W(H,\tilde{T})=Norm_H(\tilde{T})(F)/\big(U(\zeta'')(F)\times T'(F)\big)$$

\noindent c'est un groupe fini. Pour $\tilde{t}\in \tilde{T}(F)$, on a un développement local

$$\displaystyle \Theta_{\tilde{\pi}}(\tilde{t}exp(X))=\sum_{\mathcal{O}\in Nil(\mathfrak{g}_{\tilde{x}})} c_{\tilde{\pi},\mathcal{O}}(\tilde{t}) \hat{j}(\mathcal{O},X)$$

\noindent pour tout $X$ assez proche de $0$ dans $\mathfrak{g}_{\tilde{x}}(F)$, où $Nil(\mathfrak{g}_{\tilde{x}})$ désigne l'ensemble des orbites nilpotentes dans $\mathfrak{g}_{\tilde{x}}(F)$ et $\hat{j}(\mathcal{O},.)$ est la transformée de Fourier de l'intégrale orbitale suivant $\mathcal{O}$ (il faut pour cela préciser la transformée de Fourier et les mesures, on renvoie pour cela au corps du texte). On pose alors

$$\displaystyle c_{\tilde{\pi}}(\tilde{t})=\frac{1}{|Nil_{reg}(\mathfrak{g}_{\tilde{x}})|}\sum_{\mathcal{O}\in Nil_{reg}(\mathfrak{g}_{\tilde{x}})} c_{\tilde{\pi},\mathcal{O}}(\tilde{t})$$

\noindent où $Nil_{reg}(\mathfrak{g}_{\tilde{x}})$ désigne l'ensemble des orbites nilpotentes régulières de $\mathfrak{g}_{\tilde{x}}(F)$. On définit de façon analogue $c_{\tilde{\sigma}}(\tilde{t})$. Fixons un ensemble $\mathcal{T}$ de représentants des classes de conjugaison de $\underline{\mathcal{T}}$ par $H(F)$ et posons

\[\begin{aligned}
\displaystyle & \epsilon_{geom,\nu}(\tilde{\pi},\tilde{\sigma})=\\ 
 & \sum_{\tilde{T}\in\mathcal{T}} |2|_F^{2r^2+r+2rdim(W'')} |W(H,\tilde{T})|^{-1} \lim\limits_{s\to 0^+} \int_{\tilde{T}(F)/\theta} c_{\tilde{\pi}}(\tilde{t}) c_{\tilde{\sigma}}(\tilde{t}) D^{\tilde{H}}(\tilde{t}) \Delta(\tilde{t})^{r+s} d\tilde{t}
\end{aligned}\]

\noindent où $D^{\tilde{H}}$ et $\Delta$ sont certaines fonctions déterminants qui seront définies dans le corps du texte. Cette expression a un sens (cf 3.4). Posons enfin

$$\epsilon_\nu(\pi,\sigma)=\omega_{\pi}((-1)^{[m/2]}2\nu)\omega_{\sigma}((-1)^{1+[d/2]}2\nu) \epsilon(1/2,\pi\times \sigma,\psi_E)$$

\noindent où $\omega_\pi$ et $\omega_\sigma$ sont les caractères centraux de $\pi$ et $\sigma$ et où $[x]$ désigne la partie entière de $x$. Le but de cet article est alors d'établir le résultat suivant (théorème 6.1.1)

\begin{theom}
On a l'égalité

$$\epsilon_\nu(\pi,\sigma)=\epsilon_{geom,\nu}(\tilde{\pi},\tilde{\sigma})$$
\end{theom}

Essayons maintenant d'expliquer, brièvement, comment on démontre le théorème $1$. On suppose pour simplifier que $m=d-1$. Notons $Hom_H(\pi,\sigma)$ l'espace des entrelacements $H(F)$-équivariant entre $\pi$ et $\sigma$. Cet espace est de dimension $1$. Soit $\ell\in Hom_H(\pi,\sigma)$ non nulle. Il existe alors une constante $c$ vérifiant $\tilde{\sigma}(\tilde{y})\circ \ell\circ \tilde{\pi}(\tilde{y})^{-1}=c\ell$ pour tout $\tilde{y}\in \tilde{H}(F)$. Cette constante $c$ mesure en quelque sorte le défaut de $\ell$ d'être un homomorphisme $\tilde{H}(F)$-équivariant. On montre en 4.1 que $c$ est égal à des facteurs élémentaires près au facteur epsilon $\epsilon(1/2,\pi\times \sigma,\psi_E)$. On peut définir de façon naturelle l'induite lisse de $\tilde{\sigma}$ à $\tilde{G}(F)$, notée $Ind_{\tilde{H}}^{\tilde{G}}(\tilde{\sigma})$. C'est une représentation lisse de $\tilde{G}(F)$ et on a un isomorphisme de Frobenius

$$Hom_{\tilde{H}}(\tilde{\pi},\tilde{\sigma})=Hom_{\tilde{G}}(\tilde{\pi},Ind_{\tilde{H}}^{\tilde{G}}(\tilde{\sigma}))$$

\noindent L'induite $Ind_{\tilde{H}}^{\tilde{G}}(\tilde{\sigma})$ se réalise naturellement sur un espace de fonctions de $G(F)$ dans $E_\sigma$ vérifiant une condition d'équivariance à gauche par $H(F)$. Soit $\tilde{f}\in C_c^\infty(\tilde{G}(F))$ une fonction localement constante à support compact sur $\tilde{G}(F)$. Cette fonction agit sur $Ind_{\tilde{H}}^{\tilde{G}}(\tilde{\sigma})$ comme un opérateur à noyau de noyau

$$K_{\tilde{f}}(x,y)=\int_{\tilde{H}(F)} \tilde{f}(y^{-1}\tilde{h}x)\tilde{\sigma}(\tilde{h}) d\tilde{h}$$

\noindent Formellement on pourrait alors écrire la trace de $Ind_{\tilde{H}}^{\tilde{G}}\tilde{\sigma}(\tilde{f})$ comme l'intégrale sur $H(F)\backslash G(F)$ de $Tr\big(K_{\tilde{f}}(x,x)\big)$. Malheureusement, en général $Ind_{\tilde{H}}^{\tilde{G}}\tilde{\sigma}(\tilde{f})$ n'est pas un opérateur à trace et l'intégrale précédente diverge. On utilise alors la même idée qu'Arthur pour la formule des traces locale: on va tronquer l'intégrale précédente. On définit en 3.1 une suite croissance exhaustive $(\Omega_N)_{N\geqslant 1}$ de sous-ensembles compacts de $H(F)\backslash G(F)$ et pour tout $N\geqslant 1$, on note $\kappa_N$ la fonction caractéristique de $\Omega_N$. On pose alors

$$\displaystyle J_N(\Theta_{\tilde{\sigma}},\tilde{f})=\int_{H(F)\backslash G(F)} \kappa_N(g) Tr\big(K_{\tilde{f}}(g,g)\big)dg$$

\noindent On va évaluer la limite de cette expression lorsque $N$ tend vers l'infini de deux façons: l'une qualifiée de géométrique, l'autre de spectrale. L'expression spectrale contient les constantes $c$ ci-dessus (pour différentes représentations $\pi$), donc les facteurs epsilon correspondants. De l'égalité entre les deux développements, on déduit le théorème $1$. Il faut pour cela , suivant un procédé dû à Arthur, rendre les distributions obtenus (d'un côté comme de l'autre) invariantes puis appliquer l'égalité entre les deux à un pseudocoefficient de $\tilde{\pi}$. Néanmoins, afin que l'expression $J_N(\Theta_{\tilde{\sigma}},\tilde{f})$ admette une limite, il nous faut imposer à $\tilde{f}$ d'être très cuspidale. On renvoie à 1.9 pour la définition de très cuspidale. Indiquons enfin que les preuves et méthodes utilisées ici suivent grandement [W3], à tel point que certaines démonstrations, très proches de celles présentées dans [W3], n'ont pas été réécrites (notamment lors du développement spectral). \\

Décrivons brièvement le plan de l'article. Dans la première section, on rassemble les définitions, notations et résultats généraux sur les groupes tordus utiles par la suite. La deuxième section est elle consacrée à l'étude des groupes tordus qui nous intéressent: les groupes tordus de changement de base des groupes unitaires. La troisième section contient la définition de $J_N(\Theta_{\tilde{\sigma}},\tilde{f})$ et établit le développement géométrique. La quatrième section prouve le lien entre la constante $c$ et le facteur epsilon éludé ci-dessus. L'expression spectrale de la limite est énoncée et démontrée dans la cinquième section. Enfin, l'application au calcul d'un facteur epsilon de paire (théorème 1) se trouve dans la sixième et dernière section. \\

\textbf{Remerciement}: Je souhaite remercier ici Jean-Loup Waldspurger pour m'avoir proposé ce sujet ainsi que pour sa disponibilité constante et ses relectures rigoureuses du présent article.

\section{Groupes tordus}

\subsection{Notations générales}

Soit $F$ un corps local non archimédien de caractéristique $0$. On notera respectivement $\mathcal{O}_F$, $\mathfrak{p}_F$, $\pi_F$, $k_F$, $val_F$, $|.|_F$ l'anneau des entiers de $F$, son idéal maximal, une uniformisante, le corps résiduel, la valuation normalisée et la valeur absolue (normalisée par $|\pi_F|_F=|k_F|^{-1}$). Fixons un caractère additif non trivial $\psi:F\to\mathbb{C}^\times$. Soit $G$ un groupe réductif connexe défini sur $F$. Sauf mention explicite du contraire, toutes les variétés, groupes, sous-variétés, sous-groupes, morphismes algébriques seront supposés définis sur $F$. On note $A_G$ la composante déployée du centre connexe de $G$, $X^*(G)$ le groupe des caractères algébriques (définis sur $F$) de $G$, $\mathcal{A}_G^*=X^*(G)\otimes_{\mathbb{Z}}\mathbb{R}$ et $\mathcal{A}_G=Hom(X^*(G),\mathbb{R})$. On note $\mathfrak{g}$ l'algèbre de Lie de $G$ et

$$G\times\mathfrak{g}\to\mathfrak{g}$$
$$(g,X)\mapsto gXg^{-1}$$

\noindent l'action adjointe.

\subsection{Définitions des groupes tordus}

 Un groupe tordu est un couple $(G,\tilde{G})$ où $G$ est un groupe réductif connexe défini sur $F$ et $\tilde{G}$ est une variété algébrique définie sur $F$ vérifiant $\tilde{G}(F)\neq \emptyset$ qui est munie de deux actions commutantes à droite et à gauche par $G$ qui font chacune de $\tilde{G}$ un espace principal homogène sous $G$. On notera $(g,\tilde{x})\in G\times \tilde{G}\mapsto g\tilde{x}$ resp. $(\tilde{x},g)\in\tilde{G}\times G\mapsto \tilde{x}g$ l'action à gauche resp. à droite de $G$ sur $\tilde{G}$. La plupart du temps on omettra le groupe $G$ et on parlera du groupe tordu $\tilde{G}$, le groupe $G$ sous-jacent étant sous-entendu. \\
 
 Soit $\tilde{G}$ un groupe tordu. Pour $\tilde{x}\in \tilde{G}$, il existe un unique automorphisme $\theta_{\tilde{x}}$ de $G$ tel que $\tilde{x}g=\theta_{\tilde{x}}(g)\tilde{x}$ pour tout $g\in G$. On déduit de $\theta_{\tilde{x}}$ un automorphisme de $X^*(G)$, de $A_G$, de $\mathcal{A}_G$ etc... qui ne dépend pas de $\tilde{x}$. On note $\theta_{\tilde{G}}$ cet automorphisme. On fait l'hypothèse suivante \\
 
 \textbf{Hypothèse}: $\theta_{\tilde{G}}$ est d'ordre fini. \\
 
 Cette hypothèse nous permet d'appliquer les résultats présents dans la littérature pour les groupes réductifs non connexes. On adopte les notations suivantes: \\
 
 $A_{\tilde{G}}=(A_G^{\theta_{\tilde{G}}=1})^0$, $\mathcal{A}_{\tilde{G}}=\mathcal{A}_G^{\theta_{\tilde{G}}=1}$, $\mathcal{A}_{\tilde{G}}^*=\left(\mathcal{A}_G^*\right)^{\theta_{\tilde{G}}=1}$, $a_{\tilde{G}}=dim(\mathcal{A}_{\tilde{G}})$. \\
 
\noindent où l'exposant $0$ indique que l'on prend la composante neutre. On définit l'homomorphisme $H_{\tilde{G}}: G(F)\to \mathcal{A}_{\tilde{G}}$ par $H_{\tilde{G}}(g)(\chi)=log \; |\chi(g)|_F$ pour tout $\chi\in X^*(G)^{\theta_{\tilde{G}}=1}$. \\
 
 Le groupe $G$ agit par conjugaison sur $\tilde{G}$: $(g,\tilde{x})\mapsto g\tilde{x}g^{-1}$. Pour $\tilde{X}$ un sous-ensemble de $\tilde{G}$, on notera $Norm_G(\tilde{X})$ resp. $Z_G(\tilde{X})$ resp. $G_{\tilde{X}}$ le normalisateur resp. le centralisateur resp. le centralisateur connexe de $\tilde{X}$ dans $G$. Si $\tilde{X}=\{\tilde{x}\}$, on adoptera plutôt les notations $Norm_G(\tilde{x})$, $Z_G(\tilde{x})$ et $G_{\tilde{x}}$ et on désignera par $\mathfrak{g}_{\tilde{x}}$ l'algèbre de Lie de $G_{\tilde{x}}$. Pour $X$ un sous-ensemble de $G$, on notera $N_{\tilde{G}}(X)$ resp. $Z_{\tilde{G}}(X)$ le normalisateur resp. le centralisateur de $X$ dans $\tilde{G}$ pour l'opération $(\tilde{x},g)\mapsto \theta_{\tilde{x}}(g)$. \\
 
 Un élément $\tilde{x}\in \tilde{G}$ sera dit semi-simple s'il existe un couple $(B,T)$ constitué d'un sous-groupe de Borel de $G$ et d'un tore maximal de $B$ tous deux définis sur $\overline{F}$ tel que $\tilde{x}$ normalise $B$ et $T$. On notera $\tilde{G}_{ss}$ l'ensemble des éléments semi-simples dans $\tilde{G}$. Pour $\tilde{x}\in \tilde{G}_{ss}$, on pose
 
$$D^{\tilde{G}}(\tilde{x})=\left|det \; (1-\theta_{\tilde{x}})_{|\mathfrak{g}/\mathfrak{g}_{\tilde{x}}} \right|_F$$

Un élément $\tilde{x}\in \tilde{G}$ sera qualifié de régulier si $Z_G(\tilde{x})$ est abélien et $G_{\tilde{x}}$ est un tore. On note $\tilde{G}_{reg}$ l'ensemble des éléments réguliers. \\

On appelle sous-groupe parabolique tordu de $\tilde{G}$ tout couple $(P,\tilde{P})$ où $P$ est un sous-groupe parabolique de $G$ (défini sur $F$) et $\tilde{P}$ est le normalisateur de $P$ dans $\tilde{G}$ avec la condition $\tilde{P}(F)\neq \emptyset$. Pour un tel couple, $\tilde{P}$ détermine entièrement $P$, on parlera donc plutôt du sous-groupe parabolique tordu $\tilde{P}$. Soit $\tilde{P}$ un sous-groupe parabolique tordu, on étend de façon naturelle le module $\delta_P$ à $\tilde{P}(F)$. Une composante de Levi tordue de $\tilde{P}$ est un couple $(M,\tilde{M})$ constitué d'une composante de Levi $M$ de $P$ (définie sur $F$) et du normalisateur $\tilde{M}$ de $M$ dans $\tilde{P}$. On a alors $\tilde{M}(F)\neq \emptyset$, donc le couple $(M,\tilde{M})$ définit un groupe tordu. On appelle Levi tordu de $\tilde{G}$ toute composante de Levi tordue d'un sous-groupe parabolique tordu de $\tilde{G}$. De la même façon, pour un Levi tordu $(M,\tilde{M})$, le deuxième terme détermine entièrement le premier, on parlera donc du Levi tordu $\tilde{M}$. Pour $\tilde{M}$ un Levi tordu, on reprend les notations d'Arthur: $\mathcal{P}(\tilde{M})$ resp. $\mathcal{F}(\tilde{M})$ resp. $\mathcal{L}(\tilde{M})$ désignera l'ensemble des sous-groupes paraboliques tordus de composante de Levi tordue $\tilde{M}$ resp. des sous-groupes paraboliques tordus contenant $\tilde{M}$ resp. des Levi tordus contenant $\tilde{M}$. Pour $\tilde{M}$, $\tilde{L}$ des Levi tordus et $\tilde{P}$ un sous-groupe parabolique tordu, on remarque que $\tilde{M}\subset \tilde{L}$ et $\tilde{M}\subset \tilde{P}$ entraînent respectivement $M\subset L$ et $M\subset P$. Ainsi il ne peut y avoir d'ambiguïté dans les définitions précédentes.
 Soit $\tilde{Q}$ un sous-groupe parabolique tordu. lorsque l'on notera $\tilde{Q}=\tilde{L}U$ cela signifiera que $\tilde{L}$ est une composante de Levi tordue de $\tilde{Q}$ et que $U$ est le radical unipotent de $Q$. Les Levi tordus se caractérisent comme les commutants dans $\tilde{G}$ de tores déployés. On a en effet:

\begin{itemize}
\item  Si $A$ est un sous-tore déployé de $G$ tel que $Z_{\tilde{G}}(A)(F)\neq \emptyset$ alors $Z_{\tilde{G}}(A)$ est un Levi tordu de $\tilde{G}$.
\item Réciproquement si $\tilde{M}$ est un Levi tordu de $\tilde{G}$, alors $\tilde{M}=Z_{\tilde{G}}(A_{\tilde{M}})$.
\end{itemize}

 Fixons un sous-groupe prabolique minimal $P_{min}$ de $G$ ainsi qu'une composante de Levi $M_{min}$ de celui-ci. On pose $W^G=Norm_{G(F)}(M_{min})/M_{min}(F)$. On appelera sous-groupe parabolique standard (resp. semistandard) tout sous-groupe parabolique $P$ de $G$ qui contient $P_{min}$ (resp. qui contient $M_{min}$). Un sous-groupe parabolique semistandard admet une unique composante de Levi contenant $M_{min}$ et lorsque l'on écrira $P=MU$ avec $P$ semistandard, on sous-entendra que $M$ est cette composante de Levi (et que $U$ est le radical unipotent de $P$). Posons $\tilde{P}_{min}=N_{\tilde{G}}(P_{min})$ et $\tilde{M}_{min}=N_{\tilde{G}}(P_{min},M_{min})$. Alors $\tilde{P}_{min}$ est un sous-groupe parabolique tordu de $\tilde{G}$ et $\tilde{M}_{min}$ est une composante de Levi tordue de celui-ci. De plus, $\tilde{P}_{min}$ est un sous-groupe parabolique tordu minimal et $\tilde{M}_{min}$ un Levi minimal en un sens évident. On adoptera la notation $\mathcal{L}^{\tilde{G}}=\mathcal{L}(\tilde{M}_{min})$. \\
 
 Fixons un sous-groupe compact spécial $K$ de $G$ en bonne position relativement à $M_{min}$. Soient $\tilde{M}\in \mathcal{L}^{\tilde{G}}$ et $\tilde{P}=\tilde{M}U\in\mathcal{P}(\tilde{M})$. On a alors $G(F)=M(F)U(F)K$. On définit une application $H_{\tilde{P}}:G(F)\to\mathcal{A}_{\tilde{M}}$ par $H_{\tilde{P}}(g)=H_{\tilde{M}}(m)$ où $g=muk$ est une décomposition de $g$ avec $m\in M(F)$, $u\in U(F)$ et $k\in K$. \\
 
 Un sous-tore maximal tordu de $\tilde{G}$ est un couple $(T,\tilde{T})$ constitué d'un sous-tore maximal $T$ de $G$ (défini sur $F$) et d'une sous-variété $\tilde{T}$ de $\tilde{G}$ (définie sur $F$) qui est l'intersection des normalisateurs de $T$ et d'un sous-groupe de Borel $B$ défini sur $\overline{F}$ contenant $T$, vérifiant $\tilde{T}(F)\neq \emptyset$. Pour un tel couple, la restriction à $T$ de l'automorphisme $\theta_{\tilde{x}}$ pour $\tilde{x}\in \tilde{T}$ ne dépend pas de $\tilde{x}$ et on notera $\theta_{\tilde{T}}$, ou simplement $\theta$ si aucune confusion n'est possible, cet automorphisme de $T$. On désignera par $T_\theta$ la composante neutre du sous-groupe des points fixes $T^\theta$. \\
 
 Fixons un produit scalaire sur $\mathcal{A}_{M_{min}}$ invariant par l'action de $W^G$. Pour tout Levi tordu $\tilde{M}$, on en déduit par conjugaison et restriction un produit scalaire sur $\mathcal{A}_{\tilde{M}}$. Fixons une extension de $H_{\tilde{G}}$ à $\tilde{G}(F)$, c'est-à-dire une application $H_{\tilde{G}}: \tilde{G}(F)\to \mathcal{A}_{\tilde{G}}$ vérifiant $H_{\tilde{G}}(g\tilde{x}g')=H_{\tilde{G}}(g)+H_{\tilde{G}}(\tilde{x})+H_{\tilde{G}}(g')$ pour tous $g,g'\in G(F)$, $\tilde{x}\in\tilde{G}(F)$. Pour tout Levi tordu $\tilde{M}$ de $\tilde{G}$, il existe une unique extension de $H_{\tilde{M}}$ à $\tilde{M}(F)$ vérifiant
 
\begin{itemize}
\item la composée de $H_{\tilde{M}}$ avec la projection orthogonale sur $\mathcal{A}_{\tilde{G}}$ coïncide avec la restriction de $H_{\tilde{G}}$ à $\tilde{M}(F)$;
\item pour tout $n\in Norm_{G(F)}(\tilde{M})$ et pour tout $\tilde{m}\in\tilde{M}(F)$ on a $H_{\tilde{M}}(n\tilde{m}n^{-1})=H_{\tilde{M}}(\tilde{m})$.
\end{itemize}

\subsection{Mesures}

On suppose fixé jusqu'en 1.9 un groupe tordu $(G,\tilde{G})$, un Levi tordu minimal $(M_{min},\tilde{M}_{min})$ de $\tilde{G}$ et un sous-groupe compact spécial $K$ de $G(F)$ en bonne position par rapport à $M_{min}$. On fixe une mesure de Haar sur $G(F)$ et on munit $K$ de la mesure de Haar de masse totale $1$. Soit $P=MU\in\mathcal{F}(M_{min})$. Il existe sur $U(F)$ une unique mesure de Haar telle que

$$\displaystyle\int_{U(F)} \delta_{\overline{P}}(m_{\overline{P}}(u)) du=1$$

\noindent où $\overline{P}$ est le sous-groupe parabolique opposé à $P$. On munit $U(F)$ de cette mesure de Haar et $M(F)$ de l'unique mesure de Haar de sorte que l'on ait l'égalité

$$\displaystyle \int_{G(F)} f(g) dg=\int_K\int_{M(F)}\int_{U(F)} f(muk) dudmdk$$

\noindent pour tout $f\in C_c^\infty(G(F))$. Soit $T$ un tore. Si $T$ est déployé, on choisit sur $T(F)$ la mesure de Haar qui donne une mesure de $1$ au sous-groupe compact maximal. Dans le cas général, on munit $A_T(F)$ de cette mesure et $T(F)$ de la mesure telle que $T(F)/A_T(F)$ soit de masse totale $1$. \\

Les mesures que l'on vient de définir sur les groupes sont celles qui seront utilisées tout au long de l'article à l'exception des sections 3.6 à 3.9. Lors du développement géométrique (sections 3.6 à 3.9), on préférera une autre normalisation qui est la suivante. On fixe une forme bilinéaire non dégénérée $<.,.>:\mathfrak{g}(F)\times\mathfrak{g}(F)\to F$ qui est invariante par l'action adjointe de $G(F)$. Soit $\mathfrak{h}$ une sous-algèbre de $\mathfrak{g}$,

\begin{itemize}
\item Si la restriction de $<.,.>$ à $\mathfrak{h}(F)$ est non dégénérée, alors on munit $\mathfrak{h}(F)$ de la mesure autoduale pour le bicaractère $(X,Y)\mapsto \psi(<X,Y>)$;
\item Sinon, on munit $\mathfrak{h}(F)$ d'une mesure de Haar quelconque.
\end{itemize}

\noindent Si $H$ est un sous-groupe de $G$, on relève alors la mesure fixée sur $\mathfrak{h}(F)$ en une mesure de Haar sur $H(F)$ via l'exponentielle. Si $T$ est un tore, notons $d_c t$ la mesure de Haar fixée sur $T(F)$ dans notre premier choix de normalisation, alors il existe un réel $\nu(T)$ de sorte que $d_c t=\nu(T)dt$. \\

Supposons maintenant fixées nos mesures sur les groupes dans l'une ou l'autre des normalisations précédentes. Si $(G',\tilde{G}')$ est un groupe tordu avec $G'$ un sous-groupe de $G$, la mesure sur $G'(F)$ en détermine une sur l'espace principal homogène $\tilde{G}'(F)$. Soit $(T,\tilde{T})$ un sous-tore maximal tordu de $(G,\tilde{G})$. Notons $\tilde{T}(F)/\theta$ l'ensemble des orbites pour l'action de $T(F)$ par conjugaison sur $\tilde{T}(F)$. C'est naturellement une variété $F$-analytique et pour tout $\tilde{t}\in \tilde{T}(F)/\theta$, l'application

$$T_\theta(F)\to\tilde{T}(F)/\theta$$
$$t\mapsto t\tilde{t}$$

\noindent est un isomorphisme local. Il existe une unique mesure sur $\tilde{T}(F)/\theta$ telle que l'application précédente préserve localement les mesures pour tout $\tilde{t}\in \tilde{T}(F)/\theta$. On munit $\tilde{T}(F)/\theta$ de cette mesure. On a alors la formule d'intégration de Weyl

$$\displaystyle \int_{\tilde{G}(F)} \tilde{f}(\tilde{x}) d\tilde{x}=\sum_{\tilde{T}\in \mathcal{T}(\tilde{G})} |W(G,\tilde{T})|^{-1} [T^\theta(F):T_\theta(F)]^{-1} \int_{\tilde{T}(F)/\theta} D^{\tilde{G}}(\tilde{t})\int_{T_\theta(F)\backslash G(F)} f(g^{-1}\tilde{t}g) dgd\tilde{t}$$

\noindent pour tout $\tilde{f}\in C_c^\infty(\tilde{G}(F))$. \\

Posons $\mathcal{A}_{\tilde{G},F}=H_{\tilde{G}}(G(F))$, $\mathcal{A}_{A_{\tilde{G}},F}=H_{\tilde{G}}(A_{\tilde{G}}(F))$, $\mathcal{A}_{\tilde{G},F}^\vee=Hom(\mathcal{A}_{\tilde{G},F},2\pi\mathbb{Z})$, $\mathcal{A}_{A_{\tilde{G}},F}^\vee=Hom(\mathcal{A}_{A_{\tilde{G}},F},2\pi\mathbb{Z})$. Alors $\mathcal{A}_{\tilde{G},F}$, $\mathcal{A}_{A_{\tilde{G}},F}$ sont des réseaux dans $\mathcal{A}_{\tilde{G}}$ et $\mathcal{A}_{\tilde{G},F}^\vee$, $\mathcal{A}_{A_{\tilde{G}},F}^\vee$ sont des réseaux dans $\mathcal{A}_{\tilde{G}}^*$. On munit ces réseaux des mesures de comptage. On choisit sur $\mathcal{A}_{\tilde{G}}$ (resp. $\mathcal{A}_{\tilde{G}}^*$) la mesure de Haar telle que $\mathcal{A}_{\tilde{G}}/\mathcal{A}_{A_{\tilde{G}},F}$ (resp. $\mathcal{A}_{\tilde{G}}^*/\mathcal{A}_{A_{\tilde{G}},F}^\vee$) soit de masse totale $1$.

\subsection{$(\tilde{G},\tilde{M})$-familles}

Soit $\tilde{M}$ un Levi tordu de $\tilde{G}$. On peut étendre au cas tordu les notions de $(\tilde{G},\tilde{M})$-familles et de familles de points $(\tilde{G},\tilde{M})$-orthogonaux cf 2.3 [W4]. A une $(\tilde{G},\tilde{M})$-famille $(\varphi_{\tilde{P}})_{\tilde{P}\in\mathcal{P}(\tilde{M})}$ est associé un nombre (ou un opérateur si la $(\tilde{G},\tilde{M})$-famille est à valeurs opérateurs) que l'on note en général $\varphi_{\tilde{M}}(0)$. Cela dépend du choix d'une mesure sur $\mathcal{A}_{\tilde{M}}^{\tilde{G}}=\mathcal{A}_{\tilde{M}}/\mathcal{A}_{\tilde{G}}$, choix qui a été fixé au paragraphe précédent.

\subsection{Représentations de groupes tordus}

On appelle représentation de $\tilde{G}(F)$ tout triplet $(\pi,\tilde{\pi},E_\pi)$ où $\pi$ est une représentation lisse de $G(F)$ d'espace $E_\pi$ et $\tilde{\pi}$ est une application $\tilde{\pi}: \tilde{G}(F)\to Aut_{\mathbb{C}}(E_\pi)$ vérifiant $\tilde{\pi}(g\tilde{x}g')=\pi(g)\tilde{\pi}(\tilde{x})\pi(g')$ pour tous $g,g'\in G(F)$, $\tilde{x}\in \tilde{G}(F)$. Deux représentations $(\pi_1,\tilde{\pi}_1,E_{\pi_1})$ et $(\pi_2,\tilde{\pi}_2,E_{\pi_2})$ sont dites équivalentes s'il existe des isomorphismes linéaires $A:E_{\pi_1}\to E_{\pi_2}$ et $B:E_{\pi_1}\to E_{\pi_2}$ qui entrelacent $\pi_1$ et $\pi_2$ et qui vérifient $B\tilde{\pi}_1(\tilde{x})=\tilde{\pi}_2(\tilde{x})A$ pour tout $\tilde{x}\in \tilde{G}(F)$. On dit d'une représentation $(\pi,\tilde{\pi},E_\pi)$ de $\tilde{G}(F)$ qu'elle est admissible si $\pi$ l'est, qu'elle est unitaire s'il existe un produit hermitien défini positif sur $E_\pi$ invariant par l'image de $\tilde{\pi}$ et qu'elle est tempérée si elle est unitaire, $\pi$ est de longueur finie et toutes les sous-représentations irréductibles de $\pi$ sont tempérées. On définit la représentation duale de $(\pi,\tilde{\pi}, E_\pi)$: c'est la représentation $(\pi^\vee,\tilde{\pi}^\vee,E_{\pi^\vee})$ où $\pi^\vee$ est la représentation duale de $\pi$ et $\tilde{\pi}^\vee$ est telle que $<\tilde{\pi}^\vee(\tilde{x})v^\vee,\tilde{\pi}(\tilde{x})v>=<v^\vee,v>$ pour tous $v\in E_\pi$, $v^\vee\in E_{\pi^\vee}$ et $\tilde{x}\in \tilde{G}(F)$. Une représentation $(\pi,\tilde{\pi},E_\pi)$ est dite $G(F)$-irréductible si $\pi$ est irréductible. Pour $\lambda\in \mathcal{A}_{\tilde{G}}^*\otimes_{\mathbb{R}} \mathbb{C}$, on définit $(\pi_\lambda,\tilde{\pi}_\lambda,E_{\pi_\lambda})$ par $E_{\pi_\lambda}=E_\pi$ et $\tilde{\pi}_\lambda(\tilde{x})=e^{\lambda(H_{\tilde{G}}(\tilde{x}))} \tilde{\pi}(\tilde{x})$ pour $\tilde{x}\in\tilde{G}(F)$. On omettra en général les termes $\pi$ et $E_\pi$ et on parlera de la représentation $\tilde{\pi}$ de $\tilde{G}(F)$, on notera alors sans plus de commentaire $(\pi,E_\pi)$ la représentation de $G(F)$ sous-jacente. \\

Soient $\tilde{P}=\tilde{M}U$ un sous-groupe parabolique tordu de $\tilde{G}$ et $\tilde{\tau}$ une représentation de $\tilde{M}(F)$. On définit l'induite normalisée notée $i^{G}_{P}(\tilde{\tau})$ comme suit. La représentation de $G(F)$ sous-jacente est l'induite normalisée $i^G_P(\tau)$ qui se réalise sur l'espace $E^G_{P,\tau}$ des fonctions $e: G(F)\to E_{\tau}$ vérifiant

\begin{itemize}
\item $e$ est invariante à droite par un sous-groupe compact-ouvert de $G(F)$;
\item $e(mug)=\delta_P(m)^{1/2}\tau(m)e(g)$ pour tous $m\in M(F),u\in U(F), g\in G(F)$,
\end{itemize}
\noindent L'action de $\tilde{G}(F)$ sur $E^G_{P,\tau}$ est définie comme suit. Soit $\tilde{x}\in \tilde{G}(F)$ et fixons $\tilde{m}\in\tilde{M}$. Notons $\gamma\in G(F)$ l'unique élément tel que $\tilde{x}=\gamma\tilde{m}$. On a alors
 $$\left(i^{G}_{P}(\tilde{\tau},\tilde{x})e\right)(g)=\delta_P(\tilde{m})^{1/2}\tilde{\tau}(\tilde{m})e(\theta_{\tilde{m}}^{-1}(g\gamma))$$
 
\noindent pour tous $e\in E^G_{P,\tau}$, $g\in G(F)$. Supposons que $M$ est semistandard. On peut alors réaliser la représentation induite sur l'espace $\mathcal{K}_{P,\tau}^G$ des restrictions à $K$ des éléments de $E^G_{P,\tau}$. Supposons de plus $\tilde{\tau}$ tempérée. On définit alors une $(\tilde{G},\tilde{M})$-famille qui prend ses valeurs dans l'espace des opérateurs sur $\mathcal{K}^G_{P,\tau}$, par

$$\mathcal{R}_{P'}(\tau,\lambda)= R_{P'|P}(\tau)^{-1}R_{P'|P}(\tau_\lambda)$$

\noindent pour tout $\tilde{P}'\in\mathcal{P}(\tilde{M})$ et pour tout $\lambda\in i\mathcal{A}_{\tilde{M}}^*$, où $R_{P'|P}(\tau)$ est l'opérateur d'entrelacement normalisé comme en [A2] théorème 2.1. On déduit de cette $(\tilde{G},\tilde{M})$-famille un opérateur $\mathcal{R}_{\tilde{M}}(\tau)$ qui ne dépend pas de $\tilde{P}$. On définit alors le caractère pondéré de $\tilde{\tau}$ par

$$\displaystyle J^{\tilde{G}}_{\tilde{M}}(\tilde{\tau},\tilde{f})=Tr(\mathcal{R}_{\tilde{M}}(\tau)i^G_P(\tilde{\tau},\tilde{f}))$$

\noindent pour tout $\tilde{f}\in C_c^\infty(\tilde{G}(F))$. Dans le cas où $\tilde{M}=\tilde{G}$, c'est le caractère classique de $\tilde{\tau}$ et on le note plutôt $\tilde{f}\mapsto \Theta_{\tilde{\tau}}(\tilde{f})$. D'après [C] théorème 2, ce caractère est une distribution localement intégrable.

\subsection{Représentations tempérées et elliptiques des groupes tordus}

Soit $M$ un Levi de $G$ et $\tau\in\Pi_2(M)$. Supposons la condition suivante vérifiée

\begin{center}
(1) $\pi=i^G_P(\tau)$ ($P\in\mathcal{P}(M)$) est irréductible
\end{center}

Posons

$$W^{\tilde{G}}(\tau)=\{\tilde{x}\in Norm_{\tilde{G}(F)}(M); \; \tau\circ \theta_{\tilde{x}}\simeq \tau\}/ M(F)$$

Alors la représentation $\pi$ s'étend en une représentation $G(F)$-irréductible et tempérée $\tilde{\pi}$ de $\tilde{G}(F)$ si et seulement si $W^{\tilde{G}}(\tau)\neq \emptyset$. De plus, la représentation $\tilde{\pi}$ est alors unique à équivalence près. Supposons que l'hypothèse (1) est vérifiée pour tout Levi $M$ et pour tout $\tau\in\Pi_2(M)$. Alors en laissant varier $M$ et $\tau$, on obtient ainsi toutes les représentations tempérées $G(F)$-irréductibles de $\tilde{G}(F)$. Soient $L$ un autre Levi et $\tau'\in \Pi_2(L)$ tels que $W^{\tilde{G}}(\tau')\neq\emptyset$ et $\tilde{\pi'}$ une représentation $G(F)$-irréductible et tempérée de $\tilde{G}(F)$ qui étend $\pi'=i^G_Q(\tau')$ ($Q\in\mathcal{P}(L)$). Alors $\tilde{\pi}$ et $\tilde{\pi'}$ sont équivalentes si et seulement si il existe $g\in G(F)$ tel que $gLg^{-1}=M$ et $\tau\circ ad_g\simeq \tau'$. Revenons à notre représentation $\tilde{\pi}$ obtenue à partir de $M$ et $\tau$. Notons $W^{\tilde{G}}(\tau)_{reg}$ l'ensemble des $\tilde{w}\in W^{\tilde{G}}(\tau)$ tels que $\mathcal{A}_M^{\tilde{w}}=\mathcal{A}_{\tilde{G}}$. On dit que $\tilde{\pi}$ est elliptique si $W^{\tilde{G}}(\tau)=W^{\tilde{G}}(\tau)_{reg}$. Cela implique notamment que $W^{\tilde{G}}(\tau)$ est réduit à un élément que l'on note $\tilde{t}$. On pose alors $\iota(\tilde{\pi})=|det (1-\tilde{t})_{|\mathcal{A}^{\tilde{G}}_M}|^{-1}$. On définit $\Pi_{ell}(\tilde{G})$ comme l'ensemble des représentations elliptiques de $\tilde{G}(F)$ à équivalence près. On a une action naturelle de $i\mathcal{A}_{\tilde{G},F}^*$ sur $\Pi_{ell}(\tilde{G})$ et on note $\{\Pi_{ell}(\tilde{G})\}$ l'ensemble des orbites pour cette action. Soient $\mathcal{O}\in \{\Pi_{ell}(\tilde{G})\}$ et $\tilde{\pi}\in \mathcal{O}$. Alors $\iota(\tilde{\pi})$ et le stabilisateur de $\tilde{\pi}$ dans $i\mathcal{A}_{\tilde{G}}^*$ ne dépendent pas du choix de $\tilde{\pi}$, on les note respectivement $\iota(\mathcal{O})$ et $i\mathcal{A}_{\mathcal{O}}^\vee$.

\subsection{Intégrales orbitales pondérées non invariantes et invariantes}

Soit $\tilde{M}\in\mathcal{L}^{\tilde{G}}$. Soit $g\in G(F)$, la famille $\left(H_{\tilde{P}}(g)\right)_{\tilde{P}\in\mathcal{P}(\tilde{M})}$ est $(\tilde{G},\tilde{M})$-orthogonale et on peut donc lui associer une $(\tilde{G},\tilde{M})$-famille $(v_{\tilde{P}}(g))_{\tilde{P}\in\mathcal{P}(\tilde{M})}$ par la formule $v_{\tilde{P}}(g,\lambda)=e^{-\lambda(H_{\tilde{P}}(g))}$ pour $\lambda\in i\mathcal{A}_{\tilde{M}}^*$. De cette $(\tilde{G},\tilde{M})$-famille se déduit un nombre $v_{\tilde{M}}(g)$. Soient $\tilde{f}\in C_c^\infty(\tilde{G}(F))$ et $\tilde{x}\in \tilde{M}(F)\cap \tilde{G}_{reg}(F)$, on pose

$$\displaystyle J_{\tilde{M}}(\tilde{x},\tilde{f})=D^{\tilde{G}}(\tilde{x})^{1/2}\int_{G_{\tilde{x}}(F)\backslash G(F)} \tilde{f}(g^{-1}\tilde{x}g) v_{\tilde{M}}(g) dg$$

\noindent C'est l'intégrale orbitale pondérée de $\tilde{f}$ en $\tilde{x}$. \\

On peut aussi définir, grâce aux caractères pondérés et suivant un procédé d'Arthur, des intégrales orbitales pondérées invariantes. Définissons $\mathcal{H}_{ac}(\tilde{G}(F))$ comme l'espace des fonctions $\tilde{f}: \tilde{G}(F)\to\mathbb{C}$ vérifiant

\begin{itemize}
\item $\tilde{f}$ est biinvariante par un sous-groupe compact-ouvert de $G(F)$;
\item pour tout $X\in \mathcal{A}_{\tilde{G}}$ la fonction $\tilde{f}\mathbf{1}_{H_{\tilde{G}}=X}$ est à support compact.
\end{itemize}

Deux fonctions $\tilde{f},\tilde{f}'\in C_c^\infty(\tilde{G}(F))$ sont dites équivalentes si $D(\tilde{f})=D(\tilde{f}')$ pour toute distribution invariante $D$. On a alors

\begin{prop}
Soient $\tilde{f}\in C_c^\infty(\tilde{G}(F))$ et $\tilde{M}\in \mathcal{L}^{\tilde{G}}$. Il existe $\phi_{\tilde{M}}(\tilde{f})\in \mathcal{H}_{ac}(\tilde{M}(F))$ telle que pour tout $\tilde{\pi}\in Temp(\tilde{M})$ et pour tout $X\in H_{\tilde{M}}(\tilde{M}(F))$ on ait l'égalité

$$\displaystyle \int_{i\mathcal{A}_{\tilde{M},F}^*} J_{\tilde{M}}(\tilde{\pi}_\lambda,\tilde{f}) exp(-\lambda(X)) d\lambda= mes(i\mathcal{A}_{\tilde{M},F}^*) \Theta_{\tilde{\pi}}(\phi_{\tilde{M}}(\tilde{f}) \mathbf{1}_{H_{\tilde{M}}=X})$$

\noindent De plus la fonction $\phi_{\tilde{M}}(\tilde{f}) \mathbf{1}_{H_{\tilde{M}}=X}$ est bien définie à équivalence près.
\end{prop}

\ul{Preuve}: Cela découle de la proposition 6.4 alliée au théorème 5.5 et aux remarques 5.5(1) et (2) de [W4] $\blacksquare$

Soient $\tilde{f}\in C_c^\infty(\tilde{G}(F))$, $\tilde{M}\in \mathcal{L}^{\tilde{G}}$ et $\tilde{x}\in \tilde{M}(F)\cap \tilde{G}_{reg}(F)$. On définit l'intégrale orbitale pondérée invariante $I_{\tilde{M}}(\tilde{x},\tilde{f})$ par récurrence sur $a_{\tilde{M}}-a_{\tilde{G}}$ par la formule

$$\displaystyle I_{\tilde{M}}(\tilde{x},\tilde{f})=J_{\tilde{M}}(\tilde{x},\tilde{f})-\sum_{\tilde{L}\in \mathcal{L}(\tilde{M}); \tilde{L}\neq \tilde{G}} I_{\tilde{M}}^{\tilde{L}}(\tilde{x},\phi_{\tilde{L}}(\tilde{f})\mathbf{1}_{H_{\tilde{L}}-H_{\tilde{L}}(\tilde{x})})$$

\subsection{Quasi-caractères}

 Soit $\tilde{x}\in\tilde{G}_{ss}(F)$. En adaptant [W1] 3.1 au cas tordu, on définit la notion de bon voisinage $\omega\subset \mathfrak{g}_{\tilde{x}}(F)$. On note $Nil(\mathfrak{g}_{\tilde{x}})$ l'ensemble des orbites nilpotentes dans $\mathfrak{g}_{\tilde{x}}(F)$ et $Nil(\mathfrak{g}_{\tilde{x}})_{reg}$ le sous-ensemble des orbites régulières. Pour $\mathcal{O}\in Nil(\mathfrak{g}_{\tilde{x}})$ la transformée de Fourier de l'intégrale orbitale sur $\mathcal{O}$ est une fonction localement intégrable $X\mapsto \hat{j}(\mathcal{O},X)$. Elle dépend évidemment des choix de mesures et de la transformée de Fourier, mais si ces choix sont faits de façon cohérente comme en [B] 1.4, alors la fonction $\hat{j}(\mathcal{O},.)$ est bien définie. On suppose donc que ces choix ont été effectués comme en [B] 1.4. \\
 
 Soit $\Theta$ une fonction à valeurs complexes définie presque partout sur $\tilde{G}(F)$ qui est invariante par conjugaison par $G(F)$. On dit que c'est un quasi-caractère si pour tout $\tilde{x}\in\tilde{G}_{ss}(F)$, il existe un bon voisinage $\omega\subset \mathfrak{g}_{\tilde{x}}(F)$ et des nombres complexes $c_{\Theta,\mathcal{O}}(\tilde{x})$ pour $\mathcal{O}\in Nil(\mathfrak{g}_{\tilde{x}})$, tels que
 
$$\displaystyle\Theta(\tilde{x}exp(X))=\sum_{\mathcal{O}\in Nil(\mathfrak{g}_{\tilde{x}})}c_{\Theta,\mathcal{O}}(\tilde{x}) \hat{j}(\mathcal{O},X)$$

\noindent pour presque tout $X\in \omega$. D'après [C] théorème 3, pour toute représentation $\tilde{\pi}$ de $\tilde{G}(F)$ qui est $G(F)$-irréductible, le caractère $\Theta_{\tilde{\pi}}$ est un quasi-caractère. 

\subsection{Fonctions cuspidales et très cuspidales}

\begin{duf}
Soit $\tilde{f}\in C_c^\infty(\tilde{G}(F))$.
\begin{enumerate}
\item On dit que $\tilde{f}$ est très cuspidale si pour tout sous-groupe parabolique tordu propre $\tilde{P}=\tilde{M}U$ de $\tilde{G}$ et pour tout $\tilde{m}\in\tilde{M}(F)$, on a

$$\displaystyle\int_{U(F)} \tilde{f}(\tilde{m}u) du=0$$

\item On dit que $\tilde{f}$ est cuspidale si pour tout Levi tordu propre $\tilde{M}$ de $\tilde{G}$ et pour tout $\tilde{x}\in \tilde{G}_{reg}(F)\cap \tilde{M}(F)$ on a $J_{\tilde{G}}(\tilde{x},\tilde{f})=0$.
\end{enumerate}
\end{duf}

Soit $\tilde{f}\in C_c^\infty(\tilde{G}(F))$ une fonction très cuspidale. On peut lui associer un quasi-caractère $\Theta_{\tilde{f}}$ de la façon suivante. Soit $\tilde{x}\in \tilde{G}_{reg}(F)$ et notons $\tilde{M}(\tilde{x})$ le commutant de $A_{G_{\tilde{x}}}$ dans $\tilde{G}$. C'est un Levi de $\tilde{G}$ et on peut trouver $g\in G(F)$ de sorte que $g\tilde{M}(\tilde{x})g^{-1}$ soit un Levi semistandard. On pose alors

$$\Theta_{\tilde{f}}(\tilde{x})=(-1)^{a_{\tilde{M}(\tilde{x})}-a_{\tilde{G}}}D^{\tilde{G}}(\tilde{x})^{-1/2} J_{g\tilde{M}(\tilde{x})g^{-1}}(g\tilde{x}g^{-1},{}^g\tilde{f})$$

\noindent où ${}^g\tilde{f}(\tilde{x})=\tilde{f}(g^{-1}\tilde{x}g)$. Le résultat ne dépend pas du choix de $g$ (c'est une conséquence du lemme 5.2 de [W1] adapté au cas tordu). La fonction $\Theta_{\tilde{f}}$ ainsi obtenue est un quasi-caractère: il suffit de reprendre la preuve de [W1] corollaire 5.9 et de l'adapter au cas tordu. \\

En remplaçant dans ce qui précède intégrale orbitale pondérée par intégrale orbitale pondérée invariante, on peut aussi associer à toute fonction $\tilde{f}\in C_c^\infty(\tilde{G}(F))$ une fonction $I\Theta_{\tilde{f}}$: pour $\tilde{x}\in \tilde{G}_{reg}(F)$, on choisit $g\in G(F)$ tel que $g\tilde{M}(\tilde{x})g^{-1}\in\mathcal{L}^{\tilde{G}}$ et on pose

$$I\Theta_{\tilde{f}}(\tilde{x})=(-1)^{a_{\tilde{M}(\tilde{x})}-a_{\tilde{G}}}D^{\tilde{G}}(\tilde{x})^{-1/2} I_{g\tilde{M}(\tilde{x})g^{-1}}(g\tilde{x}g^{-1},\tilde{f})$$

\noindent le résultat ne dépend pas non plus du choix de $g$ (cela découle du fait que les intégrales orbitales pondérées invariantes sont bien invariantes et ne dépendent pas du choix de $K$). On a alors

\begin{prop}
Soit $\tilde{f}\in C_c^\infty(\tilde{G}(F))$ une fonction cuspidale. Alors $I\Theta_{\tilde{f}}$ est un quasi-caractère de $\tilde{G}(F)$
\end{prop}

\ul{Preuve}: C'est la même que celle de la proposition 2.5 de [W2], en utilisant le théorème 3 de [C] et la proposition suivante qui est le théorème 7.1 de [W4]

\begin{prop}
Soient $\tilde{f}\in C_c^\infty(\tilde{G}(F))$ une fonction cuspidale, $\tilde{M}$ un Levi tordu de $\tilde{G}$ et $\tilde{x}\in\tilde{M}(F)\cap \tilde{G}_{reg}(F)$, alors
\begin{enumerate}
\item Si $\tilde{x}\notin \tilde{M}(F)_{ell}$, $I_{\tilde{M}}(\tilde{x},\tilde{f})=0$;
\item Si $\tilde{x}\in \tilde{M}(F)_{ell}$,

$$\displaystyle D^{\tilde{G}}(\tilde{x})^{-1/2} (-1)^{a_{\tilde{M}}-a_{\tilde{G}}}I_{\tilde{M}}(\tilde{x},\tilde{f})=\sum_{\mathcal{O}\in\{\Pi_{ell}(\tilde{G})\}} \iota(\mathcal{O}) [i\mathcal{A}_{\mathcal{O}}^\vee: i\mathcal{A}_{\tilde{G},F}^\vee]^{-1} \int_{i\mathcal{A}_{\tilde{G},F}^*} \Theta_{\tilde{\pi}_\lambda}(\tilde{x}) \Theta_{\tilde{\pi}_\lambda^\vee}(\tilde{f}) d\lambda$$

\noindent où $\iota(\mathcal{O})$ sont certaines constantes qui sont celles définies en 1.6 dans le cas où $\tilde{G}$ vérifie les hypothèses de cette section.
\end{enumerate}
\end{prop}

On aura aussi besoin des propriétés suivantes :cf [W3] lemme 1.13.

\begin{lem}
\begin{enumerate}
\item Soit $\tilde{f}\in C_c^\infty(\tilde{G}(F))$ une fonction très cuspidale. Alors pour tout $\tilde{L}\in\mathcal{L}^{\tilde{G}}$, la fonction $\phi_{\tilde{L}}(\tilde{f})$ est cuspidale et on a l'égalité

$$\displaystyle \Theta_{\tilde{f}}=\sum_{\tilde{L}\in\mathcal{L}^{\tilde{G}}} |W^L||W^G|^{-1} (-1)^{a_{\tilde{L}}-a_{\tilde{G}}}Ind_{\tilde{L}}^{\tilde{G}}(I\Theta_{\phi_{\tilde{L}}(\tilde{f})})$$

\item Soit $\tilde{f}\in C_c^\infty(\tilde{G}(F))$ une fonction cuspidale. Alors il existe une fonction très cuspidale $\tilde{f}'\in C_c^\infty(\tilde{G}(F))$ telle que $\tilde{f}$ et $\tilde{f}'$ soit équivalentes (au sens défini en 1.7)
\end{enumerate}
\end{lem}

\section{Le groupe tordu de changement de base du groupe unitaire}

\subsection{Définition et description}

Soit $E/F$ une extension quadratique. On notera $N$ et $Tr$ la norme et la trace relatives à cette extension, $c:x\mapsto c(x)$ le $F$-automorphisme non trivial de $E$, $\chi_E$ le caractère quadratique de $F^\times$ de noyau $N(E^\times)$ et $\psi_E=\psi\circ Tr$ qui est un caractère additif non trivial de $E$. Soit $d\geqslant 1$ un entier et $V$ un espace vectoriel de dimension $d$ sur $E$. Notons $GL(V)$ le groupe des automorphismes $E$-linéaires de $V$ et $G=G_d=R_{E/F} GL(V)$ le groupe algébrique défini sur $F$ obtenu par restriction des scalaires à la Weil (dorénavant on notera $R_{E/F}$ la restriction des scalaires de $E$ à $F$). On note $V^c$ le groupe abélien $V$ muni de la structure de $E$-espace vectoriel définie par $\lambda.v=c(\lambda)v$. Introduisons $\tilde{G}=\tilde{G}_d=Isom(V,V^{c,*})$, la variété algébrique définie sur $F$ des isomorphismes $c$-linéaires de $V$ sur son dual. Le groupe $G$ agit à gauche et à droite sur $\tilde{G}$ de la façon suivante

$$(g,\tilde{x},g')\mapsto ({}^t g^c )^{-1}\tilde{x}g'$$

\noindent où ${}^t g$ est l'application transposée de $g$ et ${}^t g^c$ est l'application ${}^t g$ vue comme endomorphisme de $V^{*,c}$. Le couple $(G,\tilde{G})$ est alors un groupe tordu. On peut identifier $\tilde{G}(F)$ à l'ensemble des formes sesquilinéaires non dégénérées de $V$ (on adopte la convention qu'une telle forme est linéaire en la deuxième variable): à $\tilde{x}\in \tilde{G}(F)$ on associe la forme $(v,v')\mapsto <\tilde{x}(v),v'>$. \\

 Soit $\tilde{x}\in \tilde{G}(F)$ un élément semi-simple. On pose $x={}^t(\tilde{x}^c)^{-1}\tilde{x}$. On note $V''_{\tilde{x}}$ le noyau de $x-1$ dans $V$, $V'_{\tilde{x}}$ l'unique supplémentaire $x$-stable de $V''_{\tilde{x}}$ et $\tilde{G}'$, $\tilde{G}''$ les groupes tordus définis de la même façon que $\tilde{G}$ en changeant $V$ en $V'_{\tilde{x}}$ et $V''_{\tilde{x}}$ respectivement. La restriction de $\tilde{x}$ à $V''_{\tilde{x}}$ définit une forme hermitienne et $G''_{\tilde{x}}$ est le groupe unitaire de $(V''_{\tilde{x}},\tilde{x})$. On définit une fonction $\Delta$ sur $\tilde{G}_{ss}(F)$ par

$$\Delta(\tilde{x})=|N(det (1-x)_{|V'_{\tilde{x}}})|_F$$

Par restriction à $V'_{\tilde{x}}$, $\tilde{x}$ définit un élément de $\tilde{G}'(F)$ que l'on note encore $\tilde{x}$.
\begin{lem}
On a

$$D^{\tilde{G}}(\tilde{x})=|2|_F^{{d''}^2_{\tilde{x}}} \Delta(\tilde{x})^{d''_{\tilde{x}}}D^{\tilde{G}'}(\tilde{x})$$

\noindent où $d''_{\tilde{x}}=dim_E(V''_{\tilde{x}})$.
\end{lem}

\ul{Preuve}: Pour ne pas trop alourdir les notations, on posera $V'=V'_{\tilde{x}}$, $V''=V''_{\tilde{x}}$ et $d''=d''_{\tilde{x}}$. On a une décomposition $\mathfrak{g}(F)=\mathfrak{gl}_E(V'')\oplus Hom_E(V'',V')\oplus Hom_E(V',V'')\oplus \mathfrak{g}'(F)$. Posons $\mathfrak{g}_1=\mathfrak{gl}_E(V'')$ et $\mathfrak{g}_2=Hom_E(V'',V')\oplus Hom_E(V',V'')$. Alors $\mathfrak{g}_1$, $\mathfrak{g}_2$ et $\mathfrak{g}'(F)$ sont stables par $\theta_{\tilde{x}}$. Pour $i=1,2$, notons $\mathfrak{g}_{i,\tilde{x}}$ le sous-espace de $\mathfrak{g}_i$ où $\theta_{\tilde{x}}$ agit comme l'identité. On a alors

$$\mbox{(1)}\;\;\; D^{\tilde{G}}(\tilde{x})=\left(\prod_{j=1}^2|det(1-\theta_{\tilde{x}})_{|\mathfrak{g}_j/\mathfrak{g}_{j,\tilde{x}}}|_F\right) |det(1-\theta_{\tilde{x}})_{|\mathfrak{g}'/\mathfrak{g}'_{\tilde{x}}}|_F$$

\noindent Le dernier terme de ce produit vaut $D^{\tilde{G}'}(\tilde{x})$.

\begin{itemize}
\item Calcul de $|det(1-\theta_{\tilde{x}})_{|\mathfrak{g}_1/\mathfrak{g}_{1,\tilde{x}}}|_F$: sur le sous-espace $\mathfrak{g}_1$, $\theta_{\tilde{x}}$ est l'endomorphisme qui à $X\in End_E(V'')$ associe l'opposé de l'adjoint de $X$ pour la forme hermitienne $\tilde{x}_{|V''}$. Par conséquent, on a ${\theta_{\tilde{x}}}_{|\mathfrak{g}_1}^2=Id$. Notons $\mathfrak{g}_1^{\tilde{x}}$ le sous-espace de $\mathfrak{g}_1$ sur lequel $\theta_{\tilde{x}}=-Id$. On a alors $|det(1-\theta_{\tilde{x}})_{|\mathfrak{g}_1/\mathfrak{g}_{1,\tilde{x}}}|_F=|2|_F^{dim_F(\mathfrak{g}_1^{\tilde{x}})}$. Soit $\eta\in E^\times$ un élément non nul de trace nulle. Alors l'application $X\mapsto \eta X$ est un isomorphisme entre $\mathfrak{g}_1^{\tilde{x}}$ et $\mathfrak{g}_{1,\tilde{x}}$. Par conséquent on a $dim_F(\mathfrak{g}_1^{\tilde{x}})=dim_F(\mathfrak{g}_1)/2=d''^2$. D'où

$$\mbox{(2)}\;\;\; |det(1-\theta_{\tilde{x}})_{|\mathfrak{g}_1/\mathfrak{g}_{1,\tilde{x}}}|_F=|2|_F^{d''^2}$$

\item Calcul de $|det(1-\theta_{\tilde{x}})_{|\mathfrak{g}_2/\mathfrak{g}_{2,\tilde{x}}}|_F$: posons $\mathfrak{g}_3=Hom_E(V'',V')$ et $\mathfrak{g}_4=Hom_E(V',V'')$. Alors $\theta_{\tilde{x}}$ échange $\mathfrak{g}_3$ et $\mathfrak{g}_4$. On a donc

$$det(1-\theta_{\tilde{x}})_{|\mathfrak{g}_2/\mathfrak{g}_{2,\tilde{x}}}=det(1-\theta_{\tilde{x}}^2)_{|\mathfrak{g}_3/\mathfrak{g}_{3,\tilde{x}}}$$

\noindent L'action de $\theta_{\tilde{x}}^2$ sur $\mathfrak{g}_3$ est donné par $X\mapsto x\circ X$. Via le choix d'une base de $V''$, on a un isomorphisme $\mathfrak{g}_3\simeq V'^{\oplus d''}$ et via cet isomorphisme, $\theta_{\tilde{x}}^2$ agit par $x$ sur chaque facteur. Par conséquent

$$\mbox{(3)}\;\;\; |det(1-\theta_{\tilde{x}})_{|\mathfrak{g}_2/\mathfrak{g}_{2,\tilde{x}}}|_F=\Delta(\tilde{x})^{d''}$$
\end{itemize}

Les égalités (1), (2) et (3) entraînent le lemme $\blacksquare$

\vspace{2mm}

Fixons une base $(v_i)_{i=1,\ldots,d}$ de $V$, cela permet d'identifier $G$ à $R_{E/F}GL_d$. Pour $R$ une $F$-algèbre et $g\in (R_{E/F}GL_d)(R)=GL_d(R\otimes_F E)$, on note $g_{i,j}$ les coefficients de $g$ , $i,j=1,\ldots,d$ (ce sont des éléments de $R\otimes_F E$). Soient $B_d$, $U_d$ et $T_d$ les sous-groupes de $R_{E/F} GL_d$ constitués des matrices triangulaires supérieures, des matrices unipotentes triangulaires supérieures et des matrices diagonales respectivement. On pose $A_d=A_{T_d}$. On notera $K_d$ le sous-groupe de $(R_{E/F}GL_d)(F)=GL_d(E)$ des matrices à coefficients entiers dont le déterminant est de valuation nulle. On définit un caractère $\psi_d$ de $U_d(F)$ par la formule suivante

$$\displaystyle\psi_d: u\mapsto \psi_E(\sum_{i=1}^{d-1} u_{i,i+1})$$

Soit $\theta_d$ l'élément de $\tilde{G}(F)$ défini par $\theta_d(v_i)=(-1)^{i+[(d+1)/2]}v_{d-i}^*$ pour $i=1,\ldots,d$ où $(v_i^*)_{i=1,\ldots,d}$ est la base duale de $(v_i)_{i=1,\ldots,d}$. On notera aussi $\theta_d$ l'automorphisme de $G$ défini par $\theta_d g=\theta_d(g)\theta_d$. Alors $\theta_d$ laisse stable $B_d$, $U_d$, $T_d$ et préserve le caractère $\psi_d$ de $U_d(F)$. \\

Pour $\pi$ une représentation lisse de $G(F)$, on notera $\overline{\pi}$ la conjuguée complexe de $\pi$ et $\pi^c$ la composée de $\pi$ et de l'automorphisme de $G(F)$, $g=(g_{ij})\mapsto (c(g_{ij}))$. La classe d'isomorphisme de $\pi^c$ ne dépend pas du choix de la base $(v_i)_{i=1,\ldots,d}$.

\subsection{Modèle de Whittaker}

Fixons à nouveau une base $(v_i)_{i=1,\ldots,d}$ de $V$. Soit $(\pi, E_\pi)$ une représentation lisse irréductible de $G(F)$. Appelons fonctionnelle de Whittaker pour $\pi$ toute forme linéaire $\ell: E_\pi\to \mathbb{C}$ vérifiant $\ell(\pi(u)e)=\psi_d(u)\ell(e)$ pour tout $e\in E_\pi$ et pour tout $u\in U_d(F)$. L'espace des fonctionnelles de Whittaker pour $\pi$ est de dimension au plus $1$. Supposons que cet espace soit de dimension $1$. Soit $\theta_d(\pi)$ la représentation $g\mapsto \pi(\theta_d(g))$. Alors $\pi$ s'étend en une représentation tordue $\tilde{\pi}$ de $\tilde{G}(F)$ si et seulement si $\pi\simeq \theta_d(\pi)$. Supposons cette condition vérifiée et soit $\tilde{\pi}$ une $\tilde{G}(F)$-représentation tordue qui prolonge $\pi$. L'application $\phi\mapsto \phi\circ \tilde{\pi}(\theta_d)$ préserve la droite des fonctionnelles de Whittaker et agit par un scalaire non nul que l'on note $w(\tilde{\pi},\psi)$. Un calcul de changement de base montre que $w(\tilde{\pi},\psi)$ ne dépend pas du choix de la base. En revanche, il dépend de $\psi$. Pour $b\in F^\times$, on note $\psi^b$ le caractère défini par $\psi^b(x)=\psi(bx)$. Soit $D^b$ la matrice diagonale de coefficients $D^b_{ii}=b^{d-i}$ et $\phi$ une fonctionnelle de Whittaker non nulle pour $\psi$. Alors $\phi\circ \pi(D^b)$ est une fonctionnelle de Whittaker pour $\psi^b$ et on a $\theta_d(D^b)(D^b)^{-1}=b^{1-d}I_d$. On en déduit que \\

 (1) $w(\tilde{\pi},\psi^b)=\omega_\pi(b)^{1-d}w(\tilde{\pi},\psi)$ pour tout $b\in F^\times$. \\
 
\noindent où $\omega_\pi$ est le caractère central de $\pi$. \\

 (2) Supposons $\pi$ tempérée. Alors $w(\tilde{\pi}^\vee,\psi^-)w(\tilde{\pi},\psi)=1$ où $\psi^-=\psi^{-1}$. Si de plus $\tilde{\pi}$ est unitaire alors $|w(\tilde{\pi},\psi)|=1$. \\
 
\noindent En effet, soit $\phi$ une fonctionnelle de Whittaker non nulle de $\pi$ pour le caractère $\psi$ et $\phi^-$ une fonctionnelle de Whittaker non nulle de $\pi^\vee$ pour le caractère $\psi^-$. D'après [W2] lemme 3.6, il existe un nombre complexe non nul $C$ tel que pour tout $e\in E_\pi$ et pour tout $e^\vee\in E_{\pi^\vee}$ on ait

$$\displaystyle <e^\vee,e>= C\int_{U_{d-1}(F)\backslash GL_{d-1}(E)} \phi^-(\pi^\vee(g)e^\vee) \phi(\pi(g)e)dg$$

\noindent où on a identifié $R_{E/F}GL_{d-1}$ au sous-groupe des éléments de $G$ qui laissent stables $v_d$ et $v_d^*$. De l'égalité $<\tilde{\pi}^\vee(\theta_d)e^\vee,\tilde{\pi}(\theta_d)e>=<e^\vee,e>$, l'on déduit la première identité. Pour la deuxième, il suffit de remarquer que si $\tilde{\pi}$ est unitaire, on a 

$$\overline{w(\tilde{\pi},\psi)}=w(\overline{\tilde{\pi}},\overline{\psi})=w(\tilde{\pi}^\vee,\psi^-)$$

\subsection{Représentations induites}

Soit $\tilde{L}$ un Levi tordu de $\tilde{G}$. Il existe alors une décomposition

$$V=V_u\oplus\ldots\oplus V_1\oplus V_0\oplus V_{-1}\oplus\ldots\oplus V_{-u}$$

\noindent telle que $\tilde{L}$ soit l'ensemble des $\tilde{x}\in\tilde{G}$ qui vérifient $\tilde{x}(V_j)=V^*_{-j}$ pour $j=-u,\ldots,u$. Pour $\tilde{x}\in\tilde{L}$, on notera $\tilde{x}_j$ la restriction de $\tilde{x}$ à $V_j$, $j=-u,\ldots,u$. On a alors

$$\mbox{(1)}\;\;\; L=R_{E/F} GL(V_u)\times\ldots\times R_{E/F} GL(V_1)\times R_{E/F} GL(V_0)\times R_{E/F} GL(V_{-1})\times\ldots\times R_{E/F} GL(V_{-u})$$

\noindent Posons $d_j=dim(V_j)=d_{-j}$, $j=0,\ldots,u$. Pour $j=-u,\ldots,u$, fixons une base $(v^j_i)_{i=1,\ldots,d_j}$ de $V_j$. La famille $(v_i)_{i=1,\ldots,d}$ définie par $v_{d_u+\ldots+d_{j+1}+i}=v^j_i$ pour tout $j=-u,\ldots,u$ et pour tout $i=1,\ldots,d_j$, est une base de $V$. Grâce à cette base, on définit comme en 2.1 un élément $\theta_d$. On remarque que $\theta_d\in \tilde{L}(F)$. Notons $\theta_{V_j}$ la restriction de $\theta_d$ à $R_{E/F} GL(V_j)$. Soit $\tilde{\tau}\in Temp(\tilde{L})$. Conformément à (1), on peut écrire

$$\tau=\tau_u\otimes\ldots\otimes \tau_{-u}$$

\noindent où pour tout $j=-u,\ldots,u$, $\tau_j\in Temp(GL(V_j))$. Puisque $\tau$ s'étend en une représentation de $\tilde{L}(F)$, on a $\theta_{V_j}(\tau_j)\simeq \tau_{-j}$ pour $j=1,\ldots,u$. Choisissons des opérateurs unitaires $A_j:E_{\tau_j}\to E_{\tau_{-j}}$ tels que

$$A_j \tau_j(\theta_{V_j}(x_{-j}))=\tau_{-j}(x_{-j}) A_j$$

\noindent pour tout $x_{-j}\in GL(V_{-j})$. Notons $\tilde{G}_0$ l'analogue de $\tilde{G}$ lorsque l'on remplace $V$ par $V_0$. Il existe alors un unique prolongement $\tilde{\tau}_0$ de $\tau_0$ à $\tilde{G}_0(F)$ tel que pour $e=e_u\otimes \ldots\otimes e_1\otimes e_0\otimes e_{-1}\otimes\ldots\otimes e_{-u}\in E_{\tau}$, on ait

$$\tilde{\tau}(\theta_d)e=A_u^{-1}e_{-u}\otimes\ldots\otimes A_1^{-1}e_{-1}\otimes \tilde{\tau}_0(\theta_{V_0})e_0\otimes A_1e_1\otimes\ldots\otimes A_ue_u$$

\noindent et ce prolongement ne dépend pas du choix des $A_j$. Introduisons la représentation $\tilde{\pi}=i_Q^G(\tilde{\tau})$ ($\tilde{Q}\in \mathcal{P}(\tilde{L})$). On a alors

$$\mbox{(2)}\;\;\; w(\tilde{\pi},\psi)=w(\tilde{\tau}_0,\psi)$$

\noindent On peut en effet construire de façon explicite une fonctionnelle de Whittaker pour $\pi$ à partir d'une fonctionnelle de Whittaker pour $\tau$. Plus précisément, notons $w$ l'élément qui envoie $v^j_i$ sur $v^{-j}_i$ pour $j=-u,\ldots,u$ et pour $i=1,\ldots,d_j$. Remarquons que $\theta_d(w)=w$. On a fixé une base de $V_j$ pour $j=-u,\ldots,u$, ce qui permet de parler de fonctionnelle de Whittaker pour une représentation lisse de $GL(V_j)$. Fixons des fonctionnelles de Whittaker non nulles $\phi_j: E_{\tau_j}\to \mathbb{C}$ de $\tau_j$ pour $j=-u,\ldots,u$. On peut toujours supposer, ce que l'on fait, que l'on a $\phi_j\circ A_j=\phi_{-j}$ pour $j=1,\ldots,u$. Choisissons une suite exhaustive $(U(F)_N)_{N\geqslant 1}$ de sous-groupes compacts-ouverts de $U(F)$. Définissons l'application linéaire $\phi:E_\pi=E^G_{Q,\tau}\to \mathbb{C}$ par la formule

$$\displaystyle \phi(f)=\lim\limits_{N\to\infty} \int_{U(F)_N} \left(\phi_u\otimes\ldots\phi_{-u}\right)\left(f(wu)\right)\overline{\psi}(u) du$$

\noindent Cette expression a un sens et $\phi$ est une fonctionnelle de Whittaker non nulle pour $\pi$ (cf 3.1 [Sh]). Cette définition permet d'exprimer $\phi\circ\tilde{\pi}(\theta_d)$ de deux façons différentes, et on obtient alors la relation (2). \\

 Comme cela nous sera utile plus tard, déterminons le caractère de $\tilde{\tau}$ en fonction du caractère de $\tilde{\tau}_0$ et des caractères de $\tau_j$ pour $j=1,\ldots,u$.

\begin{lem}
Pour presque tout $\tilde{x}\in\tilde{L}(F)$, on a l'égalité

$$\displaystyle \Theta_{\tilde{\tau}}(\tilde{x})=\Theta_{\tilde{\tau}_0}(\tilde{x}_0)\prod_{j=1}^u \Theta_{\tau_j}((-1)^{d+1} \; {}^t(\tilde{x}^c_{-j})^{-1}\tilde{x}_j)$$
\end{lem}

\ul{Preuve}: Soient $f_j\in C_c^\infty(GL(V_j))$, $j=-u,\ldots,u$. Définissons $\tilde{f}\in C_c^\infty(\tilde{L}(F))$ par $\tilde{f}(x\theta_d)=\prod_{j=-u}^u f_j(x_j)$, où on a noté $x_j$ la restriction de $x$ à $V_j$ pour $x\in L(F)$. L'opérateur $\tilde{\tau}(\tilde{f})$ envoie alors $e=e_u\otimes\ldots\otimes e_{-u}\in E_{\tau}$ sur

$$\tau_u(f_u)A_u^{-1}e_{-u}\otimes\ldots\otimes \tau_1(f_1)A_1^{-1}e_{-1}\otimes \tau_0(f_0)\tilde{\tau}_0(\theta_{V_0}) e_0\otimes \tau_{-1}(f_{-1})A_{1}e_{1}\otimes\ldots\otimes \tau_{-u}(f_{-u})A_{u}e_{u}$$

\noindent Soit $j\in\{1,\ldots,u\}$. La trace de l'endomorphisme de $E_{\tau_j}\otimes E_{\tau_{-j}}$ qui envoie $e_j\otimes e_{-j}$ sur $\tau_j(f_j)A_j^{-1}e_{-j}\otimes \tau_{-j}(f_{-j})A_je_j$ est égale à la trace de l'opérateur $\tau_j(f_{-j}^\theta \star f_j)$ où $f_{-j}^\theta(x_j)=f_{-j}(\theta_{V_{-j}}(x_j))$ et $\star$ est la convolution. Cette trace est égale à

$$\displaystyle \int_{GL(V_j)\times GL(V_{-j})} \Theta_{\tau_j}(x_j\theta_{V_{-j}}(x_{-j})) f_{-j}(x_{-j}) f_j(x_j) dx_{-j}dx_j$$

\noindent La trace de l'opérateur $\tau_0(f_0)\tilde{\tau}_0(\theta_{V_0})$, elle vaut

$$\displaystyle \int_{G_0(F)} \Theta_{\tilde{\tau}_0}(x_0\theta_{V_0}) f_0(x_0)dx_0$$

\noindent On en déduit que la trace de l'opérateur $\tilde{\tau}(\tilde{f})$ est égale à

$$\displaystyle \int_{L(F)} \Theta_{\tilde{\tau}_0}(x_0\theta_{V_0})\prod_{j=1}^u \Theta_{\tau_j}(x_j\theta_{V_{-j}}(x_{-j})) \tilde{f}(\tilde{x}) dx$$

\noindent où on a posé $\tilde{x}=x\theta_d$ pour $x\in L(F)$. On vérifie facilement que $x_0\theta_{V_0}=\tilde{x}_0$ et $x_j\theta_{V_{-j}}(x_{-j})=(-1)^{d+1} \; {}^t(\tilde{x}^c_{-j})^{-1}\tilde{x}_j$ pour tout $x\in L(F)$. On en déduit le résultat annoncé $\blacksquare$

\subsection{Représentations elliptiques}

Soit $\tilde{L}$ un Levi tordu de $\tilde{G}$ et reprenons les notations de la section précédente. On a notamment une décomposition associée à $\tilde{L}$

$$V=V_u\oplus\ldots\oplus V_1\oplus V_0\oplus V_{-1}\oplus\ldots V_{-u}$$

\noindent et $\tilde{G}_0$ désigne le groupe tordu analogue de $\tilde{G}$ lorsque l'on remplace $V$ par $V_0$. Considérons un Levi $L_0$ de $G_0$. Ce Levi admet une décomposition

$$L_0=R_{E/F} GL_{d'_1}\times \ldots \times R_{E/F}GL_{d'_s}$$

\noindent Pour $j=1,\ldots,s$, soit $\rho_j$ une représentation irréductible de la série discrète de $GL_{d'_j}(E)$ telle que $\theta_{d'_j}(\rho_j)\simeq\rho_j$. On suppose que $\rho_i\not\simeq \rho_j$ pour $i\neq j$. Posons $\tau_0=i^{G_0}_{Q_0}(\rho_1\otimes\ldots\otimes \rho_s)$ ($Q_0\in\mathcal{P}^{G_0}(L_0)$). Pour $j=1,\ldots,u$, soit $\tau_j$ une représentation irréductible de la série discrète de $GL(V_j)$. Posons

$$\tau=\tau_u\otimes\ldots\otimes \tau_1\otimes\tau_0\otimes\theta^1(\tau_1)\otimes\ldots\otimes\theta^u(\tau_u)$$

\noindent Alors $\tau$ se prolonge en une représentation unitaire $\tilde{\tau}$ de $\tilde{L}(F)$ et c'est une représentation elliptique de $\tilde{L}(F)$. Posons $s(\tilde{\tau})=s$. Ce terme dépend seulement de l'orbite sous $i\mathcal{A}_{\tilde{L}}^*$ de $\tilde{\tau}$. Si $\mathcal{O}\in\{\Pi_{ell}(\tilde{L})\}$ est cette orbite, on peut donc poser $s(\mathcal{O})=s$. On a alors $\iota(\mathcal{O})=2^{-s(\mathcal{O})-a_{\tilde{L}}}$ (cf 1.6). Toutes les représentations elliptiques de $\tilde{L}(F)$ s'obtiennent de cette façon.

\subsection{Facteurs $\epsilon$}

Soient $m,d\geqslant 1$ des entiers, et posons $H=G_m$. Soient $\pi$, $\sigma$ des représentations irréductibles tempérées de $G(F)$ et $H(F)$ respectivement. Posons

$$\epsilon(\pi\times \sigma,\psi_E)=\epsilon(1/2,\pi\times\sigma,\psi_E)$$

\noindent où la fonction $s\mapsto \epsilon(s,\pi\times\sigma,\psi_E)$ est celle définie dans [JPSS]. On a alors les propriétés suivantes dont les preuves découlent de la définition ou alors se trouvent dans [JPSS]. Pour des représentations $\tau_i\in Temp(G_{d_i})$, $i=1,\ldots,t$, avec $d=d_1+\ldots+d_t$, on note $\tau_1\times\ldots\times \tau_t$ la représentation $i_P^{G}(\tau_1\otimes\ldots\otimes \tau_t)$ de $G(F)$ où $P$ est un sous-groupe parabolique de $G$ de Levi

$$G_{d_1}\times\ldots\times G_{d_t}$$

\begin{prop}
Notons $\omega_\pi$ et $\omega_\sigma$ leurs caractères centraux. On a les propriétés suivantes
\begin{enumerate}
\item $\epsilon(\pi\times \sigma,\psi_E^b)=\omega_\pi(b)^m \omega_\sigma(b)^d \epsilon(\pi\times \sigma,\psi_E)$ pour tout $b\in E^\times$
\item $\epsilon(\pi^c\times \sigma^c,\psi_E)=\epsilon(\pi\times \sigma,\psi_E)$
\item $\epsilon(\pi^\vee\times \sigma^\vee,\psi_E)\epsilon(\pi\times \sigma,\psi_E)=\omega_\pi(-1)^m \omega_\sigma(-1)^d$
\item Si $\pi=\tau_1\times\ldots\times \tau_t$, alors $\displaystyle \epsilon(\pi\times\sigma,\psi_E)=\prod_{j=1}^t \epsilon(\tau_j\times \sigma,\psi_E)$
\end{enumerate}
\end{prop}

Supposons maintenant que $m$ et $d$ sont de parités différentes et soient $\tilde{\pi}\in Temp(\tilde{G})$ et $\tilde{\sigma}\in Temp(\tilde{H})$. Posons

$$\mbox{(1)}\;\;\; \epsilon_{\nu}(\tilde{\pi},\tilde{\sigma})=\epsilon_{-\nu}(\tilde{\sigma},\tilde{\pi})=w(\tilde{\pi},\psi)w(\tilde{\sigma},\psi) \omega_\pi((-1)^{[m/2]}2\nu)\omega_\sigma((-1)^{1+[d/2]}2\nu)\epsilon(\pi\times\sigma,\psi_E)$$

\noindent Ce terme ne dépend pas du choix de $\psi:F\to\mathbb{C}^\times$ d'après 2.2(1) et le 1. de la proposition précédente. Soit $\tilde{L}$ un Levi tordu de $\tilde{G}$ que l'on écrit comme en 2.3. Soit $\tilde{\tau}\in Temp(\tilde{L})$. On peut écrire

$$\tau=\tau_u\otimes\ldots\otimes\tau_1\otimes\tau_0\otimes \theta^1(\tau_1)\otimes\ldots\otimes\theta^u(\tau_u)$$

\noindent Le prolongement $\tilde{\tau}$ de $\tau$ détermine un prolongement $\tilde{\tau}_0$ de $\tau_0$ comme en 2.3. Posons

$$\displaystyle \epsilon_\nu(\tilde{\tau},\tilde{\sigma})=\epsilon_\nu(\tilde{\tau}_0,\tilde{\sigma})\prod_{j=1}^u\omega_{\tau_j}(-1)^m$$

\noindent Supposons que $\tilde{\pi}=i_Q^G(\tilde{\tau})$ ($\tilde{Q}\in \mathcal{P}(\tilde{L})$). On a alors

$$\mbox{(2)}\;\;\; \epsilon_\nu(\tilde{\pi},\tilde{\sigma})=\epsilon_\nu(\tilde{\tau},\tilde{\sigma})$$

\noindent Cela découle en effet de 2.3(2) et des points 2., 3. et 4 de la proposition 2.5.1.

\section{La formule géométrique}

\subsection{La situation}

Soit $V$ un espace vectoriel de dimension $d$ sur $E$. On se donne une décomposition $V=Z\oplus W$ où $W$ et $Z$ sont de dimensions $m$ et $2r+1$ respectivement, $m,r\in\mathbb{N}$. Fixons une base $(z_i)_{i=-r,\ldots,r}$ de $Z$ et posons $V_0=W\oplus Ez_0$. Introduisons les groupes et groupes tordus $G=R_{E/F} GL(V)$, $H=R_{E/F} GL(W)$, $G_0=R_{E/F} GL(V_0)$, $\tilde{G}=Isom(V,V^{c,*})$ et $\tilde{H}=Isom(W,W^{c,*})$. On identifie $H$ au sous-groupe des éléments de $G$ qui agissent trivialement sur $Z$. Fixons $\nu\in F^\times$ et notons $\tilde{\zeta}$ l'élément de $Isom(Z,Z^{c,*})$ défini par $\tilde{\zeta}(z_i)=(-1)^i 2\nu z_{-i}^*$ où $(z_i^*)_{i=-r,\ldots,r}$ est la base duale de la base $(z_i)_{i=-r,\ldots,r}$. La décomposition $V=Z\oplus W$ induit une décomposition $V^{c,*}=Z^{c,*}\oplus W^{c,*}$. On a un plongement $\iota: \tilde{H}\hookrightarrow \tilde{G}$ défini de la façon suivante: pour $\tilde{x}\in \tilde{H}$, $\iota(\tilde{x})$ est l'élément dont la restriction à $W$ est $\tilde{x}$ et la restriction à $Z$ est $\tilde{\zeta}$. En pratique on ne fait pas de distinction entre $\tilde{x}$ et $\iota(\tilde{x})$. Considérons le sous-groupe parabolique $P$ de $G$ qui conserve le drapeau

$$Ez_r\subset \ldots\subset Ez_r\oplus\ldots\oplus Ez_1\subset Ez_r\oplus\ldots\oplus Ez_1\oplus V_0 $$
$$\subset Ez_r\oplus\ldots\oplus Ez_1\oplus V_0\oplus Ez_{-1}\subset \ldots \subset Ez_r\oplus\ldots\oplus Ez_1\oplus V_0\oplus Ez_{-1}\oplus\ldots\oplus Ez_{-r}$$

On note $A$ le sous-tore de $G$ des éléments qui préservent les droites $Ez_i$ pour $i=\pm 1,\ldots,\pm r$ et agissent comme l'identité sur $V_0$ et $U$ le radical unipotent de $P$. On pose $M=AG_0$, c'est une composante de Levi de $P$. Pour $a\in A(F)$, on note $a_i$ la valeur propre de $a$ agissant sur $z_i$ pour $i=\pm 1,\ldots, \pm r$. On note $\tilde{P}$ le normalisateur de $P$ dans $\tilde{G}$ et $\tilde{M}$ le normalisateur de $M$ dans $\tilde{P}$. Ce sont respectivement un sous-groupe parabolique tordu et un Levi tordu de $\tilde{G}$ et on a $\tilde{H}\subset \tilde{M}$. On définit un caractère $\xi$ de $U(F)$ par

$$\displaystyle\xi(u)=\psi_E(\sum_{i=-r}^{r-1} u_{i+1,i})$$

\noindent où on a posé $u_{i+1,i}=<z_{i+1}^*,uz_i>$ pour $i=-r,\ldots,r-1$. ce caractère est invariant par $\theta_{\tilde{y}}$ pour tout $\tilde{y}\in \tilde{H}(F)$. \\

Fixons un $\mathcal{O}_E$-réseau $R$ de $V$ qui est somme d'un $\mathcal{O}_E$-réseau de $W$ et du réseau engendré par les $z_i$ pour $i=0,\pm 1,\ldots, \pm r$. Notons $K$ le stabilisateur dans $G(F)$ de ce réseau. On a alors $G(F)=P(F)K$. Soit $R^\vee$ le réseau dual de $R$ dans $V^*$. Pour $N\geqslant 0$ un entier, on définit $\Omega_N$ comme le sous-ensemble de $G(F)$ des éléments qui s'écrivent $g=uag_0k$ où $u\in U(F)$, $k\in K$, $a\in A(F)$ et $g_0\in G_0(F)$ vérifient

\begin{itemize}
\item $|val_E(a_r)|\leqslant N$ et $|val_E(a_i)|\leqslant 2N$ pour $i=\pm 1,\ldots, \pm (r-1),-r$
\item $g_0^{-1}z_0\in \pi_E^{-2N}R$ et ${}^tg_0^{-1} z_0^*\in \pi_E^{-2N} R^\vee$
\end{itemize}

\noindent On note $\kappa_N$ le fonction caractéristique de l'ensemble $\Omega_N$. C'est une fonction invariante à gauche par $H(F)U(F)$ et invariante à droite par $K$.

\subsection{Un ensemble de tores tordus}

Considérons une décomposition $W=W'\oplus W''$ et notons $H'=R_{E/F} GL(W')$, $H''=R_{E/F} GL(W'')$, $\tilde{H}'=Isom(W',{W'}^{c,*})$, $\tilde{H}''=Isom_c(W'',{W''}^{c,*})$. Soient $\tilde{T}'$ un tore maximal anisotrope de $\tilde{H}'$ (i.e. tel que $A_{\tilde{T}'}=\{1\}$) et $\tilde{\zeta}_{H,T}\in \tilde{H}''(F)$ une forme hermitienne non dégénérée. On suppose que les groupes unitaires de $\tilde{\zeta}_{H,T}$ et $\tilde{\zeta}_{G,T}=\tilde{\zeta}\oplus \tilde{\zeta}_{H,T}$ sont quasi-déployés. Notons $\tilde{T}$ l'ensemble des éléments $\tilde{x}\in \tilde{H}$ qui envoient $W'$ sur ${W'}^{c,*}$, $W''$ sur ${W''}^{c,*}$ et tels que la restriction de $\tilde{x}$ à $W'$ appartienne à $\tilde{T}'$ et la restriction de $\tilde{x}$ à $W''$ soit égale à $\tilde{\zeta}_{H,T}$. On note $\underline{\mathcal{T}}$ l'ensemble des sous-ensembles $\tilde{T}$ de $\tilde{H}$ obtenus de cette façon. Cet ensemble est stable par conjugaison par $H(F)$ et on fixe un ensemble $\mathcal{T}$ de représentants des classes de conjugaison. Pour $\tilde{T}\in\underline{\mathcal{T}}$ on note $T$ l'ensemble des éléments de $T'$ dont on prolonge l'action sur $W$ en les laissant agir comme l'identité sur $W''$. C'est un sous-tore (en général non maximal) de $H$. Les éléments de $\tilde{T}$ définissent un automorphisme $\theta$ de $T$ et on note $T^\theta$, $T_\theta$ le sous-groupe des point fixe et sa composante neutre respectivement. On a $T^\theta=T_\theta$. On notera aussi $\tilde{T}(F)/\theta$ l'espace des classes de conjugaison par $T(F)$ dans $\tilde{T}(F)$. Le normalisateur de $\tilde{T}$ dans $H(F)$ contient $T(F)\times U(\tilde{\zeta}_{H,T})(F)$, on pose

$$W(H,\tilde{T})=Norm_{H(F)}(\tilde{T})/(T(F)\times U(\tilde{\zeta}_{H,T})(F))$$

Pour $\tilde{T}\in\underline{\mathcal{T}}$, on notera dorénavant avec un indice $T$ les objets définis dans ce paragraphe qui dépendent de $\tilde{T}$: $W'_T$, $W''_T$, $H'_T$, $\tilde{H}'_T$, $H''_T$, $\tilde{H}''_T$ etc...

\subsection{Quelques fonctions nécessaires à la formule géométrique}

Soit $\tilde{T}\in\underline{T}$. On note $\tilde{T}_\sharp$ l'ouvert de Zariski des éléments $\tilde{t}\in \tilde{T}$ tels que les valeurs propres de $t={}^t(\tilde{t}^c)^{-1}\tilde{t}$ sur $W'_T$ soient toutes distinctes et différentes de $1$. Soient $\Gamma$ et $\Theta$ des quasi-caractères de $\tilde{G}(F)$ et $\tilde{H}(F)$ respectivement. Soit $\tilde{t}\in\tilde{T}_\sharp(F)$. On a alors $G_{\tilde{t}}=U(\tilde{\zeta}_{G,T})\times T_\theta$ et $H_{\tilde{t}}=U(\tilde{\zeta}_{H,T})\times T_\theta$. Par conséquent, on a des identifications $Nil(\mathfrak{g}_{\tilde{t}})=Nil(\mathfrak{u}(\tilde{\zeta}_{G,T}))$ et $Nil(\mathfrak{h}_{\tilde{t}})=Nil(\mathfrak{u}(\tilde{\zeta}_{H,T}))$. Posons

$$\displaystyle c_\Gamma(\tilde{t})=\frac{1}{|Nil(\mathfrak{u}(\tilde{\zeta}_{G,T}))_{reg}|} \sum_{\mathcal{O}\in Nil(\mathfrak{u}(\tilde{\zeta}_{G,T}))_{reg}} c_{\Gamma,\mathcal{O}}(\tilde{t})$$

\noindent et

$$\displaystyle c_\Theta(\tilde{t})=\frac{1}{|Nil(\mathfrak{u}(\tilde{\zeta}_{H,T}))_{reg}|} \sum_{\mathcal{O}\in Nil(\mathfrak{u}(\tilde{\zeta}_{H,T}))_{reg}} c_{\Theta,\mathcal{O}}(\tilde{t})$$

On vérifie facilement que les fonctions $c_\Gamma$ et $c_\Theta$ sont invariantes par conjugaison par $T(F)$. Elles définissent donc des fonctions sur $\tilde{T}(F)/\theta$.

\subsection{Un critère de convergence}

Soient $\tilde{T}\in \underline{\mathcal{T}}$, $\tilde{t}\in\tilde{T}(F)$ et posons $t={}^t(\tilde{t}^c)^{-1} \tilde{t}\in H(F)$. Notons $E'(\tilde{t})$ resp. $E''(\tilde{t})$ l'image respectivement le noyau de $t-1$ dans $W'_T$. Introduisons les groupes $J'(\tilde{t})=R_{E/F} GL(E'(\tilde{t}))$ et $J''(\tilde{t})=R_{E/F} GL(E''(\tilde{t}))$. On note $\mathfrak{z}_{\tilde{t}}$ le centre de l'algèbre de Lie de $J'(\tilde{t})_{\tilde{t}}$. On a une inclusion naturelle $J'(\tilde{t})\times J''(\tilde{t})\subset H'_T$. Le tore $T_\theta$ est alors un sous-tore maximal de $H'_{T,\tilde{t}}=J'(\tilde{t})_{\tilde{t}}\times J''(\tilde{t})_{\tilde{t}}$. En particulier, on a $\mathfrak{z}_{\tilde{t}}\subset \mathfrak{t}_\theta$.

Pour $X$ un $F$-espace vectoriel et $i\in\mathbb{R}$, on note $C_i(X)$ l'espace des fonctions mesurables $\varphi: X\to \mathbb{C}$ vérifiant $\varphi(\lambda x)=|\lambda|_F^i\varphi(x)$ pour tout $x\in X$ et tout $\lambda\in F^{\times 2}$. On note $C_{\geqslant i}(X)$ le sous-espace vectoriel des fonctions à valeurs complexes sur $X$ engendré par les $C_j(X)$ pour $j\geqslant i$. Soit $\chi:F^{\times}\to \{\pm 1\}$ un caractère quadratique, pour $i\in \mathbb{R}$ on notera $C_{i,\chi}(X)$ l'espace des fonctions mesurables $\varphi:X\to \mathbb{C}$ telle que $\varphi(\lambda x)=\chi(\lambda)|\lambda|_F^i\varphi(x)$ pour tout $x\in X$ et tout $\lambda\in F^{\times}$. Enfin, $C_{\geqslant i,\chi}(X)$ désignera l'espace de fonctions sur $X$ engendré par $C_{i,\chi}(X)$ et les $C_j(X)$ pour $j>i$. \\
 
 Soit $\delta: \tilde{T}(F)/\theta\to \mathbb{R}$ une fonction, on définit alors $C_{\geqslant \delta}(\tilde{T})$ (resp. $C_{\geqslant \delta,\chi}(\tilde{T})$) comme l'espace des fonctions complexes $f$ définies presque partout sur $\tilde{T}(F)/\theta$ et telles que pour tout $\tilde{t}\in \tilde{T}(F)$, il existe un bon voisinage $\omega\subset \mathfrak{t}_\theta(F)$ de 0 et une fonction $\varphi\in C_{\geqslant \delta(\tilde{t})}(\mathfrak{t}_\theta(F)/\mathfrak{z}_{\tilde{t}}(F))$ (resp. $\varphi\in C_{\geqslant \delta(\tilde{t}),\chi}(\mathfrak{t}_\theta(F)/\mathfrak{z}_{\tilde{t}}(F))$) vérifiant

\begin{center}
$f(\tilde{t}exp(X))=\varphi(\overline{X})$ pour presque tout $X\in\omega$
\end{center}

\noindent où $\overline{X}$ désigne la projection de $X$ sur $\mathfrak{t}_\theta(F)/\mathfrak{z}_{\tilde{t}}(F)$. Rappelons que l'on a défini en 2.1 une fonction $\Delta$ sur $\tilde{G}_{ss}(F)$. On montre exactement de la même façon que le lemme 5.2.2 de [B] le critère de convergence suivant

\vspace{3mm}

\begin{lem}
Soit $\tilde{T}\in \underline{\mathcal{T}}$. Posons

$$\displaystyle \delta_0(\tilde{t})=inf\big(\frac{dim(\mathfrak{z}_{\tilde{t}})-dim(\mathfrak{t}_\theta)-dim(E''(\tilde{t}))}{2}, -\frac{1}{2}\big)$$

\noindent pour tout $\tilde{t}\in\tilde{T}(F)$. Soit $f\in C_{\geqslant\delta_0, \chi_E}(\tilde{T})$, alors pour tout $s\in\mathbb{C}$ tel que $Re(s)>0$ l'intégrale

$$\displaystyle\int_{\tilde{T}(F)/\theta} f(\tilde{t})\Delta(\tilde{t})^s d\tilde{t}$$

\noindent converge absolument et de plus la limite 

$$\lim\limits_{s\to 0^+} \displaystyle\int_{\tilde{T}(F)/\theta} f(\tilde{t})\Delta(\tilde{t})^s d\tilde{t}$$
existe
\end{lem}

Soient $\Gamma$ et $\Theta$ des quasi-caractères de $\tilde{G}(F)$ et $\tilde{H}(F)$ respectivement. Une démonstration tout à fait analogue à celle du lemme 5.3.1 de [B] montre que la fonction $\tilde{t}\mapsto c_\Gamma(\tilde{t})c_\Theta(\tilde{t})D^{\tilde{H}}(\tilde{t})^{1/2}D^{\tilde{G}}(\tilde{t})^{1/2}\Delta(\tilde{t})^{-1/2}$ est dans $C_{\geqslant \delta_0,\chi_E}(\tilde{T})$. On peut par conséquent définir

\[\begin{aligned}
\mbox{(1)}\;\;\; \displaystyle J_{geom}(\Theta,\Gamma)= \sum_{\tilde{T}\in\mathcal{T}} |2|_F^{dim(W'_T)-\big(\frac{d+m}{2}\big)}|W(H,\tilde{T})|^{-1} \lim\limits_{s\to 0^+}\int_{\tilde{T}(F)/\theta} & c_\Theta(\tilde{t}) c_\Gamma(\tilde{t}) D^{\tilde{H}}(\tilde{t})^{1/2} \\
 & D^{\tilde{G}}(\tilde{t})^{1/2}\Delta(\tilde{t})^{s-1/2} d\tilde{t}
\end{aligned}\]

\noindent Si $\tilde{f}\in C_c^\infty(\tilde{G}(F))$ est très cuspidale, on pose $c_{\tilde{f}}=c_{\Theta_{\tilde{f}}}$ et $J_{geom}(\Theta,\tilde{f})=J_{geom}(\Theta,\Theta_{\tilde{f}})$.

\subsection{Le théorème}

Soient $\tilde{f}\in C_c^\infty(\tilde{G}(F))$ une fonction très cuspidale et $\Theta$ un quasi-caractère de $\tilde{H}(F)$. Pour $g\in G(F)$, on définit une fonction ${}^g \tilde{f}^\xi$ sur $\tilde{H}(F)$ par

$$\displaystyle {}^g \tilde{f}^\xi(\tilde{h})=\int_{U(F)} \tilde{f}(g^{-1}\tilde{h}ug) \xi(u) du$$

\noindent C'est une fonction localement constante à support compact. On peut donc poser

$$\displaystyle J(\Theta,\tilde{f},g)=\int_{\tilde{H}(F)} \Theta(\tilde{h}) {}^g\tilde{f}^\xi(\tilde{h}) d\tilde{h}$$

\noindent pour tout $g\in G(F)$. La fonction $g\mapsto J(\Theta,\tilde{f},g)$ est localement constante et est invariante à gauche par $H(F)U(F)$. Soit $N\geqslant 0$ un entier. Posons

$$\displaystyle J_N(\Theta,\tilde{f})=\int_{H(F)U(F)\backslash G(F)} \kappa_N(g) J(\Theta,\tilde{f},g) dg$$

\begin{theo}
L'expression $J_N(\Theta,\tilde{f})$ admet une limite lorsque $N$ tend vers l'infini et on a

$$\lim\limits_{N\to \infty} J_N(\Theta,\tilde{f})=J_{geom}(\Theta,\tilde{f})$$
\end{theo}

La démonstration de ce théorème occupera les quatre paragraphes suivants. Comme on l'a dit en 1.3, on utilise dans ces deux sections une normalisation différente des mesures. Le théorème 3.5.1 ne dépend que du choix de mesures de Haar sur les tores compacts $T_\theta(F)$, $\tilde{T}\in\mathcal{T}$, qui ont jusqu'à présent été normalisées par $mes(T_\theta(F))=1$. Ces mesures sont maintenant celles notées $\nu(T_\theta)dt$. Il faut donc dans la définition 3.4(1) rajouter un facteur $\nu(T_\theta)$ dans chaque terme de la somme.

\subsection{Descente à l'algèbre de Lie}

Par un procédé de partition de l'unité, on se ramène à la situation suivante: il existe $\tilde{x}\in \tilde{G}_{ss}(F)$ et $\omega\subset \mathfrak{g}_{\tilde{x}}(F)$ un bon voisinage de $0$ aussi petit que l'on veut de sorte que $Supp(\tilde{f})\subset (\tilde{x}exp(\omega))^G=\{g\tilde{x}exp(X)g^{-1};\; g\in G(F)\; X\in\omega\}$. \\

Supposons dans un premier temps que $\tilde{x}$ n'est conjugué à aucun élément de $\tilde{H}(F)$. Cela revient à dire que $(V''_{\tilde{x}},\tilde{x})$ ne contient pas de sous espace hermitien isomorphe à $(Z,\tilde{\zeta})$. Si tel est le cas et si $\omega$ est assez petit, il en va de même de $(V''_{\tilde{x}exp(X)},\tilde{x}exp(X))$ pour tout $X\in\omega$. On en déduit que $(\tilde{x}exp(\omega))^G\cap \tilde{H}(F)=\emptyset$. D'après une propriété des bons voisinages, l'ensemble $(\tilde{x}exp(\omega))^G$ est stable par l'opération qui consiste à passer à la partie semi-simple. On a donc aussi $(\tilde{x}exp(\omega))^G\cap \tilde{H}(F)U(F)=\emptyset$. D'après la définition, on a alors ${}^g\tilde{f}^\xi=0$ pour tout $g\in G(F)$ puis $J_N(\Theta,\tilde{f})=0$ pour tout $N$. D'autre part, le quasicaractère $\Theta_{\tilde{f}}$ est alors nul au voisinage de $\tilde{H}(F)$ et, par conséquent, $J_{geom}(\Theta,\tilde{f})=0$. Ainsi, le théorème 3.5.1 est trivialement vérifié dans ce cas. \\
 
On suppose donc maintenant que la classe de conjugaison de $\tilde{x}$ coupe $\tilde{H}(F)$. Quitte à conjuguer, on peut même supposer que $\tilde{x}\in\tilde{H}_{ss}(F)$. Posons $x={}^t (\tilde{x}^c)^{-1} \tilde{x}$ et on note $V''$ resp. $W''$ resp. $W'=V'$ le noyau de $x-1$ dans $V$ resp. le noyau de $x-1$ dans $W$ resp. l'image de $x-1$ dans $W$. On a alors $V=V''\oplus V'$, $W=W''\oplus W'$ et la restriction de $\tilde{x}$ à $V''$ et $W''$ est une forme hermitienne. On désignera par $G''$, $H''$ et $H'=G'$ les groupes $R_{E/F}GL(V'')$, $R_{E/F} GL(W'')$, $R_{E/F}GL(W')=R_{E/F}GL(V')$ respectivement. On a alors $G_{\tilde{x}}=G'_{\tilde{x}}G''_{\tilde{x}}$, $H_{\tilde{x}}=H'_{\tilde{x}}H''_{\tilde{x}}$ et les groupes $G''_{\tilde{x}}$, $H''_{\tilde{x}}$ sont des groupes unitaires. Soit $g\in G(F)$, on définit les fonctions $\Theta_{\tilde{x},\omega}$ et ${}^g \tilde{f}_{\tilde{x},\omega}$ sur $\mathfrak{h}_{\tilde{x}}(F)$ et $\mathfrak{g}_{\tilde{x}}(F)$ respectivement par

$$\Theta_{\tilde{x},\omega}(X)=\left\{
    \begin{array}{ll}
        \Theta(\tilde{x}exp(X)) & \mbox{si } X\in \omega\cap\mathfrak{h}_{\tilde{x}}(F) \\
        0 & \mbox{sinon.}
    \end{array}
\right. $$

\noindent et

$${}^g \tilde{f}_{\tilde{x},\omega}(X)=\left\{
    \begin{array}{ll}
        \tilde{f}(g^{-1}\tilde{x}exp(X)g) & \mbox{si } X\in\omega\\
        0 & \mbox{sinon.}
    \end{array}
\right. $$

Pour tout $X\in\mathfrak{h}_{\tilde{x}}(F)$ et pour tout $g\in G(F)$, posons

$$\displaystyle {}^g \tilde{f}_{\tilde{x},\omega}^\xi(X)=\int_{\mathfrak{u}_{\tilde{x}}(F)} {}^g \tilde{f}_{\tilde{x},\omega}(X+N)\xi(N) dN$$

$$\displaystyle J_{\tilde{x},\omega}(\Theta,\tilde{f},g)=\int_{\mathfrak{h}_{\tilde{x}}(F)} \Theta_{\tilde{x},\omega}(X) {}^g \tilde{f}_{\tilde{x},\omega}^\xi(X) dX$$

Pour tout entier $N\geqslant 0$, définissons

$$\displaystyle J_{\tilde{x},\omega,N}(\Theta,\tilde{f})= \int_{H_{\tilde{x}}(F)U_{\tilde{x}}(F)\backslash G(F)} \kappa_N(g) J_{\tilde{x},\omega}(\Theta,\tilde{f},g) dg$$

Notons $\underline{\mathcal{T}}'$ l'ensemble des tores maximaux anisotropes de $H'_{\tilde{x}}$. La paire d'espaces hermitiens $(V'',W'')$ vérifient les hypothèses de la section 4 de [B], on peut donc lui associer un ensemble de tores $\underline{\mathcal{T}}''$ comme dans cette référence. Notons $\underline{\mathcal{T}}_{\tilde{x}}$ l'ensemble des éléments de $\underline{\mathcal{T}}$ qui contiennent $\tilde{x}$. On vérifie que l'application $\tilde{T}\mapsto T_\theta$ est une bijection entre $\underline{\mathcal{T}}_{\tilde{x}}$ et $\underline{\mathcal{T}}'\times \underline{\mathcal{T}}''$. Pour $\tilde{T}\in\underline{\mathcal{T}}$, on a défini deux fonctions $c_\Theta$ et $c_{\tilde{f}}$ sur $\tilde{T}(F)/\theta$. On définit les fonctions $c_{\Theta,\tilde{x},\omega}$ et $c_{\tilde{f},\tilde{x},\omega}$ sur $\mathfrak{t}_\theta(F)$ comme étant nulles hors de $\omega$ et dont les valeurs sur $X\in\omega\cap \mathfrak{t}_\theta(F)$ sont

\begin{center}
$c_{\Theta,\tilde{x},\omega}(X)=c_\Theta(\tilde{x} exp(X))$ et $c_{\tilde{f},\tilde{x},\omega}(X)=c_{\tilde{f}}(\tilde{x}exp(X))$
\end{center}

\noindent On note $\Delta''$ la fonction définie sur $\mathfrak{h}_{\tilde{x},ss}(F)$ par

$$\Delta''(X)=|N(det X''_{|W''/W''(X)})|_F$$

\noindent où $X''$ désigne la restriction de $X$ à $W''$ et $W''(X)$ le noyau de $X$ dans $W''$. Soit $\mathcal{T}_{\tilde{x}}$ un ensemble de représentants des classes de conjugaison de $\underline{\mathcal{T}}_{\tilde{x}}$ par $H_{\tilde{x}}(F)$. Soit $I\in \underline{\mathcal{T}}'\times \underline{\mathcal{T}}''$ et notons $\tilde{T}_I\in \underline{\mathcal{T}}_{\tilde{x}}$ le tore tordu lui correspondant. D'après le paragraphe 3.2, il est associé à $\tilde{T}_I$ une décomposition $W=W'_{T_I}\oplus W''_{T_I}$. Notons $H''_{T_I}=R_{E/F} GL(W''_{T_I})$ et $Norm_{H_{\tilde{x}}}(I)$ le normalisateur de $I$ dans $H_{\tilde{x}}$. On a alors une inclusion naturelle $H''_{T_I,\tilde{x}}I\subset Norm_{H_{\tilde{x}}}(I)$. On pose

$$W(H_{\tilde{x}},I)=Norm_{H_{\tilde{x}}}(I)(F)/H''_{T_I,\tilde{x}}(F)I(F)$$

\noindent C'est un groupe fini. Posons

\[\begin{aligned}
\displaystyle J_{geom,\tilde{x},\omega}(\Theta,\tilde{f})=\sum_{I\in \mathcal{T}_{\tilde{x}}} |W(H_{\tilde{x}},I)|^{-1} \nu(I) \lim\limits_{s\to 0^+} \int_{\mathfrak{i}(F)} &  c_{\Theta,\tilde{x},\omega}(X) c_{\tilde{f},\tilde{x},\omega}(X) D^{H_{\tilde{x}}}(X)^{1/2} \\
 & D^{G_{\tilde{x}}}(X)^{1/2}\Delta''(X)^{s-1/2} dX
\end{aligned}\]

\noindent On voit comme en 3.4 que cette expression est bien définie. Soit

$$\displaystyle C(\tilde{x})=|2|_F^{\frac{m-d}{2}-dim(W'')} D^{\tilde{H}}(\tilde{x})^{1/2} D^{\tilde{G}}(\tilde{x})^{1/2} \Delta(\tilde{x})^{-1/2}$$

\begin{prop}
On a les égalités

\begin{center}
$J_N(\Theta,\tilde{f})=C(\tilde{x})J_{N,\tilde{x},\omega}(\Theta,\tilde{f})$ et $J_{geom}(\Theta,\tilde{f})=C(\tilde{x})J_{geom,\tilde{x},\omega}(\Theta,\tilde{f})$
\end{center}
\end{prop}

\ul{Preuve}:

Il suffit de reprendre la preuve du lemme 3.4 de [W3]. Rappelons pour la commodité du lecteur comment l'on procède. \\

Soit $g\in G(F)$. Par la formule d'intégration de Weyl, on a:

$$\mbox{(1)}\;\;\; J(\Theta,\tilde{f},g)=\displaystyle\sum_{\tilde{T}\in\mathcal{T}(\tilde{H})} |W(H,\tilde{T})|^{-1} \displaystyle\int_{\tilde{T}(F)/\theta} \Theta(\tilde{t}) D^{\tilde{H}}(\tilde{t}) \int_{T_\theta(F)\backslash H(F)}{}^{hg}\tilde{f}^\xi(\tilde{t}) dh d\tilde{t}$$

Pour $\tilde{T}$ et $\tilde{T}'$ deux sous-tores maximaux tordus de $\tilde{H}$ on pose $W(\tilde{T},\tilde{T}')=\{h\in H(F): h\tilde{T}h^{-1}=\tilde{T}'\}/T(F)$. Pour tout sous-tore maximal $I$ de $H_{\tilde{x}}$, on note $T_I$ son commutant dans $H$ et $\tilde{T}_I=T_I\tilde{x}$ qui est un sous-tore maximal tordu de $\tilde{H}$. On a alors:

\vspace{2mm}

(2) Soient $\tilde{T}\in \mathcal{T}(\tilde{H})$ et $\tilde{t}\in \tilde{T}(F)\cap \tilde{H}_{reg}(F)$ tel que $\displaystyle\int_{T_\theta(F)\backslash H(F)}{}^{hg}\tilde{f}^\xi(\tilde{t}) dh\neq 0$ alors

$$\displaystyle \tilde{t}\in\bigcup_{I\in \mathcal{T}(H_{\tilde{x}})}\bigcup_{w\in W(\tilde{T}_I,\tilde{T})} w(\tilde{x}exp(\mathfrak{i}(F)\cap \omega))w^{-1}$$

\vspace{2mm}

 En effet, il existe alors $u\in U(F)$ tel que $\tilde{t}u\in(\tilde{x}exp(\omega))^G$. La partie semisimple de $\tilde{t}u$ étant conjuguée à $\tilde{t}$ on a aussi $\tilde{t}\in (\tilde{x}exp(\omega))^G$. Soient donc $X\in \omega$ et $y\in G(F)$ tels que $y\tilde{t}y^{-1}=\tilde{x}exp(X)$. Par restriction, $y$ réalise alors un isomorphisme entre les espaces hermitiens $(V''_{\tilde{t}},\tilde{t})$ et $(V''_{\tilde{x}exp(X)},\tilde{x} exp(X))$. Le premier de ces espaces contient comme sous-espace hermitien $(Z,\tilde{\zeta})$ et le deuxième est inclus dans $(V'',\tilde{x})$ qui contient aussi $(Z,\tilde{\zeta})$. D'après le théorème de Witt, quitte à multiplier $y$ par un élément de $G''_{\tilde{x}}(F)$ et à conjuguer $X$ par ce même élément, on peut supposer que $y$ agit trivialement sur $Z$. On a alors l'égalité $y{}^t(\tilde{t}^c)^{-1}\tilde{t}y^{-1}=x exp(2X)$. On en déduit que $exp(2X)$ agit comme l'identité sur $Z$ donc la restriction de $X$ à $Z$ est nulle. Un élément de $\mathfrak{g}_{\tilde{x}}(F)$ dont la restriction à $Z$ est nulle appartient à $\mathfrak{h}_{\tilde{x}}(F)$, par conséquent $X\in \mathfrak{h}_{\tilde{x}}(F)$. L'élément $y$ envoie l'orthogonal de $Z$ dans $(V''_{\tilde{t}},\tilde{t})$, qui est $W''_{\tilde{t}}$, sur l'orthogonal de $Z$ dans $V''_{\tilde{x}exp(X)}$, qui est $W''_{\tilde{x}exp(X)}$. Il envoie aussi $W'_{\tilde{t}}$ sur $W'_{\tilde{x}exp(X)}$. Puisque $W=W'_{\tilde{t}}\oplus W''_{\tilde{t}}=W'_{\tilde{x}exp(X)}\oplus W''_{\tilde{x}exp(X)}$, on en déduit que $y$ conserve $W$ donc $y\in H(F)$. Quitte à conjuguer $X$, on peut supposer qu'il existe $I\in \mathcal{T}(H_{\tilde{x}})$ tel que $X\in \mathfrak{i}(F)\cap\omega$. Puisque $\tilde{t}$ est supposé régulier, la conjugaison par $y^{-1}$ induit alors un élément $w$ de $W(\tilde{T}_I,\tilde{T})$ d'où le résultat.

\vspace{2mm}

(3) Pour $\tilde{T}\in \mathcal{T}(\tilde{H})$, $I_1,I_2\in \mathcal{T}(H_{\tilde{x}})$ et $w_1\in W(\tilde{T}_{I_1},\tilde{T})$, $w_2\in W(\tilde{T}_{I_2},\tilde{T})$ les deux ensembles $w_1(\tilde{x}exp(\mathfrak{i}_1(F)\cap\omega))w_1^{-1}$ et $w_2(\tilde{x}exp(\mathfrak{i}_2(F)\cap\omega))w_2^{-1}$ sont disjoints ou égaux et s'ils sont égaux on a $I_1=I_2$.

\vspace{2mm}

Soient $y_1,y_2\in H(F)$ qui relèvent $w_1$ et $w_2$ on pose $y=y_2^{-1}y_1$. Si les deux ensembles en question ne sont pas disjoints, on a $y(\tilde{x}exp(\omega))y^{-1}\cap \tilde{x}exp(\omega)\neq \emptyset$. Donc par une propriété des bons voisinages, $y\in Z_H(\tilde{x})(F)=H_{\tilde{x}}(F)$. Alors la conjugaison par $y$ envoie $\tilde{T}_{I_1}$ sur $\tilde{T}_{I_2}$, donc $I_1$ sur $I_2$, et laisse stable $\omega$. Ce qui entraîne aussi $I_1=I_2$.

\vspace{2mm}

(4) Soient $\tilde{T}\in\mathcal{T}(\tilde{H})$, $I\in\mathcal{T}(H_{\tilde{x}})$ et $w_1\in W(\tilde{T}_{I},\tilde{T})$. Le nombre d'éléments $w_2\in W(\tilde{T}_{I},\tilde{T})$ tels que $w_2(\tilde{x}exp(\mathfrak{i}(F)\cap\omega))w_2^{-1}=w_1(\tilde{x}exp(\mathfrak{i}(F)\cap\omega))w_1^{-1}$ est $|W(H_{\tilde{x}},I)|$.

\vspace{2mm}

En effet d'après ce qu'on vient de voir et le fait que $T_I\cap H_{\tilde{x}}=I$, l'ensemble de ces éléments est en bijection avec $\{y\in H_{\tilde{x}}(F): yIy^{-1}=I\}/I(F)$ qui est précisément $W(H_{\tilde{x}},I)$.

\vspace{2mm}

(5) Soit $I\in\mathcal{T}(H_{\tilde{x}})$, alors il existe un unique tore $\tilde{T}\in \mathcal{T}(\tilde{H})$ tel que $W(\tilde{T}_I,\tilde{T})\neq \emptyset$ et on a alors $|W(\tilde{T}_I,\tilde{T})|=|W(H,\tilde{T})|$.

\vspace{2mm}

L'existence et l'unicité de $\tilde{T}$ sont évidentes puisque $\mathcal{T}(\tilde{H})$ est précisément un ensemble de représentants des classes de conjugaison de sous-tores maximaux tordus de $\tilde{H}$. Soit donc $\tilde{T}$ l'unique élément de $\mathcal{T}(\tilde{H})$ tel que $W(\tilde{T}_I,\tilde{T})\neq\emptyset$ et fixons $w_1\in W(\tilde{T}_I,\tilde{T})$. On a alors une bijection entre $W(\tilde{T}_I,\tilde{T})$ et $W(H,\tilde{T})$ donnée par $w\mapsto ww_1^{-1}$.

\vspace{2mm}

Rappelons que d'après nos choix de mesures, pour $I\in \mathcal{T}(H_{\tilde{x}})$, l'application

$$\mathfrak{i}(F)\cap \omega \to \tilde{T}_I(F)/\theta$$
$$X\mapsto \tilde{x} exp(X)$$

\noindent préserve les mesures. Les points 2 à 5 ci-dessus permettent alors de transformer l'expression $(1)$ en l'égalité suivante

\[\begin{aligned}
\displaystyle J(\Theta,\tilde{f},g) =\sum_{I\in\mathcal{T}(H_{\tilde{x}})} |W(H_{\tilde{x}},I)|^{-1} \int_{\mathfrak{i}(F)}\int_{I(F)\backslash H(F)} & {}^{hg}\tilde{f}^\xi(\tilde{x} exp(X))dh \\
 & D^{\tilde{H}}(\tilde{x}exp(X)) \Theta_{\tilde{x},\omega}(X) dX
\end{aligned}\]

\noindent On a $D^{\tilde{H}}(\tilde{x} exp(X))=D^{\tilde{H}}(\tilde{x}) D^{H_{\tilde{x}}}(X)$. Par la formule d'intégration de Weyl appliquée au groupe $H_{\tilde{x}}$, on a

$$\displaystyle J(\Theta,\tilde{f},g)=D^{\tilde{H}}(\tilde{x}) \int_{H_{\tilde{x}}(F)\backslash H(F)} \int_{\mathfrak{h}_{\tilde{x}}(F)} \Theta_{\tilde{x},\omega}(X) {}^{hg}\tilde{f}^\xi (\tilde{x}exp(X)) dX dh$$

Notons $\mathfrak{u}^{\tilde{x}}$ l'image de $\theta_{\tilde{x}}-1$ dans $\mathfrak{u}$. On a alors $\mathfrak{u}=\mathfrak{u}_{\tilde{x}}\oplus \mathfrak{u}^{\tilde{x}}$. On a fixé une mesure sur $\mathfrak{u}(F)$ et on suppose que cette mesure est produit d'une mesure sur $\mathfrak{u}_{\tilde{x}}(F)$ et d'une mesure sur $\mathfrak{u}^{\tilde{x}}(F)$. Soit $U^{\tilde{x}}(F)=exp(\mathfrak{u}^{\tilde{x}}(F))$ que l'on munit de la mesure image par l'exponentielle (signalons qu'en général, il ne s'agit pas d'un sous-groupe de $U(F)$). Pour $h\in H(F)$ et $X\in \omega$, on a

$$\displaystyle {}^{hg}\tilde{f}^\xi(\tilde{x} exp(X))=\int_{U^{\tilde{x}}(F)}\int_{U_{\tilde{x}}(F)}{}^{hg}\tilde{f}(\tilde{x}exp(X)uv)\xi(uv) dudv$$

\noindent Pour $u\in U_{\tilde{x}}(F)$, l'application $U^{\tilde{x}}\to U_{\tilde{x}}\backslash U$, $v\mapsto \theta_{\tilde{x}exp(X)u}^{-1}(v^{-1}) v$ est un isomorphisme de variétés algébriques et le morphisme sur les $F$-points a un jacobien constant de valeur $|det (1-\theta_{\tilde{x}exp(X)u}^{-1})_{\mathfrak{u}/\mathfrak{u}_{\tilde{x}}}|_F$. Par une propriété des bons voisinages, ce jacobien est aussi la valeur absolue du déterminant de $1-\theta_{\tilde{x}}^{-1}$ agissant sur $\mathfrak{u}/\mathfrak{u}_{\tilde{x}}$. Notons $d(\tilde{x})$ ce jacobien. Remarquons que l'image de l'application précédente est incluse dans $Ker \; \xi$. On a alors

$$\displaystyle {}^{hg}\tilde{f}^\xi(\tilde{x}exp(X))=d(\tilde{x})\int_{U_{\tilde{x}}(F)\backslash U(F)} \int_{U_{\tilde{x}}(F)} {}^{vhg}\tilde{f}(\tilde{x}exp(X)u)\xi(u)dudv$$

L'application $\mathfrak{u}_{\tilde{x}}(F)\to U_{\tilde{x}}(F)$, $N\mapsto exp(-X)exp(X+N)$ est un isomorphisme de jacobien constant égale à $1$. On en déduit l'égalité

$$\displaystyle {}^{hg}\tilde{f}^\xi(\tilde{x}exp(X))=d(\tilde{x})\int_{U_{\tilde{x}}(F)\backslash U(F)} {}^{vhg}\tilde{f}_{\tilde{x},\omega}^\xi(X) dv$$

\noindent D'où

$$\displaystyle J(\Theta,\tilde{f},g)=D^{\tilde{H}}(\tilde{x})d(\tilde{x})\int_{U_{\tilde{x}}(F)H_{\tilde{x}}(F)\backslash U(F)H(F)}  J_{\tilde{x},\omega}(\Theta,\tilde{f},vhg) dhdv$$

\noindent Puis $J_N(\Theta,\tilde{f})=D^{\tilde{H}}(\tilde{x})d(\tilde{x}) J_{\tilde{x},\omega,N}(\Theta,\tilde{f})$. Il ne nous reste plus qu'à déterminer la valeur de $d(\tilde{x})$. \\

Définissons les sous-espaces suivants de $\mathfrak{u}$

\begin{center}
$\mathfrak{u}_1=\{N\in\mathfrak{u}; \; N(V')=0 \;\mbox{et} \; N(V'')\subset V''\}$ et $\mathfrak{u}_2=\{N\in\mathfrak{u}; \; N(V'')\subset V'\}$
\end{center}
Ces deux sous-espaces sont stables par $\theta_{\tilde{x}}$ et on a $\mathfrak{u}=\mathfrak{u}_1\oplus\mathfrak{u}_2$. On note $\mathfrak{u}_{1,\tilde{x}}$ et $\mathfrak{u}_{2,\tilde{x}}$ les noyaux de $\theta_{\tilde{x}}-1$ dans $\mathfrak{u}_1$ et $\mathfrak{u}_2$ respectivement. \\

\begin{itemize}
\item Calcul de $|det(1-\theta_{\tilde{x}}^{-1})_{|\mathfrak{u}_1/\mathfrak{u}_{1,\tilde{x}}}|_F$. L'espace $\mathfrak{u}_1$ est inclus dans $\mathfrak{gl}(V'')$. Sur cet espace $\theta_{\tilde{x}}$ n'est autre que l'application qui à $X$ associe l'opposé de l'adjoint de $X$ pour la forme hermitienne $\tilde{x}_{|V''}$. On a par conséquent $\theta_{\tilde{x}}^2=Id$. Notons $\mathfrak{u}_1^{\tilde{x}}$ le sous-espace de $\mathfrak{u}_1$ sur lequel $\theta_{\tilde{x}}=-Id$. On a alors $|det(1-\theta_{\tilde{x}}^{-1})_{|\mathfrak{u}_1/\mathfrak{u}_{1,\tilde{x}}}|_F=|2|_F^{dim(\mathfrak{u}_1^{\tilde{x}})}$. Soit $\eta\in E$ un élément non nul de trace nulle. L'application $N\mapsto \eta N$ est alors un isomorphisme entre $\mathfrak{u}_{1,\tilde{x}}$ et $\mathfrak{u}_1^{\tilde{x}}$. On en déduit que $dim(\mathfrak{u}_1^{\tilde{x}})=dim(\mathfrak{u}_1)/2$. Or, un calcul indolore nous donne que $dim(\mathfrak{u}_1)=4r^2+2r+4rdim(W'')$. On en déduit que

$$\mbox{(6)}\;\;\; |det(1-\theta_{\tilde{x}}^{-1})_{|\mathfrak{u}_1/\mathfrak{u}_{1,\tilde{x}}}|_F=|2|_F^{2r^2+r+2rdim(W'')}$$

\item Calcul de $|det(1-\theta_{\tilde{x}}^{-1})_{|\mathfrak{u}_2/\mathfrak{u}_{2,\tilde{x}}}|_F$. Posons $\mathfrak{u}_3=\{N\in\mathfrak{u}_2;\; N(V'')=0\}$ et $\mathfrak{u}_4=\{N\in\mathfrak{u}_2;\; N(V')=0\}$. On a alors $\mathfrak{u}_2=\mathfrak{u}_3\oplus\mathfrak{u}_4$, $\theta_{\tilde{x}}(\mathfrak{u}_3)=\mathfrak{u}_4$ et $\theta_{\tilde{x}}(\mathfrak{u}_4)=\mathfrak{u}_3$. On en déduit que 

$$det(1-\theta_{\tilde{x}}^{-1})_{|\mathfrak{u}_2/\mathfrak{u}_{2,\tilde{x}}}=det(1-\theta_{\tilde{x}}^{-2})_{|\mathfrak{u}_3/\mathfrak{u}_{3,\tilde{x}}}$$

\noindent où $\mathfrak{u}_{3,\tilde{x}}$ désigne le noyau de $1-\theta_{\tilde{x}}$ dans $\mathfrak{u}_3$. Notons $Z_+$ le sous-espace de $V$ engendré par les $z_i$, $i=1,\ldots,r$. On a alors un isomophisme naturel $\mathfrak{u}_3\simeq Hom_E(V',Z_+)$. Via cet isomorphisme, l'action de $\theta_{\tilde{x}}^2$ sur $\mathfrak{u}_3$ est donné par $X\mapsto X\circ x^{-1}$. On a donc un isomorphisme de $F[\theta_{\tilde{x}}^2]$-modules $\displaystyle\mathfrak{u}_3\simeq \bigoplus_{i=1}^r V'$ où l'action de $\theta_{\tilde{x}}^2$ sur le deuxième espace est donné par l'action de $x^{-1}$ sur chaque facteur. On en déduit que

$$\mbox{(7)}\;\;\; |det(1-\theta_{\tilde{x}}^{-1})_{|\mathfrak{u}_2/\mathfrak{u}_{2,\tilde{x}}}|_F=|N(det (1-x)_{|V'})|_F^r=\Delta(\tilde{x})^r$$
\end{itemize}

\noindent Des relations (6) et (7), on tire l'égalité $d(\tilde{x})=|2|_F^{2r^2+r+2rdim(W'')}\Delta(\tilde{x})^r$. On vérifie facilement que $D^{\tilde{H}}(\tilde{x})d(\tilde{x})=C(\tilde{x})$, cf lemme 2.1.1, d'où la première assertion de l'énoncé.

\vspace{2mm}

Pour $I\in \mathcal{T}_{\tilde{x}}$, notons $\tilde{T}_I$ l'unique élément de $\underline{\mathcal{T}}_{\tilde{x}}$ tel que $T_{I,\theta}=I$. Pour $\tilde{T}_1,\tilde{T}_2\in\underline{\mathcal{T}}$, posons

$$W(\tilde{T}_1,\tilde{T}_2)=\{h\in H(F);\; h\tilde{T}_1h^{-1}=\tilde{T}_2\}/(T_1(F)\times U(\tilde{\zeta}_{H,T_1})(F))$$

Les propriétés suivantes se montrent exactement de la même manière que les propriétés (2) à (5): \\

 (8) Soient $\tilde{T}\in\mathcal{T}$ et $\tilde{t}\in\tilde{T}_\sharp(F)$ tel que $c_{\tilde{f}}(\tilde{t})\neq 0$ alors
 
 $$\displaystyle\tilde{t}\in \bigcup_{I\in \mathcal{T}_{\tilde{x}}}\bigcup_{w\in W(\tilde{T}_I,\tilde{T})} w(\tilde{x}exp(\mathfrak{i}(F)\cap\omega))w^{-1}$$
 
 (9) Soient $\tilde{T}\in\mathcal{T}$, $I_i\in \mathcal{T}_{\tilde{x}}$ et $w_i\in W(\tilde{T}_{I_i},\tilde{T})$ pour $i=1,2$. Alors les ensembles $w_i(\tilde{x}exp(\mathfrak{i}_i(F)\cap\omega))w_i^{-1}$, $i=1,2$, sont disjoints ou confondus et s'ils sont confondus on a $I_1=I_2$. \\
 
 (10) Soient $\tilde{T}\in\mathcal{T}$, $I\in\mathcal{T}_{\tilde{x}}$ et $w_1\in W(\tilde{T}_I,\tilde{T})$. Alors le nombre d'éléments $w_2\in W(\tilde{T}_I,\tilde{T})$ tels que
 
 $$w_1(\tilde{x}exp(\mathfrak{i}(F)\cap \omega))w_1^{-1}=w_2(\tilde{x}exp(\mathfrak{i}(F)\cap \omega))w_2^{-1}$$
 
\noindent est égal à $|W(H_{\tilde{x}},I)|$. \\
 
 (11) Pour tout $I\in \mathcal{T}_{\tilde{x}}$, il existe un unique $\tilde{T}\in \mathcal{T}$ tel que $W(\tilde{T}_I,\tilde{T})\neq \emptyset$ et pour ce $\tilde{T}$, on a $|W(\tilde{T}_I,\tilde{T})|=|W(H,\tilde{T})|$. \\
 
Les points (8) à (11) permettent de transformer l'expression initiale de $J_{geom}(\Theta,\tilde{f})$ et d'obtenir l'égalité

\[\begin{aligned}
\mbox{(12)}\;\;\; \displaystyle J_{geom}(\Theta,\tilde{f}) & =D^{\tilde{H}}(\tilde{x})^{1/2} D^{\tilde{G}}(\tilde{x})^{1/2} \sum_{I\in \mathcal{T}_{\tilde{x}}}|2|_F^{dim(W'_{T_I})-\big(\frac{d+m}{2}\big)} |W(H_{\tilde{x}},I)|^{-1} \nu(I) \\
 & \lim\limits_{s\to 0^+}\int_{\mathfrak{i}(F)} c_{\Theta,\tilde{x},\omega}(X) c_{\tilde{f},\tilde{x},\omega}(X) D^{H_{\tilde{x}}}(X)^{1/2}D^{G_{\tilde{x}}}(X)^{1/2} \Delta(\tilde{x} exp(X))^{s-1/2} dX
\end{aligned}\]

Soit $I\in\mathcal{T}_{\tilde{x}}$. Pour tout $X\in\omega$, on a ${}^t[(\tilde{x}exp(X))^c]^{-1}\tilde{x}exp(X)=xexp(2X)$, d'où

$$\Delta(\tilde{x}exp(X))=|2|_F^{2dim(W'')-2dim(W''_{T_I})}\Delta''(X)\Delta(\tilde{x})$$

\noindent pour presque tout $X\in\omega\cap \mathfrak{i}(F)$. Puisque $dim(W'_{T_I})+dim(W''_{T_I})=m$ et $\lim\limits_{s\to 0}\Delta(\tilde{x})^s=1$, on déduit facilement de (12) la deuxième égalité de l'énoncé $\blacksquare$

\subsection{Première transformation de $J_{\tilde{x},\omega,N}(\Theta,\tilde{f})$}

D'après la proposition précédente, on est ramené à montrer l'égalité

$$\lim\limits_{N\to\infty} J_{\tilde{x},\omega,N}(\Theta,\tilde{f})=J_{\tilde{x},\omega,geom}(\Theta,\tilde{f})$$

On a des décompositions en produits $G_{\tilde{x}}=G''_{\tilde{x}}G'_{\tilde{x}}$ et $H_{\tilde{x}}=H''_{\tilde{x}}H'_{\tilde{x}}$. Les groupes $G''_{\tilde{x}}$ et $H''_{\tilde{x}}$ sont les groupes unitaires des espaces hermitiens $(V'',\tilde{x})$ et $(W'',\tilde{x})$ respectivement. Pour $X\in\mathfrak{h}_{\tilde{x}}(F)=\mathfrak{h}'_{\tilde{x}}(F)\oplus\mathfrak{h}''_{\tilde{x}}(F)$ (resp. $X\in\mathfrak{g}_{\tilde{x}}(F)=\mathfrak{g}'_{\tilde{x}}(F)\oplus\mathfrak{g}''_{\tilde{x}}(F)$) on notera $X=X'+X''$ l'unique écriture de $X$ avec $X'\in\mathfrak{h}'_{\tilde{x}}(F)$ et $X''\in\mathfrak{h}''_{\tilde{x}}(F)$ (resp. $X'\in\mathfrak{g}'_{\tilde{x}}(F)$ et $X''\in\mathfrak{g}''_{\tilde{x}}(F)$). Quitte à rétrécir $\omega$, on peut supposer que $\Theta_{\tilde{x},\omega}$ est combinaison linéaire sur $\omega$ de transformées de Fourier d'intégrales orbitales nilpotentes. D'après "la conjecture de Howe", $\Theta_{\tilde{x},\omega}$ est alors aussi combinaison linéaire de fonctions $\hat{j}^{H_{\tilde{x}}}(Y,.)$ pour des éléments $Y\in\mathfrak{h}_{\tilde{x},reg}(F)$ génériques. Par linéarité de l'égalité que l'on cherche à montrer, on peut donc supposer qu'il existe $S\in\mathfrak{h}''_{\tilde{x},reg}(F)$ de noyau nul dans $W''$ et un quasi-caractère $\hat{j}_S$ de $\mathfrak{h}'_{\tilde{x}}(F)$ tels que

$$\Theta_{\tilde{x},\omega}(X)=\hat{j}_S(X')\hat{j}^{H''_{\tilde{x}}}(S,X'')$$

\noindent pour presque tout $X\in \omega\cap\mathfrak{h}_{\tilde{x}}(F)$. On a alors pour tout $g\in G(F)$,

$$\displaystyle J_{\tilde{x},\omega}(\Theta,\tilde{f},g)= \int_{\mathfrak{h}'_{\tilde{x}}(F)\times \mathfrak{h}''_{\tilde{x}}(F)} \hat{j}_S(X')\hat{j}^{H''_{\tilde{x}}}(S,X'') {}^g\tilde{f}_{\tilde{x},\omega}^\xi(X) dX''dX'$$

Appliquons la formule d'intégration de Weyl à l'intégrale sur $\mathfrak{h}'_{\tilde{x}}(F)$. On obtient alors pour $N\geqslant 0$,

\[\begin{aligned}
\displaystyle J_{\tilde{x},\omega,N}(\Theta,\tilde{f}) & =\sum_{I'\in \mathcal{T}(H'_{\tilde{x}})} |W(H'_{\tilde{x}},I')|^{-1} \int_{\mathfrak{i}'(F)} \hat{j}_S(X') D^{H'_{\tilde{x}}}(X') \\
 & \int_{I'(F)H''_{\tilde{x}}(F)U_{\tilde{x}}(F)\backslash G(F)} \int_{\mathfrak{h}''_{\tilde{x}}(F)} \hat{j}^{H''_{\tilde{x}}}(S,X'') {}^g\tilde{f}_{\tilde{x},\omega}^\xi(X'+X'')dX'' \kappa_N(g) dg dX'
\end{aligned}\]

Le terme intérieur peut se réécrire

$$\displaystyle \int_{I'(F)G''_{\tilde{x}}(F)\backslash G(F)} \int_{H''_{\tilde{x}}(F)U_{\tilde{x}}(F)\backslash G''_{\tilde{x}}(F)} \int_{\mathfrak{h}''_{\tilde{x}}(F)} \hat{j}^{H''_{\tilde{x}}}(S,X'') {}^{g''g}\tilde{f}_{\tilde{x},\omega}^\xi(X'+X'') dX'' \kappa_N(g''g) dg'' dg$$

D'après 7.4 [B], on peut exprimer différemment les deux intégrales intérieures de cette expression. Pour retouver les notations de cette référence, on pose $\varphi(X'')={}^g \tilde{f}_{\tilde{x},\omega}(X'+X'')$, $\kappa''(g'')=\kappa_N(g''g)$, $\theta''=\hat{j}^{G''_{\tilde{x}}}(S,.)$, $v_i=z_i$ pour $i=0,\ldots,r$ et $v_{-i}=z_{-i}/((-1)^i 2\nu)$ pour $i=1,\ldots,r$. On doit aussi fixer des constantes $\xi_i\in F^\times$ vérifiant

$$\displaystyle \sum_{i=-r}^{r-1} u_{i+1,i}=\sum_{i=0}^{r-1} \xi_i <\tilde{x}v_{-i-1}, uv_i>$$

\noindent pour tout $u\in U_{\tilde{x}}(F)$. Dans la section 7 de [B], il est défini des éléments $\Xi$ et $\Sigma$. L'élément $\Xi\in\mathfrak{g}_{\tilde{x}}''(F)$ est défini ainsi: il annule $W''$ et vérifie $\Xi(v_{i+1})=\xi_iv_i$ pour $i=0,\ldots,r-1$, $\Xi v_0=-2\nu\xi_0 v_{-1}$ et $\Xi v_{-i}=-\xi_i v_{-i-1}$ pour $i=1,\ldots,r-1$. Le sous-espace $\Sigma$ de $\mathfrak{g}_{\tilde{x}}''$ est défini comme l'orthogonal de $\mathfrak{h}_{\tilde{x}}''\oplus\mathfrak{u}_{\tilde{x}}$. Soit $I''\in \mathcal{T}(G''_{\tilde{x}})$. On note $\mathfrak{i}''(F)^S$ l'ensemble des éléments de $\mathfrak{i}''(F)$ dont la classe de conjugaison par $G''_{\tilde{x}}(F)$ coupe $\Xi+S+\Sigma$. Comme en 7.3 [B], on peut fixer:

\vspace{3mm}

\begin{itemize}
\item Deux fonctions localement analytiques $X''\in\mathfrak{i}''(F)^S\mapsto \gamma_{X''}\in G''_{\tilde{x}}(F)$ et $X''\in\mathfrak{i}''(F)^S\mapsto X''_\Sigma\in \Xi+S+\Sigma$ vérifiant $X''_\Sigma=\gamma_{X''}^{-1}X''\gamma_{X''}$ pour tout $X''\in\mathfrak{i}''(F)^S$;
\item Une fonction polynomiale non nulle $Q_S$ sur $\mathfrak{i}''(F)$
\end{itemize}

\vspace{3mm}

\noindent telles que pour tout compact $\omega_{I''}\subset \mathfrak{i}''(F)$, il existe $\alpha\in F^\times$ et $c>0$ vérifiant

\vspace{3mm}

\begin{itemize}
\item $\alpha Q_S(X'')R''\subset R_{X''_\Sigma,v_r}\subset \alpha^{-1} R''$ où on a posé $R''=R\cap V''$ et $R_{X''_\Sigma,v_r}=\bigoplus_{i=0}^{d''-1} {X''_\Sigma}^i v_r$;
\item $\sigma(\gamma_{X''})\leqslant c(1+log \; |D^{G''_{\tilde{x}}}(X'')|)$
\end{itemize}

\vspace{3mm}

\noindent pour tout $X''\in\mathfrak{i}''(F)^S\cap \omega_{I''}$. On peut maintenant appliquer 7.4 [B] à nos deux intégrales intérieures. On obtient l'égalité

\[\begin{aligned}
\mbox{(1)}\;\;\; \displaystyle J_{\tilde{x},\omega,N}(\Theta,\tilde{f})=\sum_{I\in\mathcal{T}(G_{\tilde{x}})} \nu(A_{I''})^{-1} |W(G_{\tilde{x}},I)|^{-1} \int_{\mathfrak{i}'(F)\times \mathfrak{i}''(F)^S} \hat{j}_S(X') D^{G'_{\tilde{x}}}(X')D^{G''_{\tilde{x}}}(X'')^{1/2} \\
\int_{I'(F)A_{I''}(F)\backslash G(F)} {}^g \tilde{f}_{\tilde{x},\omega}^\sharp(X) \kappa_{N,X''}(g) dg dX''dX'
\end{aligned}\]

\noindent où on a noté ${}^g \tilde{f}_{\tilde{x},\omega}^\sharp$ la transformée de Fourier de ${}^g \tilde{f}_{\tilde{x},\omega}$ par rapport à $X''\in \mathfrak{g}''_{\tilde{x}}(F)$ et on a posé

$$\displaystyle \kappa_{N,X''}(g)=\nu(A_{I''}) \int_{A_{I''}(F)} \kappa_N(\gamma_{X''}^{-1}ag) da$$

Fixons $I=I'I''\in\mathcal{T}(G_{\tilde{x}})$. On a le lemme suivant dont la démonstration est tout à fait analogue à celle de 8.1.1 [B]

\begin{lem}
Soit $\omega_{I''}\subset \mathfrak{i}''(F)$ un sous-ensemble compact. Il existe un entier $k\in \mathbb{N}$ tel que

$$\kappa_{N,X''}(g)<< N^k (1+ |log \; D^{G''_{\tilde{x}}}(X'')|)^k (1+ |log \; |Q_S(X'')|_F|)^k \sigma(g)^k$$

\noindent pour tous $N\geqslant 1$, $g\in G(F)$, $X''\in\mathfrak{i}''(F)^S\cap \omega_{I''}$.
\end{lem}

Pour $\epsilon>0$, notons $\mathfrak{i}(F)[> \epsilon]$ l'ensemble des $X\in \mathfrak{i}(F)$ tels que $|det \; ad(X'')_{|\mathfrak{g}''/\mathfrak{i}''}|_F>\epsilon$, $|det \; ad(X')_{|\mathfrak{g}'/\mathfrak{i}'}|_F>\epsilon$, $|Q_S(X'')|_F> \epsilon$. Soit $\mathfrak{i}(F)[\leqslant \epsilon]$ le complémentaire de $\mathfrak{i}(F)[> \epsilon]$ dans $\mathfrak{i}(F)$.

\begin{lem}
Il existe un entier $b\geqslant 1$ tel que

\[\begin{aligned}
\int_{\mathfrak{i}'(F)\times \mathfrak{i}''(F)^S[\leqslant N^{-b}]} |\hat{j}_S(X')| D^{G'_{\tilde{x}}}(X')D^{G''_{\tilde{x}}}(X'')^{1/2} \int_{I'(F)A_{I''}(F)\backslash G(F)} & |{}^g \tilde{f}_{\tilde{x},\omega}^\sharp(X)| \\
 & \kappa_{N,X''}(g) dg dX''dX'<< N^{-1}
\end{aligned}\]

\noindent pour tout $N\geqslant 1$.
\end{lem}

\ul{Preuve}: Considérons l'intégrale intérieure et décomposons la en une double intégrale sur $G_{\tilde{x}}(F)\backslash G(F)$ et sur $I'(F)A_{I''}(F)\backslash G_{\tilde{x}}(F)$. L'intégrale sur $G_{\tilde{x}}(F)\backslash G(F)$ est localement constante et d'après une propriété des bons voisinages, elle est à support compact. Il existe un sous-ensemble compact $\omega_I\subset \mathfrak{i}(F)$ tel que l'intégrale intérieure soit nulle pour $X\notin\omega_I$. D'après ces deux remarques et le lemme 3.7.1, il existe un entier $k$ tel que notre expression soit essentiellement majorée par

\[\begin{aligned}
\mbox{(2)}\;\;\; \displaystyle N^k & \int_{\mathfrak{i}'(F)\times \mathfrak{i}''(F)^S[\leqslant N^{-b}]\cap \omega_I} |\hat{j}_S(X')| D^{G'_{\tilde{x}}}(X')D^{G''_{\tilde{x}}}(X'')^{1/2} \\
 & \int_{I'(F)A_{I''}(F)\backslash G_{\tilde{x}}(F)} \varphi(g^{-1} Xg) (1+ |log \; D^{G''_{\tilde{x}}}(X'')|)^k (1+ |log \; |Q_S(X'')|_F|)^k \sigma_I(g)^k dg dX''dX'
\end{aligned}\]

\noindent où $\varphi\in C_c^\infty(\mathfrak{g}_{\tilde{x}}(F))$ est une fonction à valeurs réelles positives et où on a posé

$$\sigma_I(g)=inf\{\sigma(tg);\; t\in I'(F)A_{I''}(F)\}$$

\noindent D'après le lemme 7.2.1 de [B] (qui est dû à Arthur), pour tout $X\in \omega_I$ et pour tout $g\in G_{\tilde{x}}(F)$ tels que $\varphi(g^{-1} Xg)\neq 0$ on a une majoration

$$\sigma_I(g)<< (1+ |log\; D^{G_{\tilde{x}}}(X)|)$$

\noindent D'après Harish-Chandra les fonctions

$$X'\mapsto |\hat{j}_S(X')| D^{G'_{\tilde{x}}}(X')^{1/2}$$

\noindent et

$$X\mapsto D^{G_{\tilde{x}}}(X)^{1/2}\int_{I'(F)A_{I''}(F)\backslash G_{\tilde{x}}(F)} \varphi(g^{-1} Xg) dX$$

\noindent sont localement bornées. Par conséquent, notre expression (2) est essentiellement majorée par

$$N^k\int_{\mathfrak{i}'(F)\times \mathfrak{i}''(F)^S[\leqslant N^{-b}]\cap \omega_I}(1+ |log \; D^{G''_{\tilde{x}}}(X'')|)^{2k} (1+ |log \; |Q_S(X'')|_F|)^k(1+|log \; D^{G'_{\tilde{x}}}(X')|)^k dX'' dX'$$

Il est alors aisé de montrer que cette intégrale vérifie la majoration de l'énoncé pour $b$ assez grand $\blacksquare$

\subsection{Changement de fonction de troncature}

Rappelons que l'on a fixé $I=I'I''\in\mathcal{T}(G_{\tilde{x}})$. Notons $\tilde{M}_\natural$ le commutant de $A_{I''}$ dans $\tilde{G}$. Il contient $I''G'\tilde{x}$ et c'est un Levi tordu de $\tilde{G}$. On a $A_{\tilde{M}_\natural}=A_{I''}$. Quitte à conjuguer $\tilde{x}$ et $I''$, on peut supposer que $\tilde{M}_{\natural}$ est semi-standard. Soit $Y\in\mathcal{A}_{M_{min}}$ tel que $\alpha(Y)>0$ pour tout $\alpha\in \Delta$. Pour $P'\in\mathcal{P}(M_{min})$, il existe un unique $w\in W^G$ tel que $wP_{min}=P'$ et on pose alors $Y_{P'}=wY$. Pour $\tilde{Q}\in\mathcal{P}(\tilde{M}_\natural)$, on définit $Y_{\tilde{Q}}$ comme le projeté sur $\mathcal{A}_{\tilde{M}_\natural}$ de $Y_{P'}$ pour n'importe quel $P'\in\mathcal{P}(M_{min})$ vérifiant $P'\subset Q$. Pour $g\in G(F)$, on pose $Y_{\tilde{Q}}(g)=Y_{\tilde{Q}}-H_{\overline{\tilde{Q}}}(g)$ pour tout $\tilde{Q}\in\mathcal{P}(\tilde{M}_\natural)$ où $\overline{\tilde{Q}}$ désigne le parabolique tordu opposé à $\tilde{Q}$. La famille $\mathcal{Y}(g)=(Y_{\tilde{Q}}(g))_{\tilde{Q}\in\mathcal{P}(\tilde{M}_\natural)}$ est $(\tilde{G},\tilde{M}_\natural)$-orthogonale et il existe $c_1>0$ tel que si

$$inf\{\sigma(mg);\; m\in M_\natural(F)\}\leqslant c_1 inf\{\alpha(Y);\; \alpha\in\Delta\}$$

\noindent la famille $\mathcal{Y}(g)$ soit $(\tilde{G},\tilde{M}_\natural)$-orthogonale positive. Supposons que ce soit le cas. On définit alors $\sigma^{\tilde{G}}_{\tilde{M}_\natural}(.,\mathcal{Y}(g))$ comme la fonction caractéristique dans $\mathcal{A}_{\tilde{M}_\natural}$ de l'enveloppe convexe des $Y_{\tilde{Q}}(g)$, $\tilde{Q}\in\mathcal{P}(\tilde{M}_\natural)$, et on pose

$$\displaystyle v(g,Y)=\nu(A_{I''})\int_{A_{I''}(F)} \sigma^{\tilde{G}}_{\tilde{M}_\natural}(H_{\tilde{M}_\natural},\mathcal{Y}(g))da$$

On a le fait suivant \\

 (1) Il existe $c_2>0$ et un sous-ensemble compact $\omega_I\subset \mathfrak{i}(F)$ tels que pour tous $X\in\mathfrak{i}(F)[>N^{-b}]$, $g\in G(F)$ vérifiant ${}^g\tilde{f}^\sharp_{\tilde{x},\omega}(X)\neq 0$, on ait $X\in\omega_I$ et
 $$\sigma_I(g)<c_2 log(N)$$
 
 En effet, on a vu lors de la preuve du lemme 3.7.2 qu'il existe un sous-ensemble compact $\omega_I\subset \mathfrak{i}(F)$ tel que pour tous $X$ et $g$ comme précédemment, on ait $X\in\omega_I$ et $\sigma_I(g)<<1+|log D^{G_{\tilde{x}}}(X)|$. Or, pour $X\in\omega_I\cap\mathfrak{i}(F)[>N^{-b}]$, on a $|log D^{G_{\tilde{x}}}(X)|<<log(N)$. D'où (1). \\
 
 Ceci permet de vérifier que le membre de droite de l'égalité ci-dessous est bien défini.
 
\begin{prop}
Il existe un réel $c_3>c_2/c_1$ et un entier $N_0\geqslant 1$ tels que, si $N\geqslant N_0$ et

$$inf\{\alpha(Y);\; \alpha\in\Delta\}>c_3 log(N)$$
on ait l'égalité

$$\displaystyle\int_{I'(F)A_{I''}(F)\backslash G(F)} {}^g\tilde{f}_{\tilde{x},\omega}^\sharp(X) \kappa_{N,X''}(g) dg=\int_{I'(F)A_{I''}(F)\backslash G(F)} {}^g\tilde{f}_{\tilde{x},\omega}^\sharp(X) v(g,Y) dg$$

\noindent pour tout $X\in\mathfrak{i}'(F)\times\mathfrak{i}''(F)^S[>N^{-b}]$.
\end{prop}

\ul{Preuve}: Soit $N\geqslant 1$ un entier. Posons

$$\mathcal{X}_N=\omega_I\cap\mathfrak{i}'(F)\times\mathfrak{i}''(F)^S[>N^{-b}]$$

\noindent et

$$\Omega_N=\{g\in G(F);\; \sigma_I(g)<c_2 log(N)\}$$

Soient $c_3>c_2/c_1$ et $Y\in\mathcal{A}_{M_{min}}$ vérifiant

$$inf\{\alpha(Y);\; \alpha\in\Delta\}>c_3 log(N)$$

Ainsi $v(Y,g)$ est bien défini pour tout $g\in \Omega_N$. On va montrer que l'égalité de l'énoncé est vérifiée si $N$ et $c_3$ sont assez grands. Notons $J(N,X)$ resp. $J(Y,X)$ le membre de gauche resp. de droite de l'égalité de l'énoncé. D'après la remarque (1), il suffit de montrer l'égalité pour $X\in \mathcal{X}_N$ et on peut restreindre les intégrales à $I'(F)A_{I''}(F)\backslash \Omega_N$. Soit $Z\in\mathcal{A}_{M_{min}}$ vérifiant

$$inf\{\alpha(Z);\; \alpha\in\Delta\}\geqslant c_2/c_1 log(N)$$

Pour tout $g\in\Omega_N$, on définit comme pour $Y$ une $(\tilde{G},\tilde{M}_\natural)$-famille orthogonale positive $\mathcal{Z}(g)=(Z_{\tilde{Q}}(g))_{\tilde{Q}\in\mathcal{P}(\tilde{M}_\natural)}$. Pour $\tilde{Q}=\tilde{L}U_Q\in\mathcal{F}(\tilde{M}_\natural)$, notons $\sigma^{\tilde{Q}}_{\tilde{M}_\natural}(.,\mathcal{Z}(g))$ la fonction caractéristique dans $\mathcal{A}_{\tilde{M}_\natural}$ de la somme de $\mathcal{A}_{\tilde{L}}$ et de l'enveloppe convexe des $Z(g)_{\tilde{P}'}$ pour $\tilde{P}'\in\mathcal{P}(\tilde{M}_\natural)$ tel que $\tilde{P}'\subset \tilde{Q}$. On note aussi $\tau_{\tilde{Q}}$ la fonction caractéristique de $\mathcal{A}^{\tilde{L}}_{\tilde{M}_\natural}\oplus\mathcal{A}^+_{\tilde{Q}}$ dans $\mathcal{A}_{\tilde{M}_\natural}$. On a l'égalité

$$\displaystyle \sum_{\tilde{Q}\in\mathcal{F}(\tilde{M}_\natural)} \sigma^{\tilde{Q}}_{\tilde{M}_\natural}(\lambda,\mathcal{Z}(g))\tau_{\tilde{Q}}(\lambda-Z(g)_{\tilde{Q}})=1$$

\noindent pour tout $\lambda\in\mathcal{A}_{\tilde{M}_\natural}$. On a par conséquent pour tous $g\in\Omega_N$, $X\in\mathcal{X}_N$,

$$\displaystyle\mbox{(2)}\;\;\; v(g,Y)=\sum_{\tilde{Q}\in\mathcal{F}(\tilde{M}_\natural)} v(\tilde{Q},Y,g)$$
$$\displaystyle\mbox{(3)}\;\;\; \kappa_{N,X''}(g)=\sum_{\tilde{Q}\in\mathcal{F}(\tilde{M}_\natural)} \kappa_{N,X''}(\tilde{Q},g)$$

\noindent où on a posé

$$\displaystyle v(\tilde{Q},Y,g)=\nu(A_{I''})\int_{A_{I''}(F)} \sigma^{\tilde{G}}_{\tilde{M}_\natural}(H_{\tilde{M}_\natural}(a),\mathcal{Y}(g)) \sigma^{\tilde{Q}}_{\tilde{M}_\natural}(H_{\tilde{M}_\natural}(a),\mathcal{Z}(g))\tau_{\tilde{Q}}(H_{\tilde{M}_\natural}(a)-Z(g)_{\tilde{Q}}) da$$

\noindent et

$$\displaystyle \kappa_{N,X''}(\tilde{Q},g)=\nu(A_{I''})\int_{A_{I''}(F)} \kappa_N(\gamma_{X''}^{-1}ag) \sigma^{\tilde{Q}}_{\tilde{M}_\natural}(H_{\tilde{M}_\natural}(a),\mathcal{Z}(g))\tau_{\tilde{Q}}(H_{\tilde{M}_\natural}(a)-Z(g)_{\tilde{Q}}) da$$

Les fonctions $\kappa_{N,X''}(\tilde{Q},.)$ et $v(\tilde{Q},Y,.)$ sont invariantes à gauche par $I'(F)A_{I''}(F)$. Par conséquent, on peut décomposer $J(Y,X)$ et $J(N,X)$ suivant les égalités (2) et (3). On obtient

$$\displaystyle J(Y,X)=\sum_{\tilde{Q}\in\mathcal{F}(\tilde{M}_\natural)} J(\tilde{Q},Y,X)$$
$$\displaystyle J(N,X)=\sum_{\tilde{Q}\in\mathcal{F}(\tilde{M}_\natural)} J(\tilde{Q},N,X)$$

\noindent avec une définition plus ou moins évidente des termes $J(\tilde{Q},Y,X)$ et $J(\tilde{Q},N,X)$. Montrons le fait suivant \\

 (4) Si $Z$ vérifie les inégalités
 
$$sup\{\alpha(Z);\; \alpha\in\Delta\}\leqslant \left\{
    \begin{array}{ll}
        inf\{\alpha(Y);\alpha\in \Delta\}\\
        log(N)^2
    \end{array}
\right. $$

\noindent on a $J(\tilde{G},Y,X)=J(\tilde{G},N,X)$ pour tout $X\in\mathcal{X}_N$ et pour $N$ assez grand. \\

La première inégalité entraîne que l'enveloppe convexe des $Z(g)_{\tilde{Q}}$, $\tilde{Q}\in\mathcal{P}(\tilde{M}_\natural)$, est inclus dans l'enveloppe convexe des $Y(g)_{\tilde{Q}}$, $\tilde{Q}\in\mathcal{P}(\tilde{M}_\natural)$, pour tout $g\in\Omega_N$. Par conséquent $v(\tilde{G},Y,g)=v(g,Z)$ pour tout $g\in \Omega_N$. \\
Il existe $c>0$ tel que $\sigma(g)<cN$ entraîne $\kappa_N(g)=1$ pour tout $g\in G(F)$. Pour $X\in \mathcal{X}_N$, on a $\sigma(\gamma_{X''})<<log(N)$. D'après la deuxième inégalité, pour $a\in A_{I''}(F)$ et $g\in \Omega_N$, l'égalité $\sigma_{\tilde{M}_\natural}^{\tilde{G}}(H_{\tilde{M}_\natural}(a),\mathcal{Z}(g))=1$ entraîne $\sigma(ag)<<log(N)^2$ donc entraîne aussi $\kappa_N(\gamma_{X''}^{-1}ag)=1$ si $N$ est assez grand. On a donc aussi pour $N$ assez grand $\kappa_{N,X''}(\tilde{G},g)=v(g,Z)$ pour tous $g\in\Omega_N$, $X\in\mathcal{X}_N$. On en déduit (4). \\

Fixons à présent $\tilde{Q}=\tilde{L}U_Q\in \mathcal{F}(\tilde{M}_\natural)$, $\tilde{Q}\neq\tilde{G}$. Etablissons la propriété suivante \\

 (5) Il existe $c_4>0$ tel que si
 
$$inf\{\alpha(Z);\; \alpha\in\Delta\}\geqslant c_4 log(N)$$

\noindent et si $c_3$ et $N$ sont assez grands alors $J(\tilde{Q},N,X)=J(\tilde{Q},Y,X)=0$ pour tout $X\in\mathcal{X}_N$. \\

Remarquons tout d'abord que les points (4) et (5) entraînent la proposition: pour $c_3$ assez grand on peut trouver un $Z$ qui vérifie les inégalités (4) et (5). En sommant alors les égalités $J(\tilde{Q},N,X)=J(\tilde{Q},Y,X)$ pour $\tilde{Q}\in\mathcal{F}(\tilde{M}_\natural)$, on obtient le résultat voulu. \\
 On a les égalités

$$\displaystyle J(\tilde{Q},N,X)=\int_K \int_{I'(F)A_{I''}(F)\backslash L(F)}\int_{U_{\overline{Q}}(F)} {}^{\overline{u}lk} \tilde{f}^\sharp_{\tilde{x},\omega}(X) \kappa_{N,X''}(\tilde{Q},\overline{u}lk)d\overline{u} \delta_Q(l) dl dk$$

$$\displaystyle J(\tilde{Q},Y,X)=\int_K \int_{I'(F)A_{I''}(F)\backslash L(F)}\int_{U_{\overline{Q}}(F)} {}^{\overline{u}lk} \tilde{f}^\sharp_{\tilde{x},\omega}(X) v(\tilde{Q},Y,\overline{u}lk)d\overline{u} \delta_Q(l) dl dk$$

Il existe $c>0$ tel que pour $\overline{u}\in U_{\overline{Q}}(F)$, $l\in L(F)$ et $k\in K$ vérifiant $\overline{u}lk\in\Omega_N$, on ait $\sigma(\overline{u}),\sigma_I(l),\sigma_I(\overline{u}lk)\leqslant c log(N)$. Supposons que pour de tels éléments on ait $\kappa_{N,X''}(\tilde{Q},\overline{u}lk)=\kappa_{N,X''}(\tilde{Q},lk)$ pour tout $X\in \mathcal{X}_N$ et $v(\tilde{Q},Y,\overline{u}lk)=v(\tilde{Q},Y,lk)$. Alors il apparaît dans les deux intégrales précédentes l'intégrale intérieure

$$\displaystyle \int_{U_{\overline{Q}}(F)} {}^{\overline{u}lk} \tilde{f}^\sharp_{\tilde{x},\omega}(X) d\overline{u}$$

\noindent qui est nulle d'après [W3] lemme 5.5(i): c'est ici que l'on utilise l'hypothèse que $\tilde{f}$ est très cuspidale. Il suffit donc pour établir (5), de montrer la propriété suivante. \\

 (6) Soit $c>0$, il existe $c_4>0$ tel que si
 
 $$inf\{\alpha(Z);\; \alpha\in\Delta\}\geqslant c_4 log(N)$$
 
\noindent et si $c_3$ et $N$ sont assez grands alors $\kappa_{N,X''}(\tilde{Q},\overline{u}g)=\kappa_{N,X''}(\tilde{Q},g)$ et $v(\tilde{Q},Y,\overline{u}g)=v(\tilde{Q},Y,g)$ pour tout $X\in\mathcal{X}_N$ et pour tous $\overline{u}\in U_{\overline{Q}}(F)$, $g\in G(F)$ vérifiant $\sigma(\overline{u}),\sigma(g),\sigma(\overline{u}g)\leqslant c log(N)$. \\
 
 Fixons $c_4>0$ et montrons que (6) est vérifié si $c_4$, $c_3$ et $N$ sont assez grands. Pour $c_3$ et $c_4$ assez grands toutes les $(\tilde{G},\tilde{M}_\natural)$-familles orthogonales $(T_{\tilde{P}'})_{\tilde{P}'\in\mathcal{P}(\tilde{M}_\natural)}$ qui apparaissent par la suite vérifient la propriété suivante: pour tout $\tilde{P}'\in\mathcal{P}(\tilde{M}_\natural)$, on a $T_{\tilde{P}'}\in\mathcal{A}_{\tilde{P}'}^+$. Montrons d'abord l'assertion sur la fonction $v(\tilde{Q},Y,.)$. Soit $g\in G(F)$ tel que $\sigma(g)\leqslant clog(N)$. Par définition on a
 
$$\displaystyle v(\tilde{Q},Y,g)=\nu(A_{I''})\int_{A_{I''}(F)} \sigma^{\tilde{G}}_{\tilde{M}_\natural}(H_{\tilde{M}_\natural}(a),\mathcal{Y}(g)) \sigma^{\tilde{Q}}_{\tilde{M}_\natural}(H_{\tilde{M}_\natural}(a),\mathcal{Z}(g))\tau_{\tilde{Q}}(H_{\tilde{M}_\natural}(a)-Z(g)_{\tilde{Q}}) da$$

Les fonctions $\sigma^{\tilde{Q}}_{\tilde{M}_\natural}(.,\mathcal{Z}(g))$ et $\tau_{\tilde{Q}}(.-Z(g)_{\tilde{Q}})$ ne dépendent que de $Z(g)_{\tilde{P}'}$ pour $\tilde{P}'\in\mathcal{P}(\tilde{M}_\natural)$ vérifiant $\tilde{P}'\subset \tilde{Q}$. Par conséquent, elles sont invariantes par translation à gauche de $g$ par $U_{\overline{Q}}(F)$. De plus, pour $\lambda\in\mathcal{A}_{\tilde{M}_\natural}$ vérifiant $\sigma^{\tilde{Q}}_{\tilde{M}_\natural}(\lambda,\mathcal{Z}(g))\tau_{\tilde{Q}}(\lambda-Z(g)_{\tilde{Q}})=1$, il existe $\tilde{P}'\in\mathcal{P}(\tilde{M}_\natural)$, $\tilde{P}'\subset \tilde{Q}$ tel que $\lambda\in \mathcal{A}_{\tilde{P}'}^+$. Or, la restriction de $\sigma^{\tilde{G}}_{\tilde{M}_\natural}(.,\mathcal{Y}(g))$ à $\mathcal{A}_{\tilde{P}'}^+$ ne dépend que de $Y(g)_{\tilde{P}'}$ donc est aussi invariante par translation à gauche de $g$ par $U_{\overline{Q}}(F)$. Cela établit le résultat pour $v(\tilde{Q},Y,.)$. \\

Montrons maintenant (6) pour $\kappa_{N,X''}(\tilde{Q},.)$. La fonction

$$g\mapsto \sigma^{\tilde{Q}}_{\tilde{M}_\natural}(.,\mathcal{Z}(g))\tau_{\tilde{Q}}(.-Z(g)_{\tilde{Q}})$$

\noindent définie pour $\sigma(g)\leqslant c log(N)$ est toujours invariante à gauche par $U_{\overline{Q}}(F)$. Soit $c'>0$. Pour $c_4$ assez grand l'égalité $\sigma^{\tilde{Q}}_{\tilde{M}_\natural}(\lambda,\mathcal{Z}(g))\tau_{\tilde{Q}}(\lambda-Z(g)_{\tilde{Q}})=1$, $\lambda\in\mathcal{A}_{\tilde{M}_\natural}$, entraîne $\alpha(\lambda)\geqslant c' log(N)$ pour toute racine $\alpha$ de $A_{I''}$ dans $\mathfrak{u}_Q$. Le résultat est alors une conséquence du lemme suivant

\begin{lem}
Soit $c>0$, il existe un entier $N_1\geqslant 1$ et $c'>0$ tels que pour tous $\overline{u}\in U_{\overline{Q}}(F)$, $g\in G(F)$ vérifiant $\sigma(\overline{u}),\sigma(g),\sigma(\overline{u}g)\leqslant c log(N)$, pour tout $a\in A_{I''}$ vérifiant $\alpha(H_{\tilde{M}_\natural}(a))\geqslant c' log(N)$ pour toute racine $\alpha$ de $A_{I''}$ dans $\mathfrak{u}_Q$ et pour tout $N\geqslant N_1$, on ait l'égalité

$$\kappa_N(\gamma_{X''}^{-1}ag)=\kappa_N(\gamma_{X''}^{-1}a\overline{u}g)$$
pour tout $X\in\mathcal{X}_N$.
\end{lem}

\ul{Preuve}: Pour $x$ un réel quelconque, on pose $\pi_E^x=\pi_E^{E(x)}$ où $E(x)$ désigne la partie entière de $x$. Il existe un réel $C>0$ qui vérifie les conditions suivantes

\begin{enumerate}
\item $\pi_E^{C\sigma(g)}R\subset gR\subset \pi_E^{-C\sigma(g)}R$ pour tout $g\in G(F)$
\item $\sigma(\gamma_{X''})\leqslant Clog(N)$ pour tout $N\geqslant 2$ et pour tout $X\in\mathcal{X}_N$
\item ${X''_{\Sigma}}^iR\subset \pi_E^{-C}R$ pour tout $X''\in \omega_I\cap \mathfrak{i}''(F)^S$ et pour tout $i=0,\ldots,d''-1$
\item $\pi_E^{Clog(N)} R''\subset R_{X''_{\Sigma},z_r}$ pour tout $N\geqslant 2$ et pour tout $X\in\mathcal{X}_N$
\item $\pi_E^C R\subset R''\oplus R'$ où $R'=R\cap V'$
\item $\pi_E^C R'' \subset {R''}^*\subset \pi_E^{-C}R''$ où ${R''}^*=\tilde{x}^{-1} {R''}^\vee$.
\end{enumerate}
\noindent Pour les conditions 2, 3 et 4 on utilise les propriétés des applications $X''\mapsto \gamma_{X''}$ et $X''\mapsto X''_{\Sigma}$ rappelées en 3.7. Fixons $N_1\geqslant 2$ un entier et $c'>0$ et montrons que la conclusion du lemme est vérifiée si $N_1$ et $c'$ sont assez grands. Considérons des éléments $a,\overline{u},g,N$ et $X$ comme dans le lemme. Supposons que $\kappa_N(\gamma_{X''}^{-1}ag)=1$ et montrons que $\kappa_N(\gamma_{X''}^{-1}a\overline{u}g)=1$ (cela suffit à montrer le lemme: il suffit de changer $g$ en $\overline{u}g$ et $\overline{u}$ en $\overline{u}^{-1}$ pour obtenir l'autre inégalité). On a alors $a^{-1}\gamma_{X''}z_r\in \pi_E^{-N}g(R)$ ce qui entraîne $\gamma_{X''}^{-1}a^{-1}\gamma_{X''}z_r\in \pi_E^{-N-(c+C)Clog(N)}R$. Comme $X''_{\Sigma}=\gamma_{X''}^{-1}X''\gamma_{X''}$ commute à $\gamma_{X''}^{-1}a^{-1}\gamma_{X''}$, on a $\gamma_{X''}^{-1}a^{-1}\gamma_{X''}R_{X''_{\Sigma},z_r}\subset \pi_E^{-N-(c+C)Clog(N)-C}R$ d'où

$$\gamma_{X''}^{-1}a^{-1}\gamma_{X''}R''\subset \pi_E^{-N-(c+C+1)Clog(N)-C}R''$$

\noindent Puisque $\gamma_{X''}^{-1}a^{-1}\gamma_{X''}\in G''_{\tilde{x}}$, on a aussi

$$\pi_E^{N+(c+C+1)Clog(N)+C}{R''}^* \subset \gamma_{X''}^{-1}a^{-1}\gamma_{X''}{R''}^*$$
On en déduit que

$$\mbox{(7)}\;\;\; \pi_E^{N+(c+3C+1)Clog(N)+4C}R \subset a^{-1}R\subset \pi_E^{-N-(c+3C+1)Clog(N)-2C}R$$

\noindent Soit $C'>0$. Si $c'$ est assez grand, on a $(a\overline{u}^{-1}a^{-1}-Id)R\subset \pi_E^{C' log(N)}R$. Combiné à (7), cela entraîne

\[\begin{aligned}
g^{-1}a^{-1}\gamma_{X''}z_r-g^{-1}\overline{u}^{-1}a^{-1}\gamma_{X''}z_r\in \pi_E^{-N}R 
\end{aligned}\]

\noindent si $C'$ est assez grand. Par conséquent, on a aussi

$$\mbox{(8)}\;\;\; g^{-1}\overline{u}^{-1}a^{-1}\gamma_{X''}z_r\in \pi_E^{-N}R$$

\noindent si $c'$ est assez grand. De plus, d'après (7) si $N_1$ est assez grand on a

$$\mbox{(9)} \;\;\; \pi_E^{2N}R \subset g^{-1}\overline{u}^{-1}a^{-1}\gamma_{X''}R\subset \pi_E^{-2N}R$$

Par définition de $\kappa_N$, (8) et (9) entraînent $\kappa_N(\gamma_{X''}^{-1}a\overline{u}g)=1$. D'où le résultat $\blacksquare$

\begin{prop}
Il existe un entier $N_0\geqslant 1$ tel que, pour tout $N\geqslant N_0$ et tout $X\in (\mathfrak{i}'(F)\times\mathfrak{i}''(F)^S)[>N^{-b}]$, on ait l'égalité

$$\displaystyle\int_{I'(F)A_{I''}(F)\backslash G(F)} {}^g\tilde{f}_{\tilde{x},\omega}^\sharp(X)\kappa_{N,X''}(g) dg= \left\{
    \begin{array}{ll}
        0 \;\mbox{si } A_{I'}\neq \{1\}\\
        \nu(I')\nu(A_{I''}) \Theta_{\tilde{f},\tilde{x},\omega}^\sharp(X) \;\mbox{sinon}
    \end{array}
\right. $$
\end{prop}

\ul{Preuve}: D'après la proposition 3.8.1, on peut en tout cas remplacer $\kappa_{N,X''}(g)$ par $v(g,Y)$ pour $Y\in \mathcal{A}_{M_{min}}$ vérifiant une minoration

$$inf\{\alpha(Y);\;\;\alpha\in \Delta\}>>log(N)$$

Soit $\mathcal{R}\subset \mathcal{A}_{M_{min},F}\otimes \mathbb{Q}$ un $\mathbb{Z}$-réseau. Notons $\mathcal{R}^\vee$ l'ensemble des éléments $\lambda\in \mathcal{A}_{M_{min}}^*$ tels que $\lambda(\zeta)\in 2\pi\mathbb{Z}$ pour tout $\zeta\in \mathcal{R}$. Alors d'après le lemme 4.7 de [W4], il existe un sous-ensemble fini $\mathbf{\Xi}\subset i\mathcal{A}_{M_{min}}^*/i\mathcal{R}^\vee$, un entier $M\geqslant 1$ ainsi qu'une famille de fonctions $g\mapsto c_{\mathcal{R}}(\Lambda,k,g)$ pour $\Lambda\in\mathbf{\Xi}$ et $0\leqslant k\leqslant M$ telles que, pour tout $g\in G(F)$ et pour tout $Y\in\mathcal{R}$ vérifiant une majoration $\alpha(Y)>>\sigma(g)$ pour tout $\alpha\in\Delta$, on ait l'égalité

$$\displaystyle v(g,Y)=\sum_{\Lambda\in \mathbf{\Xi},\; 0\leqslant k\leqslant M} c_{\mathcal{R}}(\Lambda,k,g)e^{\Lambda(Y)} Y^k$$

\noindent Les fonctions $Y\mapsto e^{\Lambda(Y)}Y^k$, $k\in\mathbb{N}$, $\Lambda\in i\mathcal{A}_{M_{min}}^*/i\mathcal{R}^\vee$ sont linéairement indépendantes sur le cône $\{Y\in\mathcal{R}: \; \sigma(g)\leqslant c_1\alpha(Y) \; \forall \alpha\in \Delta\}$. Puisque l'intégrale que l'on considère ne dépend pas de $Y$ on peut aussi bien remplacer la fonction $v(g,Y)$ par le terme constant $c_{\mathcal{R}}(0,0,g)$. On peut bien sûr appliquer ce qui précède au réseau $\frac{1}{k}\mathcal{R}$ pour tout $k\in\mathbb{N}^*$. Toujours d'après le lemme 4.7 de [W4], on a une inégalité

$$\displaystyle |c_{\frac{1}{k}\mathcal{R}}(0,0,g)-(-1)^{a_{\tilde{M}_\natural}}v_{\tilde{M}_\natural}(g)|<<k^{-1}\sigma_I(g)^{a_{\tilde{M}_\natural}}$$

\noindent pour tout $g\in G(F)$ et pour tout $k\in\mathbb{N}^*$. Puisque l'intégrale

$$\displaystyle \int_{I'(F)A_{I''}(F)\backslash G(F)} |{}^g\tilde{f}_{\tilde{x},\omega}(X)| \sigma_I(g)^{a_{\tilde{M}_\natural}} dg$$

\noindent est convergente, on obtient par passage à la limite l'égalité

$$\displaystyle \int_{I'(F)A_{I''}(F)\backslash G(F)} {}^g\tilde{f}_{\tilde{x},\omega}^\sharp(X)\kappa_{N,X''}(g) dg=(-1)^{a_{\tilde{M}_\natural}}\int_{I'(F)A_{I''}(F)\backslash G(F)} {}^g\tilde{f}_{\tilde{x},\omega}^\sharp(X)v_{\tilde{M}_\natural}(g) dg$$

Si $A_{I'}\neq \{1\}$, cette dernière intégrale vaut $0$ d'après [W3] 1.8(2), sinon d'après la proposition 1.8 de [W3], elle vaut $\nu(I)mes(I(F)/I'(F)A_{I''}(F))\Theta_{\tilde{f},\tilde{x},\omega}^\sharp(X)$. Il ne reste plus qu'à vérifier l'égalité $\nu(I)mes(I(F)/I'(F)A_{I''}(F))=\nu(I')\nu(A_{I''})$ pour obtenir le résultat $\blacksquare$

\subsection{Preuve du théorème 3.5.1}

D'après l'égalité 3.7(1), le lemme 3.7.2 et la proposition 3.8.2, $J_{N,\tilde{x},\omega}(\Theta,\tilde{f})$ admet une limite lorsque $N$ tend vers l'infini et cette limite vaut

\[\begin{aligned}
\displaystyle K_{\tilde{x},\omega}(\Theta,\tilde{f})=\sum_{(I',I'')\in \mathcal{T}_{ell}(H'_{\tilde{x}})\times \mathcal{T}(G''_{\tilde{x}})} \nu(I') |W(G_{\tilde{x}},I)|^{-1} \int_{\mathfrak{i}'(F)\times \mathfrak{i}''(F)^S} & \hat{j}_S(X') D^{G'_{\tilde{x}}}(X') D^{G''_{\tilde{x}}}(X'')^{1/2} \\
 & \Theta_{\tilde{f},\tilde{x},\omega}^\sharp(X) dX''dX'
\end{aligned}\]

\noindent La fonction $\Theta_{\tilde{f},\tilde{x},\omega}$ est un quasi-caractère de $\mathfrak{g}_{\tilde{x}}(F)$ à support dans $\omega=\omega'\times \omega''$. Un tel quasi-caractère se décompose comme combinaison linéaire de produits d'un quasi-caractère sur $\omega'$ et d'un quasi-caractère sur $\omega''$. On a donc

$$\displaystyle \Theta_{\tilde{f},\tilde{x},\omega}(X)=\sum_{b\in B} \Theta_{\tilde{f},b}'(X') \Theta_{\tilde{f},b}''(X'')$$

\noindent pour tout $X\in \mathfrak{g}_{\tilde{x}}(F)$ où $\Theta_{\tilde{f},b}',\Theta_{\tilde{f},b}''$ sont des quasi-caractères dans $\omega'$ et $\omega''$ respectivement. Rappelons que l'on a supposé $\Theta(X)=\hat{j}_S(X')\hat{j}^{H''_{\tilde{x}}}(S,X'')$ pour tout $X\in\omega$. Posons $\Theta'=\hat{j}_S$ et $\Theta''=\hat{j}^{H''_{\tilde{x}}}(S,.)$. On a alors

$$\displaystyle K_{\tilde{x},\omega}(\Theta,\tilde{f})=\sum_{b\in B} J'_b K''_b$$

\noindent où on a posé

$$\displaystyle J'_b=\sum_{I'\in \mathcal{T}_{ell}(H'_{\tilde{x}})} \nu(I') |W(G'_{\tilde{x}},I')|^{-1} \int_{\mathfrak{i}'(F)} D^{G'_{\tilde{x}}}(X')\Theta'(X')\Theta_{\tilde{f},b}'(X') dX'$$

\noindent et

$$\displaystyle K''_b=\sum_{I''\in \mathcal{T}(G''_{\tilde{x}})} |W(G''_{\tilde{x}},I'')|^{-1} \int_{\mathfrak{i}''(F)^S} D^{G''_{\tilde{x}}}(X'')^{1/2} \hat{\Theta}_{\tilde{f},b}''(X'') dX''$$

\noindent D'après la section 3.6, on a aussi

$$\displaystyle J_{geom,\tilde{x},\omega}(\Theta,\tilde{f})= \sum_{b\in B} J'_b J''_b$$

\noindent où

$$\displaystyle J''_b=\sum_{I''\in\mathcal{T}''} |W(H''_{\tilde{x}},I'')|^{-1} \lim\limits_{s\to 0^+} \int_{\mathfrak{i}''(F)} c_{\Theta''}(X'') c_{\Theta_{\tilde{f},b}''}(X'') D^{H''_{\tilde{x}}}(X'')^{1/2}D^{G''_{\tilde{x}}}(X'')^{1/2} \Delta''(X'')^{s-1/2} dX''$$

D'après le lemme 9.1.1 et le théorème 5.5.1 de [B], alliés à la proposition 6.4 de [W1], on a $J''_b=K''_b$. D'où l'égalité $K_{\tilde{x},\omega}(\Theta,\tilde{f})=J_{geom,\tilde{x},\omega}(\Theta,\tilde{f})$.

\section{Entrelacements tempérés et facteurs $\epsilon$}

Pour $c\in\mathbb{Z}$, on définit $U(F)_c$ comme le sous-groupe constitué des éléments $u\in U(F)$ qui vérifient $val_E(<uz_i,z_{i+1}^*>)\geqslant -c$ pour $i=-r,\ldots,r-1$. Soient $(\pi,E_\pi)\in Temp(G)$ et $(\sigma,E_\sigma)\in Temp(H)$. Pour $e,e'\in E_\pi$, $\epsilon,\epsilon'\in E_\sigma$ et $c\in\mathbb{Z}$, on pose

$$\displaystyle \mathcal{L}_{\pi,\sigma,c}(\epsilon'\otimes e',\epsilon\otimes e)=\int_{H(F)U(F)_c} (\sigma(h)\epsilon',\epsilon)(e',\pi(hu)e) \overline{\xi}(u) dudh$$

\noindent C'est une expression absolument convergente et elle ne dépend pas de $c$ pour $c$ assez grand (cf 5.1 [W3]). Posons

$$\mathcal{L}_{\pi,\sigma}(\epsilon'\otimes e',\epsilon\otimes e)=\lim\limits_{c\to\infty} \mathcal{L}_{\pi,\sigma,c}(\epsilon'\otimes e',\epsilon\otimes e)$$

\noindent pour tous $e,e'\in E_\pi$ et $\epsilon,\epsilon'\in E_\sigma$. C'est une forme sesquilinéaire symétrique sur $E_\sigma\otimes E_\pi$ et on vérifie facilement que l'on a

$$\mathcal{L}_{\pi,\sigma}(\sigma(h)\epsilon'\otimes e',\epsilon\otimes \pi(hu)e)=\xi(u)\mathcal{L}_{\pi,\sigma}(\epsilon'\otimes e',\epsilon\otimes e)$$

\noindent pour tous $e,e'\in E_\pi$, $\epsilon,\epsilon'\in E_\sigma$, $h\in H(F)$ et $u\in U(F)$. Définissons $Hom_{H,\xi}(\pi,\sigma)$ comme l'espace des homomorphismes $\ell: E_\pi\to E_\sigma$ vérifiant $\ell\circ \pi(hu)=\xi(u)\sigma(h)\circ\ell$ pour tout $h\in H(F)$ et pour tout $u\in U(F)$. D'après [AGRS] théorème 1 et [GGP] théorème 15.1, cet espace est de dimension au plus $1$. Par conséquent, il existe $\ell\in Hom_{H,\xi}(\pi,\sigma)$ et $C\in\mathbb{C}$ tels que

$$\mathcal{L}_{\pi,\sigma}(\epsilon'\otimes e',\epsilon\otimes e)=C(\epsilon',\ell(e))(\ell(e'),\epsilon)$$

\noindent pour tous $e,e'\in E_\pi$ et $\epsilon,\epsilon'\in E_\sigma$. \\

 Fixons une base $(v_i)_{i=1,\ldots,m}$ de $W$ et complétons la en une base $(v_i)_{i=1,\ldots,d}$ de $V$ par $v_{m+i}=z_{r+1-i}$ pour $i=1,\ldots,2r+1$. On identifie ainsi $G$ à $R_{E/F} GL_d$ et $H$ à $R_{E/F} GL_m$. Rappelons que l'on avait noté $B_d$ le sous-groupe de Borel standard de $R_{E/F} GL_d$ et $U_d$ son radical unipotent. Soient $\overline{B}_d$ le sous-groupe de Borel opposé à $B_d$ et $\overline{U}_d$ son radical unipotent. Pour $h,k,l\in\mathbb{N}$, $h<k\leqslant l\leqslant d$, on note $\overline{U}_{h,k,l}$ le sous-groupe de $\overline{U}_d$ consitué des éléments $\overline{u}\in\overline{U}_d$ tels que pour $i,j=1,\ldots,d$, $i\neq j$, $\overline{u}_{i,j}$ n'est non nul que si $(i,j)$ appartient au rectangle $1\leqslant j\leqslant h$ $k\leqslant i\leqslant l$. Fixons des modèles de Whittaker $\phi$ et $\phi'$ de $\pi$ et $\sigma$ comme en 2.2. Pour $\epsilon\in E_\sigma$, on définit $W_\epsilon: H(F)\to\mathbb{C}$ par
 
$$W_\epsilon(h)=\phi'(\sigma(h)\epsilon)$$
 
\noindent On définit de même $W_e$ pour $e\in E_\pi$. Pour $e\in E_\pi$ et $\epsilon\in E_\sigma$, on pose

$$\displaystyle L_{\pi,\sigma}(\epsilon,e)=\int_{\overline{U}_{m,m+1,m+r}(F)}\int_{U_m(F)\backslash H(F)} \overline{W_\epsilon}(h)W_e(h\overline{u}) |det(h)|_E^{-r} dh d\overline{u}$$

\noindent On a alors

\begin{prop}
\begin{enumerate}
\item L'expression ci-dessus est absolument convergente.
\item Il existe $C>0$ tel que, pour tous $e,e'\in E_\pi$ et $\epsilon,\epsilon'\in E_\sigma$, on ait

$$\mathcal{L}_{\pi,\sigma}(\epsilon'\otimes e',\epsilon\otimes e)=C L_{\pi,\sigma}(\epsilon',e) \overline{L_{\pi,\sigma}(\epsilon,e')}$$
\item Les formes sesquilinéaires $L_{\pi,\sigma}$ et $\mathcal{L}_{\pi,\sigma}$ sont non nulles.
\end{enumerate}
\end{prop}

C'est le contenu des lemmes 5.3 et 5.4 de [W3]. On en déduit en particulier que \\
\noindent $Hom_{H,\xi}(\pi,\sigma)\neq\{0\}$. On aura aussi besoin du lemme 5.6 de [W3] dont on rappelle ici l'énoncé.

\begin{lem}
Soit $Q=LU_Q$ un sous-groupe parabolique de $G$ et $\tau\in Temp(L)$. Pour tout $\lambda\in\mathcal{A}_{L,F}^*$, on réalise $\pi_\lambda=i_Q^G(\tau_\lambda)$ sur l'espace commun $\mathcal{K}_{Q,\tau}^G$. Cet espace est muni d'un produit scalaire qui est invariant par $\pi_\lambda$ pour tout $\lambda$. Soit aussi $\sigma\in Temp(H)$.
\begin{enumerate}
\item Soient $e,e'\in\mathcal{K}^G_{Q,\tau}$, et $\epsilon,\epsilon'\in E_\sigma$. La fonction $\lambda\mapsto \mathcal{L}_{\pi_\lambda,\sigma}(\epsilon'\otimes e',\epsilon\otimes e)$ est $C^\infty$ sur $i\mathcal{A}_{L,F}^*$.
\item Il existe des familles finies $(\epsilon_j)_{j=1,\ldots,n}$, $(e_j)_{j=1,\ldots,n}$, $(\varphi_j)_{j=1,\ldots,n}$, où $\epsilon_j\in E_\sigma$, $e_j\in \mathcal{K}^G_{Q,\tau}$ et $\varphi_j\in C^\infty(i\mathcal{A}_{L,F}^*)$ pour tout $j$, de sorte que

$$\displaystyle \sum_{j=1}^n \varphi_j(\lambda)\mathcal{L}_{\pi_\lambda,\sigma}(\epsilon_j\otimes e_j,\epsilon_j\otimes e_j)=1$$

\noindent pour tout $\lambda\in i\mathcal{A}_{L,F}^*$.
\end{enumerate}
\end{lem}

\subsection{Apparition des facteurs $\epsilon$}

\begin{prop}
Soient $\tilde{\pi}\in Temp(\tilde{G})$ et $\tilde{\sigma}\in Temp(\tilde{H})$. On a l'égalité

$$\mathcal{L}_{\pi,\sigma}(\epsilon'\otimes \tilde{\pi}(\tilde{y}) e',\tilde{\sigma}(\tilde{y}) \epsilon\otimes e)=\epsilon_\nu(\tilde{\pi}^\vee,\tilde{\sigma}) \mathcal{L}_{\pi,\sigma}(\epsilon'\otimes e', \epsilon\otimes e)$$

\noindent pour tous $e,e'\in E_\pi$, $\epsilon,\epsilon'\in E_\sigma$ et pour tout $\tilde{y}\in\tilde{H}(F)$.
\end{prop}

\ul{Preuve}: Puisque l'on a $\mathcal{L}_{\pi,\sigma}(\epsilon'\otimes \pi(h) e',\sigma(h) \epsilon\otimes e)= \mathcal{L}_{\pi,\sigma}(\epsilon'\otimes e', \epsilon\otimes e)$ pour tout $h\in H(F)$, il suffit de montrer le résultat pour $\tilde{y}=\theta_m$. D'après la proposition 4.0.1, il suffit de montrer que

$$L_{\pi,\sigma}(\tilde{\sigma}(\theta_m)\epsilon,\tilde{\pi}(\theta_m)e)=\overline{\epsilon_\nu(\tilde{\pi}^\vee,\tilde{\sigma})} L_{\pi,\sigma}(\epsilon,e)$$

\noindent pour tous $e\in E_\pi$ et $\epsilon\in E_\sigma$. Soit $\gamma\in G(F)$ tel que $\theta_m=\theta_d \gamma$. On a alors

\[\begin{aligned}
\displaystyle L_{\pi,\sigma}(\tilde{\sigma}(\theta_m)\epsilon,\tilde{\pi}(\theta_m)e)=\int_{\overline{U}_{m,m+1,m+r}(F)}\int_{U_m(F)\backslash H(F)} \overline{W_\epsilon}(J_m{}^t(h^c)^{-1}J_m^{-1}) & W_e(J_d{}^t(h^c\overline{u}^c)^{-1} J_d^{-1}\gamma) \\
 & |det(h)|_E^{-r} dh d\overline{u}
\end{aligned}\]

\noindent Après les changements de variables $h\mapsto h^c$ et $\overline{u}\mapsto \overline{u}^c$, on retrouve une expression qui est calculée en 5.5 [W3]. Le résultat est que

$$L_{\pi,\sigma}(\tilde{\sigma}(\theta_m)\epsilon,\tilde{\pi}(\theta_m)e)=\overline{w(\tilde{\sigma},\psi)}w(\tilde{\pi},\psi) \overline{\omega_\sigma((-1)^{d-1}z)}\omega_\pi(z') \epsilon(\pi\times \overline{\sigma},\psi) L_{\pi,\sigma}(\epsilon,e)$$

\noindent où $z=(-1)^{r+[(m+1)/2]}2\nu$ et $z'=(-1)^{r+1+[(d+1)/2]}2\nu$. Pour obtenir le résultat, il reste à vérifier que

$$\mbox{(1)}\;\; \overline{w(\tilde{\sigma},\psi)}w(\tilde{\pi},\psi) \overline{\omega_\sigma((-1)^{d-1}z)}\omega_\pi(z') \epsilon(\pi\times \overline{\sigma},\psi)=\overline{\epsilon_\nu(\tilde{\pi}^\vee,\tilde{\sigma})}$$

Puisque $\pi$ est tempérée donc unitaire, on a $\overline{\pi}^\vee\simeq \pi$. D'après le 1. de la proposition 2.5.1, on a donc

$$\overline{\epsilon(\pi^\vee\times \sigma,\psi)}=\epsilon(\overline{\pi}^\vee\times \overline{\sigma},\psi^-)=\omega_\pi(-1)^m \omega_\sigma(-1)^d\epsilon(\pi\times \overline{\sigma},\psi)$$

Puisque $\tilde{\pi}$ est tempérée donc unitaire, on a $\overline{\tilde{\pi}}^\vee\simeq \tilde{\pi}$. D'après 2.2(1), on a donc

$$\overline{w(\tilde{\pi}^\vee,\psi)}=w(\overline{\tilde{\pi}}^\vee,\psi^-)=\omega_\pi(-1)^{d-1}w(\tilde{\pi},\psi)$$

\noindent Comme $\pi\simeq \pi^\vee\circ c$, le caractère $\omega_\pi$ est trivial sur $N(E^\times)$, donc $\omega_{\pi| F^\times}^2=1$. Alors l'égalité (1) est équivalente à

$$\omega_\sigma((-1)^{2d-1+r+[(m+1)/2]}2\nu)\omega_\pi((-1)^{m+r+d+[(d+1)/2]}2\nu)=\omega_\sigma((-1)^{1+[d/2]}2\nu)\omega_\pi((-1)^{[m/2]}2\nu)$$

\noindent identité qu'il est facile de vérifier $\blacksquare$

\section{Le développement spectral}

On conserve la situation et les notations fixées en 3.1. On fixe une base $(v_i)_{i=1,\ldots,d}$ de $V$ qui est aussi une base de $R$ comme $\mathcal{O}_E$-réseau et on prend pour Levi tordu minimal $\tilde{M}_{min}=\tilde{T}_d$.

\subsection{La formule}

Posons

$$\displaystyle J_{spec}(\Theta_{\tilde{\sigma}},\tilde{f})=\sum_{\tilde{L}\in\mathcal{L}^{\tilde{G}}} (-1)^{a_{\tilde{L}}}|W^L||W^G|^{-1}\sum_{\mathcal{O} \in\{\Pi_{ell}(\tilde{L})\}} [i\mathcal{A}_{\mathcal{O}}^\vee,i\mathcal{A}_{\tilde{L},F}^\vee]^{-1}$$
$$\displaystyle 2^{-s(\mathcal{O})-a_{\tilde{L}}} \epsilon_\nu(\tilde{\rho}^\vee,\tilde{\sigma}) \int_{i\mathcal{A}_{\tilde{L},F}^*} J^{\tilde{G}}_{\tilde{L}}(\tilde{\rho}_\lambda,\tilde{f}) d\lambda$$

\noindent où pour toute orbite $\mathcal{O}\in \{\Pi_{ell}(\tilde{L})\}$, on a fixé un point base $\tilde{\rho}\in\mathcal{O}$.

\begin{theo}
On a

$$\displaystyle \lim\limits_{N\to\infty} J_N(\Theta_{\tilde{\sigma}},\tilde{f})=J_{spec}(\Theta_{\tilde{\sigma}},\tilde{f})$$
\end{theo}

\subsection{Utilisation de la formule de Plancherel}

Fixons $\tilde{y}\in\tilde{H}(F)$. On définit $f\in C_c^\infty(G(F))$ par $f(g)=\tilde{f}(g\tilde{y})$. Soit $g\in G(F)$. Fixons un sous-groupe ouvert et distingué $K_f$ de $K$ tel que $f$ soit biinvariante par $K_f$ et un sous-groupe compact ouvert $K'_g\subset H(F)$ tel que $gK'_gg^{-1}\subset K_f$ et $\theta_{\tilde{y}}(g)K'_g\theta_{\tilde{y}}(g)^{-1}\subset K_f$. Fixons aussi une base orthonormée $\mathcal{B}_\sigma^{K'_g}$ de l'espace des invariants $E_{\sigma}^{K'_g}$. On a alors

$$\displaystyle J(\Theta_{\tilde{\sigma}},\tilde{f},g)=\int_{H(F)}\Theta_{\tilde{\sigma}}(h\tilde{y}) \int_{U(F)} f(g^{-1} hu\theta_{\tilde{y}}(g)) \xi(u) du dh$$
$$\displaystyle =\sum_{\epsilon\in \mathcal{B}_\sigma^{K'_g}}\int_{H(F)U(F)} (\epsilon,\tilde{\sigma}(h\tilde{y})\epsilon) f(g^{-1} hu\theta_{\tilde{y}}(g)) \xi(u) du dh$$

\noindent La formule de Plancherel telle que rappelée en 1.7 [B], appliquée à $f$, s'écrit

$$\displaystyle f(g)=\sum_{L\in\mathcal{L}(A_d)} |W^L||W^G|^{-1}\sum_{\mathcal{O}\in\{\Pi_2(L)\}} f_{\mathcal{O}}(g)$$

\noindent où $f_{\mathcal{O}}\in\mathcal{S}(G(F))$ est définie par

$$\displaystyle f_{\mathcal{O}}(g)=[i\mathcal{A}_{\mathcal{O}}^\vee: i\mathcal{A}_{L,F}]^{-1}\int_{i\mathcal{A}_{L,F}^*} m(\tau_\lambda) Tr(i_Q^G(\tau_\lambda,g^{-1})i_Q^G(\tau_\lambda,f)) d\lambda$$

\noindent (on a comme toujours fixé un point base $\tau\in\mathcal{O}$ et un parabolique $Q\in\mathcal{P}(L)$). On en déduit que

$$\displaystyle \mbox{(1)} \;\; J(\Theta_{\tilde{\sigma}},\tilde{f},g)=\sum_{L\in\mathcal{L}(A_d)}|W^L||W^G|^{-1}\sum_{\mathcal{O}\in\{\Pi_2(L)\}} \sum_{\epsilon\in\mathcal{B}_{\sigma}^{K'_g}}J_{L,\mathcal{O}}(\epsilon,f,g)$$

\noindent où on a posé

$$\displaystyle J_{L,\mathcal{O}}(\epsilon,f,g)=\int_{H(F)U(F)} (\epsilon,\tilde{\sigma}(h\tilde{y})\epsilon) f_{\mathcal{O}}(g^{-1}hu\theta_{\tilde{y}}(g)) \xi(u) du dh$$

\noindent On justifie la manipulation formelle ayant permis d'obtenir (1) par l'absolue convergence de $J_{L,\mathcal{O}}(\epsilon,f,g)$: puisque $f_{\mathcal{O}}\in\mathcal{S}(G(F))$, il suffit de montrer la convergence de l'intégrale

$$\displaystyle \int_{H(F)U(F)} \Xi^H(h)\Xi^G(hu) \sigma(hu)^{-d} du dh$$

\noindent pour un entier $d$ assez grand. Cela découle de la proposition II.4.5 ce [W5] et de 4.1(3) [W3]. Fixons $L\in\mathcal{L}(A_d)$, $\mathcal{O}\in\{\Pi_2(L)\}$. \\

 (2) Il existe un entier $c_0\geqslant 0$ tel que l'on ait l'égalité
 
 $$\displaystyle \int_{U(F)}f_{\mathcal{O}}(g^{-1}hu\theta_{\tilde{y}}(g)) \xi(u) du=\int_{U(F)_c}f_{\mathcal{O}}(g^{-1}hu\theta_{\tilde{y}}(g)) \xi(u) du$$
 
\noindent pour tout $c\geqslant c_0$, pour tout $h\in H(F)$ et pour tout $g\in M(F)K$. \\

La fonction $f_{\mathcal{O}}$ est biinvariante par $K_f$. Posons $K_A=A(F)\cap K_f\cap \theta_{\tilde{y}}(K_f)$. Puisque $K_f$ est distingué dans $K$, que $A$ commute à $M$ et que $\theta_{\tilde{y}}(M)=M$, on a pour tout $a\in K_A$ l'égalité

$$\displaystyle \int_{U(F)}f_{\mathcal{O}}(g^{-1}hu\theta_{\tilde{y}}(g)) \xi(u) du=\int_{U(F)}f_{\mathcal{O}}(g^{-1}haua^{-1}\theta_{\tilde{y}}(g)) \xi(u) du$$

\noindent D'où

$$\displaystyle \int_{U(F)}f_{\mathcal{O}}(g^{-1}hu\theta_{\tilde{y}}(g)) \xi(u) du=\frac{1}{mes(K_A)}\int_{U(F)}f_{\mathcal{O}}(g^{-1}hu\theta_{\tilde{y}}(g))\int_{K_A} \xi(aua^{-1})da du$$

\noindent Or, il existe un entier $c_0\geqslant 0$ tel que pour $u\in U(F)-U(F)_{c_0}$, on ait

$$\displaystyle \int_{K_A} \xi(aua^{-1})da=0$$

\noindent On en déduit (2). \\

 Fixons une base orthonormée $\mathcal{B}_{\mathcal{O}}^{K_f}$ de $(\mathcal{K}_Q^G)^{K_f}$. On a l'égalité

$$\displaystyle f_{\mathcal{O}}(g')=[i\mathcal{A}_{\mathcal{O}}^\vee:i\mathcal{A}_{L,F}^\vee]^{-1}\sum_{e\in \mathcal{B}_{\mathcal{O}}^{K_f}} \int_{i\mathcal{A}_{L,F}^*} m(\tau_\lambda) (\pi_\lambda(g')e,\pi_\lambda(f)e) d\lambda$$

\noindent pour tout $g'\in G(F)$, où on a posé $\pi_\lambda=i_Q^G(\tau_\lambda)$. De (2) et de l'égalité précédente, l'on déduit

\[\begin{aligned}
\displaystyle J_{L,\mathcal{O}}(\epsilon,f,g)=[i\mathcal{A}_{\mathcal{O}}^\vee:i\mathcal{A}_{L,F}^\vee]^{-1} \sum_{e\in\mathcal{B}_{\mathcal{O}}^{K_f}} & \int_{H(F)U(F)_c} (\epsilon,\tilde{\sigma}(h\tilde{y}) \epsilon) \\
 & \int_{i\mathcal{A}_{L,F}^*} m(\tau_\lambda) \left(\pi_\lambda(hu\theta_{\tilde{y}}(g))e, \pi_\lambda(g)\pi_\lambda(f)e \right) d\lambda \xi(u) du dh
\end{aligned}\]

\noindent pour tout $c\geqslant c_0$, pour tout $\epsilon\in \mathcal{B}_\sigma^{K'_g}$ et pour tout $g\in M(F)K$. A $g\in M(F)K$ fixé, on a une majoration  $|(\pi_\lambda(hu\theta_{\tilde{y}}(g))e, \pi_\lambda(g)\pi_\lambda(f)e)|<<\Xi^G(hu)$ pour tous $\lambda$, $h$, $u$, $e\in\mathcal{B}_{\mathcal{O}}^{K_f}$. D'après 4.1(3) [W3], l'expression ci-dessus est absolument convergente. Cela justifie la permutation de la somme et de l'intégrale et permet aussi d'intervertir les deux intégrales intérieures. Après les changements de variables $h\mapsto h^{-1}$ et $u\mapsto hu^{-1}h^{-1}$, on obtient

\[\begin{aligned}
\displaystyle J_{L,\mathcal{O}}(\epsilon,f,g)=[i\mathcal{A}_{\mathcal{O}}^\vee:i\mathcal{A}_{L,F}^\vee]^{-1} \sum_{e\in\mathcal{B}_{\mathcal{O}}^{K_f}} \int_{i\mathcal{A}_{L,F}^*} & m(\tau_\lambda) \int_{H(F)U(F)_c} (\sigma(h)\epsilon,\tilde{\sigma}(\tilde{y}) \epsilon) \\
 & \left(\pi_\lambda(\theta_{\tilde{y}}(g))e, \pi_\lambda(hug)\pi_\lambda(f)e \right) \overline{\xi}(u) du dh d\lambda
\end{aligned}\]

\noindent On reconnaît l'intégrale intérieure: c'est $\mathcal{L}_{\pi_\lambda,\sigma,c}\left(\epsilon\otimes \pi_\lambda(\theta_{\tilde{y}}(g))e, \tilde{\sigma}(\tilde{y}) \epsilon\otimes \pi_\lambda(g)\pi_\lambda(f)e \right)$. Quitte à augmenter $c_0$, c'est aussi $\mathcal{L}_{\pi_\lambda,\sigma} \left(\epsilon\otimes \pi_\lambda(\theta_{\tilde{y}}(g))e, \tilde{\sigma}(\tilde{y}) \epsilon\otimes \pi_\lambda(g)\pi_\lambda(f)e \right)$. Fixons $\lambda\in i\mathcal{A}_{L,F}^*$ en position générale. Alors la représentation $\pi_\lambda$ est irréductible et on peut trouver un homomorphisme $\ell_\lambda\in Hom_{H,\xi}(\pi_\lambda,\sigma)$ et une constante $C_\lambda\in\mathbb{C}$ de sorte que l'on ait l'égalité

$$\mathcal{L}_{\pi_\lambda,\sigma}(\epsilon_1\otimes e_1,\epsilon_2\otimes e_2)=C_\lambda (\epsilon_1,\ell_\lambda(e_2))(\ell_\lambda(e_1),\epsilon_2)$$

\noindent pour tous $\epsilon_1,\epsilon_2\in E_\sigma$, $e_1,e_2\in\mathcal{K}_{Q,\tau}^G$. On a alors

$$\displaystyle \sum_{\epsilon\in\mathcal{B}_\sigma^{K'_g}} \mathcal{L}_{\pi_\lambda,\sigma} \left(\epsilon\otimes \pi_\lambda(\theta_{\tilde{y}}(g))e, \tilde{\sigma}(\tilde{y}) \epsilon\otimes \pi_\lambda(g)\pi_\lambda(f)e \right)$$
$$\displaystyle =\sum_{\epsilon\in\mathcal{B}_\sigma^{K'_g}} C_\lambda \left(\epsilon,\ell_\lambda(\pi_\lambda(g)\pi_\lambda(f)e) \right) \left(\ell_\lambda(\pi_\lambda(\theta_{\tilde{y}}(g))e),\tilde{\sigma}(\tilde{y})\epsilon \right)$$

\noindent Puisque $\ell_\lambda(\pi_\lambda(g)\pi_\lambda(f)e)$ est invariante par $K'_g$, on a

$$\displaystyle \ell_\lambda(\pi_\lambda(g)\pi_\lambda(f)e)= \sum_{\epsilon\in\mathcal{B}_\sigma^{K'_g}} \left( \epsilon,\ell_\lambda(\pi_\lambda(g)\pi_\lambda(f)e) \right) \epsilon$$

\noindent On en déduit que

$$\mbox{(3)}\;\; \displaystyle \sum_{\epsilon\in\mathcal{B}_\sigma^{K'_g}} \mathcal{L}_{\pi_\lambda,\sigma} \left( \epsilon\otimes \pi_\lambda(\theta_{\tilde{y}}(g))e, \tilde{\sigma}(\tilde{y}) \epsilon\otimes \pi_\lambda(g)\pi_\lambda(f)e \right)= C_\lambda \left( \ell_\lambda(\pi_\lambda(\theta_{\tilde{y}}(g))e),\tilde{\sigma}(\tilde{y})\ell_\lambda(\pi_\lambda(g)\pi_\lambda(f)e) \right)$$

\noindent Fixons des familles $(\varphi_j)_{j=1,\ldots,n}$, $(e_j)_{j=1,\ldots,n}$ et $(\epsilon_j)_{j=1,\ldots,n}$ vérifiant le 2. du lemme 4.0.1. On a alors

$$\displaystyle \sum_{j=1}^n C_\lambda \varphi_j(\lambda) (\epsilon_j,\ell_\lambda(e_j))(\ell_\lambda(e_j),\epsilon_j)=1$$

\noindent Multipliant cette égalité par (3), on obtient

$$\displaystyle \sum_{\epsilon\in\mathcal{B}_\sigma^{K'_g}} \mathcal{L}_{\pi_\lambda,\sigma} \left( \epsilon\otimes \pi_\lambda(\theta_{\tilde{y}}(g))e, \tilde{\sigma}(\tilde{y}) \epsilon\otimes \pi_\lambda(g)\pi_\lambda(f)e \right)$$
$$\displaystyle =\sum_{j=1}^n C_\lambda^2 \varphi_j(\lambda) (\epsilon_j,\ell_\lambda(e_j))(\ell_\lambda(e_j),\epsilon_j)\left( \ell_\lambda(\pi_\lambda(\theta_{\tilde{y}}(g))e),\tilde{\sigma}(\tilde{y})\ell_\lambda(\pi_\lambda(g)\pi_\lambda(f)e) \right)$$

Le produit d'un des facteurs $C_\lambda$, du premier et du troisième produits scalaires est égal à

$$\displaystyle \mathcal{L}_{\pi_\lambda,\sigma}(\epsilon_j\otimes \pi_\lambda(\theta_{\tilde{y}}(g))e,\tilde{\sigma}(\tilde{y})\ell_\lambda(\pi_\lambda(g)\pi_\lambda(f)e)\otimes e_j)$$
$$=\displaystyle \int_{H(F)U(F)_c} (\sigma(h)\epsilon_j,\tilde{\sigma}(\tilde{y})\ell_\lambda(\pi_\lambda(g)\pi_\lambda(f)e))(\pi_\lambda(\theta_{\tilde{y}}(g))e,\pi_\lambda(hu)e_j) \overline{\xi}(u) du dh$$

\noindent pour tout entier $c\geqslant 0$ assez grand. Le produit du deuxième facteur $C_\lambda$, de $(\ell_\lambda(e_j),\epsilon_j)$ et de $(\sigma(h)\epsilon_j,\tilde{\sigma}(\tilde{y})\ell_\lambda(\pi_\lambda(g)\pi_\lambda(f)e))$ est égal à

$$\mathcal{L}_{\pi_\lambda,\sigma}(\tilde{\sigma}(h^{-1}\tilde{y})^{-1}\epsilon_j\otimes e_j, \epsilon_j\otimes \pi_\lambda(g)\pi_\lambda(f)e)$$
$$\displaystyle =\int_{H(F)U(F)_{c'}} (\sigma(\theta_{\tilde{y}}(h')h)\epsilon_j, \tilde{\sigma}(\tilde{y})\epsilon_j) (e_j,\pi_\lambda(h'u'g) \pi_\lambda(f)e) \overline{\xi}(u') du'dh'$$

\noindent pour tout entier $c'\geqslant 0$ assez grand. Posons

\[\begin{aligned}
\displaystyle X_{\lambda,j,c,c'}(e,g)=\int_{H(F)U(F)_c}\int_{H(F)U(F)_{c'}} & (\sigma(h)\epsilon_j,\tilde{\sigma}(\tilde{y})\epsilon_j) \left(\pi_\lambda(\theta_{\tilde{y}}(h'u'g))e,\pi_\lambda(hu)e_j \right) \\
 & \left( e_j,\pi_\lambda(h'u'g)\pi_\lambda(f)e \right) \overline{\xi}(u) du'dh'dudh
\end{aligned}\]

\noindent Cette expression est absolument convergente d'après [W3] 4.1(4). Alors d'après ce qui précède et après les changements de variables $h\mapsto \theta_{\tilde{y}}(h')^{-1}h$ et $u\mapsto h^{-1}\theta_{\tilde{y}}(h')\theta_{\tilde{y}}(u')^{-1} \theta_{\tilde{y}}(h')^{-1}hu$, on obtient l'égalité

$$\mbox{(4)}\;\;\; \displaystyle \sum_{\epsilon\in\mathcal{B}_\sigma^{K'_g}} \mathcal{L}_{\pi_\lambda,\sigma} \left( \epsilon\otimes \pi_\lambda(\theta_{\tilde{y}}(g))e, \tilde{\sigma}(\tilde{y}) \epsilon\otimes \pi_\lambda(g)\pi_\lambda(f)e \right)$$
$$\displaystyle =\sum_{j=1}^n \varphi_j(\lambda) X_{\lambda,j,c,c'}(e,g)$$

\noindent pour tous entiers $c,c'\geqslant 0$ assez grands. De la même façon que l'on a montré (2), on montre qu'il existe un entier $c_0\geqslant 0$ tel que l'égalité précédente soit vérifiée pour tout $g\in M(F)K$, pour tout $\lambda\in i\mathcal{A}_{L,F}^*$, pour tout $e\in \mathcal{B}_{\mathcal{O}}^{K_f}$ et pour tous entiers $c,c'\geqslant c_0$. On fixe un tel entier $c_0$. Posons

$$\displaystyle J_{L,\mathcal{O}}(\Theta_{\tilde{\sigma}},\tilde{f},g)=[i\mathcal{A}_{\mathcal{O}}^\vee:i\mathcal{A}_{L,F}^\vee]^{-1}\sum_{e\in\mathcal{B}_{\mathcal{O}}^{K_f}} \sum_{j=1}^n \int_{i\mathcal{A}_{L,F}^*} m(\tau_\lambda) \varphi_j(\lambda) X_{\lambda,j,c,c'}(e,g) d\lambda$$

\noindent pour tout $g\in M(F)K$ et pour tous $c,c'\geqslant c_0$. On vient donc de montrer

$$\displaystyle \mbox{(5)}\;\;\; \sum_{\epsilon\in \mathcal{B}_{\sigma}^{K'_g}} J_{L,\mathcal{O}}(\epsilon,f,g)=J_{L,\mathcal{O}}(\Theta_{\tilde{\sigma}},\tilde{f},g)$$

\noindent pour tout $g\in M(F)K$. Soit $C>0$. Pour $N\geqslant 1$ un entier, posons

\[\begin{aligned}
\displaystyle J_{L,\mathcal{O},N,C}(\Theta_{\tilde{\sigma}},\tilde{f})& = [i\mathcal{A}_\mathcal{O}^\vee:i\mathcal{A}_{L,F}^\vee]^{-1} \sum_{e\in\mathcal{B}_\mathcal{O}^{K_f}}\sum_{j=1}^n \int_{i\mathcal{A}_{L,F}^*} m(\tau_\lambda) \phi_j(\lambda) \int_{H(F)U(F)_c} \mathbf{1}_{\sigma<Clog(N)}(hu) \\
 & (\sigma(h)\epsilon_j,\tilde{\sigma}(\tilde{y})\epsilon_j) \overline{\xi}(u)\int_{G(F)}(\pi_\lambda(\theta_{\tilde{y}}(g))e,\pi_\lambda(hu)e_j)(e_j,\pi_\lambda(g)\pi_\lambda(f)e) \kappa_N(g) dgdudhd\lambda
\end{aligned}\]

\noindent pour tout $c\geqslant c_0$. D'après (1) et (5), on a

$$\displaystyle J(\Theta_{\tilde{\sigma}},\tilde{f},g)=\sum_{L\in\mathcal{L}(A_d)}|W^L||W^G|^{-1} \sum_{\mathcal{O}\in\{\Pi_2(L)\}} J_{L,\mathcal{O}}(\Theta_{\tilde{\sigma}},\tilde{f},g)$$

\noindent pour tout $g\in M(F)K$. Par définition, $J_N(\Theta_{\tilde{\sigma}},\tilde{f})$ est l'intégrale de $J(\Theta_{\tilde{\sigma}},\tilde{f},g)\kappa_N(g)$ sur $H(F)U(F)\backslash G(F)$, ou ce qui revient au même l'intégrale de $J(\Theta_{\tilde{\sigma}},\tilde{f},mk)\kappa_N(mk)\delta_P(m)^{-1}$ sur $H(F)\backslash M(F)\times K$.

\begin{prop}
\begin{enumerate}
\item L'expression $J_{L,\mathcal{O},N,C}(\Theta_{\tilde{\sigma}},\tilde{f})$ est absolument convergente.
\item Il existe $C>0$ tel que

$$|J_N(\Theta_{\tilde{\sigma}},\tilde{f})-\sum_{L\in\mathcal{L}(A_d)}|W^L||W^G|^{-1} \sum_{\mathcal{O}\in\{\Pi_2(L)\}} J_{L,\mathcal{O},N,C}(\Theta_{\tilde{\sigma}},\tilde{f})|<<N^{-1}$$
\noindent pour tout $N\geqslant 1$.
\end{enumerate}
\end{prop}

\ul{Preuve}: C'est la même que celle du lemme 6.4 de [W2], en utilisant les majorations 4.1(5), (6) et (7) de [W3]. \\

Dans la suite, on fixe un tel réel $C$, un Levi $L\in\mathcal{L}(A_d)$ et une orbite $\mathcal{O}\in\{\Pi_2(L)\}$.

\subsection{Changement de fonction de troncature}

Notons $\Delta$ l'ensemble des racines simples de $A_d=A_{T_d}$ dans $\mathfrak{u}_d$ et $\tilde{\Delta}$ les restrictions des éléments de $\Delta$ à $A_{\tilde{T}_d}$. L'ensemble $\tilde{\Delta}$ est en bijection avec l'ensemble des orbites de $\theta_d$ dans $\Delta$. Notons $\alpha\mapsto (\alpha)$ cette bijection. Pour $\alpha\in \tilde{\Delta}$, on définit un poids $\varpi_{\alpha^\vee}$ par

$$\displaystyle \varpi_{\alpha^\vee}=\frac{1}{n(\alpha)}\sum_{\beta\in (\alpha)} \varpi_{\beta^\vee}$$

\noindent où $n(\alpha)$ désigne le cardinal de l'orbite $(\alpha)$ et les $\varpi_{\beta^\vee}\in \mathcal{A}_{A_d}^*$ sont les poids classiques. Fixons $\delta>0$. Introduisons l'ensemble $\mathcal{D}$ constitué des $Y\in \mathcal{A}_{A_d,F}\otimes \mathbb{Q}$ vérifiant

\begin{itemize}
\item $\alpha(Y)\geqslant 0$ pour tout $\alpha\in \Delta$;
\item $inf\{\alpha(Y);\; \alpha\in\Delta\}\geqslant \delta sup\{\alpha(Y);\; \alpha\in\Delta\}$.
\end{itemize}

Soit $Y\in\mathcal{D}$. On note $\varphi^+(.,Y)$ la fonction caractéristique du sous-ensemble des $\zeta\in \mathcal{A}_{A_d}$ vérifiant

\begin{itemize}
\item $\alpha(\zeta)\geqslant 0$ pour tout $\alpha\in \Delta$;
\item $\varpi_{\alpha^\vee}(Y-\zeta)\geqslant 0$ pour tout $\alpha\in\tilde{\Delta}$.
\end{itemize}

\noindent C'est un sous-ensemble compact modulo $\mathcal{A}_G$. Définissons $g\mapsto \tilde{u}(g,Y)$ comme la fonction caractéristique du sous-ensemble des élément $g$ de $G(F)$ qui s'écrivent $g=k_1ak_2$ avec $k_1,k_2\in K$ et $a\in T_d(F)$ tel que $\varphi^+(H_{T_d}(a),Y)=1$. Soit $f'\in C_c^\infty(G(F))$ définie par $f'(g)=\tilde{f}(g\theta_d)$ pour tout $g\in G(F)$. Fixons un sous-groupe ouvert-compact $K_{f'}$ de $G(F)$ tel que $f'$ soit biinvariante par $K_{f'}$ et une base orthonormée $\mathcal{B}_{\mathcal{O}}^{K_{f'}}$ de $\mathcal{K}^G_{Q,\tau}$. Fixons aussi une fonction $\varphi\in C^\infty(i\mathcal{A}_{L,F}^*)$ et deux éléments $e',e''\in \mathcal{K}^G_{Q,\tau}$. Posons

$$\displaystyle \Phi(e,g,g',\lambda)=(\pi_\lambda(\theta_d(g))e,\pi_\lambda(g')e')(e'',\pi_\lambda(g)\pi_\lambda(f') e)$$

\noindent pour tout $e\in\mathcal{B}_{\mathcal{O}}^{K_{f'}}$, pour tous $g,g'\in G(F)$ et pour tout $\lambda\in \mathcal{A}_{L,F}^*$. Posons aussi

$$\displaystyle \Phi_N(g')=\sum_{e\in \mathcal{B}_{\mathcal{O}}^{K_{f'}}} \int_{G(F)}\int_{i\mathcal{A}_{L,F}^*} m(\tau_\lambda)\varphi(\lambda)\Phi(e,g,g',\lambda)d\lambda\kappa_N(g) dg$$

$$\displaystyle \Phi_Y(g')=\sum_{e\in \mathcal{B}_{\mathcal{O}}^{K_{f'}}} \int_{G(F)}\int_{i\mathcal{A}_{L,F}^*} m(\tau_\lambda)\varphi(\lambda)\Phi(e,g,g',\lambda)d\lambda\tilde{u}(g,Y) dg$$

\noindent pour tout $g'\in G(F)$, pour tout entier $N\geqslant 1$ et pour tout $Y\in\mathcal{D}$.

\begin{prop}
\begin{enumerate}
\item Les expressions $\Phi_N(g')$ et $\Phi_Y(g')$ convergent dans l'ordre d'intégration où elles ont été définies.
\item Pour tout $0<\eta<1$, il existe $c_1>0$ et $c_2>0$ tels que

$$|\Phi_N(g')-\Phi_Y(g')|<<N^{-1}$$

\noindent pour tout entier $N\geqslant 1$, pour tout $g'\in G(F)$ vérifiant $\sigma(g')\leqslant Clog(N)$ et pour tout $Y\in\mathcal{D}$ tel que

$$c_1N^\eta\leqslant sup\{\alpha(Y);\; \alpha\in\Delta\}\leqslant c_2N$$
\end{enumerate}
\end{prop}

\ul{Preuve}: Il suffit de reprendre celle de la proposition 6.5 de [W3]. En effet, à part des considérations valables pour tous les groupes tordus, la démonstration de [W3] n'utilise que les faits particuliers suivants:

\begin{itemize}
\item La majoration 4.1(1) de [W3] qui est toujours valable ici car $G(F)$ est le groupe des $E$-points d'un groupe linéaire.
\item L'existence d'un réel $c>0$ tel que $\sigma(g)<cN$ entraîne $\kappa_N(g)=1$ pour tout $g\in G(F)$.
\item L'injectivité de l'application $\mathcal{A}_{G}\to\mathcal{A}_{G}$, $X\mapsto X-\theta(X)$
\end{itemize}

$\blacksquare$

\subsection{Une approximation de $\Phi_Y(g')$}

Soit $\tilde{L}'\in\mathcal{L}^{\tilde{G}}$ tel que $L\subset L'$ et fixons $\tilde{S}'=\tilde{L}'U'\in \mathcal{P}(\tilde{L}')$. Notons $\Lambda_{\mathcal{O},ell}^{\tilde{L}'}$ l'ensemble des éléments $\lambda\in i\mathcal{A}_L^*$ tel que $i^{L'}_{Q\cap L'}(\tau_\lambda)$ s'étende en une représentation elliptique de $\tilde{L}'(F)$. Pour un tel élément, on fixe un prolongement unitaire $\tilde{\rho}_\lambda$ de $i^{L'}_{Q\cap L'}(\tau_\lambda)$ à $\tilde{L}'(F)$. Puisque $\pi_\lambda=i^G_Q(\tau_\lambda)$ s'identifie à l'induite normalisée de $S'$ à $G$ de $\rho_\lambda$, l'on déduit de ce prolongement un prolongement $\tilde{\pi}_\lambda$ de $\pi_\lambda$ à $\tilde{G}(F)$. L'ensemble $\Lambda_{\mathcal{O},ell}^{\tilde{L}'}$ est stable par translation par $i\mathcal{A}_{L,F}^\vee+i\mathcal{A}_{\tilde{L}'}^*$. On suppose que nos choix ont été faits de sorte que $\tilde{\rho}_{\lambda+\chi}=\tilde{\rho}_\lambda$ pour tout $\chi\in i\mathcal{A}_{L,F}^\vee$ et $\tilde{\rho}_{\lambda+\chi}=(\tilde{\rho}_\lambda)_\chi$ pour tout $\chi\in i\mathcal{A}_{\tilde{L}'}^*$. Posons

$$\displaystyle \Phi(g')=\sum_{\tilde{L}'\in\mathcal{L}^{\tilde{G}}; L\subset L'}(-1)^{a_{\tilde{L}'}} \sum_{\lambda\in \Lambda^{\tilde{L}'}_{\mathcal{O},ell}/(i\mathcal{A}_{L,F}^\vee+ i\mathcal{A}_{\tilde{L}'}^*)} 2^{a_{\tilde{L}'}-a_L}$$
$$\displaystyle \int_{i\mathcal{A}_{L',F}^*} (\tilde{\pi}_{\lambda+\chi}(\theta_d)e'',\pi_{\lambda+\chi}(g')e') J_{\tilde{L}'}^{\tilde{G}}(\tilde{\rho}_{\lambda+\chi},\tilde{f})\varphi(\lambda+\chi) d\chi$$

\begin{prop}
Pour tout réel $R\geqslant 1$, il existe un entier $k\geqslant 0$ tel que l'on ait une majoration

$$|\Phi_Y(g')-\Phi(g')|<<\sigma(g')^k\Xi^G(g') |Y|^{-R}$$

\noindent pour tout $g'\in G(F)$ et pour tout $Y\in\mathcal{D}$.
\end{prop}

\ul{Preuve}: Il suffit de reprendre celle de la proposition 6.7 et celle du lemme 6.8 de [W3]. La proposition 6.7 de [W3] est valable pour tout les groupes tordus, à la seule condition d'avoir à $g'$ fixé une majoration

$$|\Phi_Y(g')-\Phi_{Y'}(g')|<<|Y|^{-1}$$

\noindent pour tous $Y,Y'\in\mathcal{D}$ vérifiant $|Y-Y'|\geqslant |Y|/2$. Mais cela est une conséquence de la proposition 5.3.1. La démonstration du lemme 6.8 de [W3] n'utilise que les propriétés suivantes du groupe tordu $\tilde{G}$:

\begin{enumerate}
\item Toute induite parabolique normalisée d'une représentation irréductible de la série discrète est aussi irréductible.
\item Pour tout $\lambda\in \Lambda_{\mathcal{O},ell}^{\tilde{L}'}$, il existe un unique $\tilde{t}\in W^{\tilde{L}'}(L)_{reg}$ tel que $\theta_{\tilde{t}}(\tau_\lambda)\simeq \tau_\lambda$ et pour ce $\tilde{t}$, on a $|det(1-\theta_{\tilde{t}})_{\mathcal{A}_L/\mathcal{A}_{\tilde{L}'}}|=2^{a_L-a_{\tilde{L}'}}$
\end{enumerate}

Dans notre cas, la première propriété est bien connue et la deuxième découle de notre description en 2.4 des représentations elliptiques des Levi tordu de $\tilde{G}$ $\blacksquare$

\subsection{Preuve du théorème 5.1.1}

Il suffit cette fois de reprendre les sections 6.9 et 6.10 de [W3]. Dans la preuve du lemme 6.9, il faut utiliser la proposition 4.1.1 à la place de la proposition 5.5 de [W3] et la relation 2.5(2) au lieu de la relation 2.5(2) de [W3]. Enfin dans la section 6.10 de [W3], il est utilisé le fait suivant \\

 Pour tout Levi tordu $\tilde{L}'\in \mathcal{L}^{\tilde{G}}$, pour tout Levi $L\in \mathcal{L}^{L'}(A_d)$ et pour toute représentation $\tau\in\Pi_2(L)$ tels que $i_Q^{L'}(\tau)$ ($Q\in\mathcal{P}^{L'}(L)$) se prolonge en une représentation elliptique de $\tilde{L}(F)$, on a $W^{L'}(\tau)=\{1\}$. \\
 
 Cela se voit directement sur la description des représentations elliptiques données en 2.4.

\section{Calcul d'un facteur $\epsilon$ par une formule intégrale}

\subsection{Enoncé}

On conserve toujours la situation et les notations de 3.1. Soient $\Gamma$ et $\Theta$ des quasi-caractères de $\tilde{G}(F)$ et $\tilde{H}(F)$ respectivement. Posons

\[\begin{aligned}
\displaystyle & \epsilon_{geom,\nu}(\Gamma,\Theta)=\epsilon_{geom,-\nu}(\Theta,\Gamma)= \\ 
 & \sum_{\tilde{T}\in\mathcal{T}} |2|_F^{dim(W'_T)-\big(\frac{d+m}{2}\big)} |W(H,\tilde{T})|^{-1} \lim\limits_{s\to 0^+} \int_{\tilde{T}(F)/\theta} c_\Gamma(\tilde{t}) c_\Theta(\tilde{t}) D^{\tilde{H}}(\tilde{t})^{1/2}D^{\tilde{G}}(\tilde{t})^{1/2} \Delta(\tilde{t})^{s-1/2} d\tilde{t}
\end{aligned}\]

\noindent Tous les termes apparaissant dans cette définition ont été introduits en 2.1, 3.2 et 3.3. La dépendance en $\nu$ de l'expression précédente est cachée dans le terme $c_\Gamma(\tilde{t})$ qui est défini à partir du plongement de $\tilde{H}$ dans $\tilde{G}$ fixé en 3.1, plongement qui dépend de $\nu$. Cette expression ne dépend vraiment que de $\Gamma$, $\Theta$ et $\nu$, c'est-à-dire qu'il ne dépend pas du choix d'une base $(z_i)_{i=-r,\ldots,r}$ d'un supplémentaire de $W$ dans $V$ qui permet de construire un plongement $\tilde{H}\hookrightarrow \tilde{G}$. Pour $\tilde{\pi}\in Temp(\tilde{G})$ et $\tilde{\sigma}\in Temp(\tilde{H})$, on adoptera les notations suivantes: $c_{\tilde{\pi}}=c_{\Theta_{\tilde{\pi}}}$, $c_{\tilde{\sigma}}=c_{\Theta_{\tilde{\sigma}}}$ et $\epsilon_{geom,\nu}(\tilde{\pi},\tilde{\sigma})=\epsilon_{geom,\nu}(\Theta_{\tilde{\pi}},\Theta_{\tilde{\sigma}})$.

\begin{theo}
Pour tout $\tilde{\pi}\in Temp(\tilde{G})$ et pour tout $\tilde{\sigma}\in Temp(\tilde{H})$, on a

$$\epsilon_{geom,\nu}(\tilde{\pi},\tilde{\sigma})=\epsilon_\nu(\tilde{\pi},\tilde{\sigma})$$
\end{theo}

\subsection{Induction du terme géométrique}

Soit $\tilde{L}$ un Levi tordu de $\tilde{G}$. Il existe alors une décomposition

$$V=V_u\oplus\ldots\oplus V_1\oplus V_0\oplus V_{-1}\oplus \ldots\oplus V_{-u}$$

\noindent de sorte que $\tilde{L}$ soit l'ensemble des $\tilde{x}\in \tilde{G}$ vérifiant $\tilde{x}(V_j)=V_{-j}^*$ pour $j=-u,\ldots,u$. Pour $\tilde{x}\in \tilde{L}$ et $j=-u,\ldots,u$, on notera $\tilde{x}_j$ la restriction de $\tilde{x}$ à $V_j$. On note $G_0=R_{E/F} GL_E(V_0)$ et $\tilde{G}_0=Isom_c(V_0,V_0^*)$, les analogues de $G$ et $\tilde{G}$ lorsque l'on remplace $V$ par $V_0$. Soit $\Gamma_0$ un quasi-caractère de $\tilde{G}_0(F)$ et pour $j=1,\ldots,u$, $\Gamma_j$ un quasi-caractère de $GL_E(V_j)$. On définit un quasi-caractère $\Gamma^{\tilde{L}}$ de $\tilde{L}(F)$ par

$$\displaystyle \Gamma^{\tilde{L}}(\tilde{x})=\Gamma_0(\tilde{x}_0)\prod_{j=1}^u \Gamma_j({}^t(\tilde{x}_{-j}^c)^{-1} \tilde{x}_j)$$

\noindent pour tout $\tilde{x}\in \tilde{L}(F)$. On peut considérer le quasi-caractère induit $\Gamma=Ind_{\tilde{L}}^{\tilde{G}}(\Gamma^{\tilde{L}})$ (cf 1.12 [W3]). Pour $d\geqslant 1$ un entier, l'algèbre de Lie $\mathfrak{gl}_{d}(E)$ de $GL_d(E)$ possède une unique orbite nilpotente régulière $\mathcal{O}_{GL}$. Si $\Gamma^{GL}$ est un quasi-caractère de $GL_d(E)$, on pose

$$m_{geom}(\Gamma^{GL})=c_{\Gamma^{GL},\mathcal{O}_{GL}}(1)$$

\begin{prop}
On a l'égalité

$$\displaystyle \epsilon_{geom,\nu}(\Gamma,\Theta)=\epsilon_{geom,\nu}(\Gamma_0,\Theta)\prod_{j=1}^u m_{geom}(\Gamma_j)$$
\end{prop}

\ul{Preuve}: Supposons d'abord $d_0<m$ où $d_0=dim(V_0)$. Quitte à conjuguer $\tilde{L}$, on peut supposer que $V_0\subset W$. Soit $Z_0$ un supplémentaire de $V_0$ dans $W$, sa dimension est $2r_0+1$ pour un certain entier naturel $r_0$. Fixons une base $(z_{0,i})_{i=-r_0,\ldots,r_0}$ de $Z_0$ et définissons un isomorphisme $c$-linéaire $\tilde{\zeta}_0: Z_0\to Z_0^*$ par $\tilde{\zeta}_0(z_{0,i})=(-1)^{i+1}2\nu z_{0,-i}$ pour $i=-r_0,\ldots,r_0$. En appliquant les mêmes constructions qu'en 3.1, on définit un plongement $\iota_0:\tilde{G}_0\hookrightarrow \tilde{H}$ et un ensemble $\underline{\mathcal{T}}_0$ de sous-variétés $\tilde{T}_0\subset \tilde{G}_0$ qui sont des espaces homogènes principaux à droite et à gauche sous certains sous-tores $T_0$ de $G_0$. Posons $Z''=Z\oplus Z_0$ et $\tilde{\zeta}''=\tilde{\zeta}\oplus\tilde{\zeta}_0$. L'espace hermitien $(Z'',\tilde{\zeta}'')$ est hyperbolique. On peut donc fixer une base $(z''_i)_{i=\pm1,\ldots,\pm r''}$ de $Z''$ telle que $\tilde{\zeta}'' z''_i={z''}_{-i}^*$ pour $i=\pm1,\ldots,\pm r''$. Quitte à conjuguer à nouveau $\tilde{L}$, on peut supposer qu'il existe des sous-ensembles $I_j\subset\{1,\ldots,r''\}$, $j=1,\ldots,u$ tels que pour $j=1,\ldots,u$, $V_j$ soit engendré par les $z''_i$, $i\in I_j$ et $V_{-j}$ par les $z''_{-i}$, $i\in I_j$. La composée des inclusions $\tilde{G}_0\subset \tilde{H}\subset \tilde{G}$ envoie alors $\tilde{G}_0$ dans $\tilde{L}$. On a par définition

\[\begin{aligned}
\displaystyle \mbox{(1)}\;\; \epsilon_{geom,\nu}(\Gamma_0,\Theta)=\epsilon_{geom,-\nu}(\Theta,\Gamma_0) & =\sum_{\tilde{T}_0\in\mathcal{T}_0} |2|_F^{dim(V'_{0,T_0})-\big(\frac{d_0+m}{2}\big)} |W(G_0,\tilde{T}_0)|^{-1} \\
 & \lim\limits_{s\to 0^+} \int_{\tilde{T}_0(F)/\theta} c_{\Gamma_0}(\tilde{t}) c_{\Theta}(\tilde{t}) D^{\tilde{G}_0}(\tilde{t})^{1/2}D^{\tilde{H}}(\tilde{t})^{1/2} \Delta(\tilde{t})^{s-1/2} d\tilde{t}
\end{aligned}\]

\[\begin{aligned}
\displaystyle \mbox{(2)}\;\; \epsilon_{geom,\nu}(\Gamma,\Theta)=\sum_{\tilde{T}\in\mathcal{T}} |2|_F^{dim(W'_T)-\big(\frac{d+m}{2}\big)} |W(H,\tilde{T})|^{-1} \lim\limits_{s\to 0^+} \int_{\tilde{T}(F)/\theta} & c_{\Gamma}(\tilde{t}) c_{\Theta}(\tilde{t}) D^{\tilde{H}}(\tilde{t})^{1/2} \\
 & D^{\tilde{G}}(\tilde{t})^{1/2}\Delta(\tilde{t})^{s-1/2} d\tilde{t}
\end{aligned}\]

\noindent Puisque l'espace hermitien $(Z'',\tilde{\zeta}'')=(Z\oplus Z_0,\tilde{\zeta}\oplus \tilde{\zeta}_0)$ est hyperbolique, on a $\underline{\mathcal{T}}_0\subset \underline{\mathcal{T}}$. Montrons le fait suivant \\

 (3) Deux éléments de $\underline{\mathcal{T}}_0$ sont conjugués par $H(F)$ si et seulement si ils sont conjugués par $G_0(F)$. \\
 
 Soient $\tilde{T}_1,\tilde{T}_2\in \underline{\mathcal{T}}_0$ et supposons qu'il existe $h\in H(F)$ tel que $h\tilde{T}_1h^{-1}=\tilde{T}_2$. On a des décompositions $V_0=V'_{0,i}\oplus V''_{0,i}$ pour $i=1,2$ relatives à $\tilde{T}_1$ et $\tilde{T}_2$. On pose $W''_i=V''_{0,i}\oplus Z_0$, $i=1,2$. Pour un élément $\tilde{t}\in\tilde{T}_i$ en position générale, $W''_i$ est le noyau de $t-1$ dans $W$ où $t={}^t(\tilde{t}^c)^{-1} \tilde{t}$. Par conséquent, la restriction de $h$ à $W''_1$ réalise un isomorphisme entre les espaces hermitiens $(W''_1,\tilde{\zeta}_{H,T_1})$ et $(W''_2,\tilde{\zeta}_{H,T_2})$. Ces deux espaces hermitiens contiennent le sous-espace hermitien $(Z_0,\tilde{\zeta}_0)$. D'après le théorème de Witt, il existe $h'\in U(W''_2)(F)\subset Z_{H(F)}(\tilde{T}_2)$ tel que la restriction de $h'h$ à $Z_0$ soit l'identité. Quitte à remplacer $h$ par $h'h$, on peut donc supposer que $h$ agit trivialement sur $Z_0$. Mais alors $h$ envoie l'orthogonal de $Z_0$ dans $W''_1$, qui est $V''_{0,1}$, dans l'orthogonal de $Z_0$ dans $W''_2$, qui est $V''_{0,2}$. Pour $\tilde{t}\in\tilde{T}_i$ en position générale, $V'_{0,i}$ est l'unique supplémentaire $t$-stable de $W''_i$ dans $W$ où $t={}^t(\tilde{t}^c)^{-1} \tilde{t}$. Donc $h$ envoie $V'_{0,1}$ dans $V'_{0,2}$. Comme $V_0=V''_{0,i}\oplus V'_{0,i}$ pour $i=1,2$, on en déduit que $h$ laisse $V_0$ stable d'où $h\in G_0(F)$. \\
 
 D'après (3) on peut supposer $\mathcal{T}_0\subset \mathcal{T}$. Montrons \\
 
 (4) Pour $\tilde{T}\in\mathcal{T}-\mathcal{T}_0$, le terme indexé par $\tilde{T}$ dans (2) est nul. \\
 
 Introduisons la décomposition $W=W''\oplus W'$ relative à $\tilde{T}$ et posons $V''=W''\oplus Z$. Soit $\tilde{t}\in\tilde{T}(F)$ en position générale et montrons que $c_\Gamma(\tilde{t})=0$, ce qui bien évidemment entraîne (4). Par définition de l'induite d'un quasi-caractère, le support semi-simple de $\Gamma$ ne contient que des éléments dont la classe de conjugaison coupe $\tilde{L}(F)$. Il suffit donc de vérifier que $\tilde{t}$ n'est pas conjugué à un élément de $\tilde{L}(F)$. Supposons que ce soit le cas. Alors $G_{\tilde{t}}$ contient un conjugué de $A_{\tilde{L}}$. Puisque $\tilde{t}$ est supposé en position générale, on a $G_{\tilde{t}}=U(V'')\times T_\theta$. Le tore $T_\theta$ est anisotrope, donc $U(V'')$ contient comme sous-groupe $gA_{\tilde{L}}g^{-1}$ pour un certain $g\in G(F)$. Le sous-espace $g(V_1\oplus\ldots\oplus V_u)$ est alors inclus dans $V''$ et totalement isotrope pour la forme hermitienne $\tilde{\zeta}_{G,T}$. Puisque $dim(V_1\oplus\ldots\oplus V_u)=(d_V-d_0)/2$, l'espace $(V'',\tilde{\zeta}_{G,T})$ contient un sous-espace hermitien isomorphe à la somme directe orthogonale de $(d_V-d_0)/2$ plans hyperboliques c'est-à-dire isomorphe à $(Z'',\tilde{\zeta}'')$. D'après le théorème de simplification de Witt, $(W'',\tilde{\zeta}_{H,T})$ contient un sous-espace hermitien isomorphe à $(Z_0,\tilde{\zeta}_0)$. Fixons un tel sous-espace $W''_1$ et notons $W''_2$ son orthogonal dans $W''$. Il existe $h\in H(F)$ qui réalise une isométrie de $(W''_1,\tilde{\zeta}_{H,T})$ vers $(Z_0,\tilde{\zeta}_0)$ et qui envoie $W''_2\oplus W'$ sur $V_0$. Quitte à conjuguer $\tilde{T}$ par $h$, on s'est ramené au cas où $W''_1=(Z_0,\tilde{\zeta}_0)$ et $W''_2\oplus W'=V_0$. On a alors $\tilde{T}\subset \tilde{G_0}$. Puisque l'espace hermitien $V''=Z\oplus^\perp Z_0\oplus^\perp W''_2$ est quasi-déployé et que $Z\oplus^\perp Z_0$ est somme directe orthogonale de plans hyperboliques, l'espace hermitien $W''_2$ est aussi quasi-déployé. On en déduit que $\tilde{T}\in\underline{\mathcal{T}}_0$, ce qui contredit l'hypothèse. \\
 
 Fixons maintenant $\tilde{T}\in\mathcal{T}_0$. On note $V_0=V_0''\oplus V_0'$ et $\tilde{\zeta}_{G_0,T}$ la décomposition et la forme hermitienne sur $V_0''$ relatives à $\tilde{T}$. On posera $H'=R_{E/F} GL(V_0')$, $W''=Z_0\oplus V_0''$, $V''=Z\oplus W''$ et $\tilde{\zeta}_{H,T}=\tilde{\zeta}_0\oplus\tilde{\zeta}_{G_0,T}$, $\tilde{\zeta}_{G,T}=\tilde{\zeta}\oplus \tilde{\zeta}_{H,T}$ qui sont des 
 formes hermitiennes sur $W''$ et $V''$ respectivement. Etablissons \\
 
 (5) on a $|W(H,\tilde{T})|=|W(G_0,\tilde{T})|$. \\
 
 En effet, un élément de $Norm_{H(F)}(\tilde{T})$ réalise forcément par restriction une isométrie de $(W'',\tilde{\zeta}_{H,T})$ sur lui-même et laisse stable $V_0'$. Par conséquent $Norm_{H(F)}(\tilde{T})=U(\tilde{\zeta}_{H,T})\times Norm_{H'(F)}(\tilde{T})$. De la même façon $Norm_{G_0(F)}(\tilde{T})=U(\tilde{\zeta}_{G_0,T})\times Norm_{H'(F)}(\tilde{T})$. On en déduit des isomorphismes
 
\begin{center}
$W(H,\tilde{T})\simeq Norm_{H'(F)}(\tilde{T})/T(F)$ et $W(G_0,\tilde{T})\simeq Norm_{H'(F)}(\tilde{T})/T(F)$
\end{center}

\noindent d'où l'égalité (5). \\

 Soit $\tilde{t}\in\tilde{T}(F)$ en position générale. D'après (1),(2),(4) et (5), il suffit pour établir la proposition de vérifier l'égalité

$$\mbox{(6)}\;\;\; \displaystyle c_\Gamma(\tilde{t})D^{\tilde{G}}(\tilde{t})^{1/2}=c_{\Gamma_0}(\tilde{t}) D^{\tilde{G}_0}(\tilde{t})^{1/2} |2|_F^a\prod_{j=1}^u m_{geom}(\Gamma_j)$$

\noindent où $a=(d-d_0)/2$. Rappelons que

$$\displaystyle c_\Gamma(\tilde{t})=\frac{1}{|Nil(\mathfrak{u}(V''))_{reg}|}\sum_{\mathcal{O}\in Nil(\mathfrak{u}(V''))_{reg}} c_{\Gamma,\mathcal{O}}(\tilde{t})$$

\noindent et

$$\displaystyle c_{\Gamma_0}(\tilde{t})=\frac{1}{|Nil(\mathfrak{u}(V_0''))_{reg}|}\sum_{\mathcal{O}\in Nil(\mathfrak{u}(V_0''))_{reg}} c_{\Gamma_0,\mathcal{O}}(\tilde{t})$$

\noindent Soit $\mathcal{O}\in Nil(\mathfrak{u}(V''))_{reg}$. La proposition 1.12 de [W3] permet d'exprimer $c_{\Gamma,\mathcal{O}}(\tilde{t})$ en fonction de $\Gamma^{\tilde{L}}$. Explicitons la formule dans notre cas particulier. On reprend pour cela les notations de loc.cit. \\

 (7) L'ensemble $\mathcal{X}^{\tilde{L}}(\tilde{t})$ est réduit à un élément que l'on peut supposer être $\tilde{t}$. \\
 
 Rappelons que l'ensemble $\mathcal{X}^{\tilde{L}}(\tilde{t})$ est un système de représentants des classes de conjugaison par $L(F)$ des éléments dans $\tilde{L}(F)$ qui sont conjugués par $G(F)$ à $\tilde{t}$. Soit $g\in G(F)$ tel que $g\tilde{t}g^{-1}\in \tilde{L}(F)$. Tout comme dans la preuve du point (4), cela entraîne $g^{-1}A_{\tilde{L}}g\subset U(V'')$. On en déduit que les espaces $g^{-1}V_j$ sont inclus dans $V''$, sont totalement isotropes et que $g^{-1}V_j$ est en dualité avec $g^{-1}V_{-j}$. Or c'est aussi le cas de la famille de sous-espaces $(V_j)_{j=\pm1,\ldots,\pm u}$ de $V''$. Il existe un élément de $U(V'')(F)\subset G_{\tilde{t}}(F)$ qui envoie la première famille sur la deuxième. Donc quitte à multiplier $g$ à droite par un élément de $G_{\tilde{t}}(F)$, on peut supposer que $gV_j=V_j$ pour $j=\pm1,\ldots,\pm u$. Puisque $g\tilde{t}g^{-1}\in\tilde{L}(F)$, on a $g\tilde{t}g^{-1}V_j=V^*_{-j}=\tilde{t}(V_j)$. On en déduit que ${}^tg(V^*_j)=V^*_j$ pour $j=\pm 1,\ldots,\pm u$ puis par dualité que $g(V_0)=V_0$. L'élément $g$ conserve donc les sous-espaces $V_j$, $j=-u,\ldots,u$ donc appartient à $L(F)$, ce qui prouve le point (7). \\
 
 L'ensemble $\Gamma_{\tilde{t}}$ de 1.12 [W3] est le groupe $Z_G(\tilde{t})(F)=U(V'')\times T(F)^\theta=U(V'')\times T_\theta(F)=G_{\tilde{t}}(F)$. La deuxième somme ne comporte donc elle aussi qu'un seul terme que l'on peut supposer être $g=Id$. La dernière somme porte sur les orbites nilpotentes $\mathcal{O}'$ de $\mathfrak{l}_{\tilde{t}}(F)$ telles que $\mathcal{O}$ soit dans l'induite de l'orbite de $\mathcal{O}'$, condition que l'on écrit $[\mathcal{O}:\mathcal{O}']=1$. Une telle orbite est forcément régulière. Le morphisme qui à $g\in L_{\tilde{t}}$ associe sa restriction à $V_u\oplus\ldots\oplus V_1\oplus V_0$ définit un isomorphisme
 
$$\mbox{(8)}\;\;\; L_{\tilde{t}}\simeq R_{E/F} GL(V_u)\times\ldots\times R_{E/F} GL(V_1)\times G_{0,\tilde{t}}$$
$$=R_{E/F} GL_{d_u}\times\ldots\times R_{E/F} GL_{d_1}\times U(V_0'')\times T_\theta$$

\noindent d'où l'on déduit une identification $Nil(\mathfrak{l}_{\tilde{t}}(F))_{reg}=Nil(\mathfrak{gl}_E(V_u))_{reg}\times\ldots\times Nil(\mathfrak{gl}_E(V_1))_{reg}\times Nil(\mathfrak{u}(V_0''))_{reg}$. Les ensembles $Nil(\mathfrak{gl}_E(V_j))_{reg}$, $j=1,\ldots,u$ sont réduits à un élément que l'on note $\mathcal{O}_j$. On vérifie alors sur la description des orbites nilpotentes régulières dans les groupes unitaires, donnée en 3.2 [B], qu'il existe une unique orbite $\mathcal{O}_0\in Nil(\mathfrak{u}(V_0''))_{reg}$ de sorte que

$$[\mathcal{O}:\mathcal{O}_u\times\ldots\times\mathcal{O}_1\times\mathcal{O}_0]=1$$

\noindent Posons $\mathcal{O}'=\mathcal{O}_u\times\ldots\times\mathcal{O}_1\times\mathcal{O}_0$. Le lemme 1.12 de [W3] nous donne alors

$$c_{\Gamma,\mathcal{O}}(\tilde{t})=D^{\tilde{G}}(\tilde{t})^{-1/2}D^{\tilde{L}}(\tilde{t})^{1/2}c_{\Gamma^{\tilde{L}},\mathcal{O}'}(\tilde{t})$$

\noindent Par définition, le terme $c_{\Gamma^{\tilde{L}},\mathcal{O}'}(\tilde{t})$ se calcule à partir du développement local de $X\mapsto \Gamma^{\tilde{L}}(\tilde{t}exp(X))$ au voisinage de $0$ dans $\mathfrak{l}_{\tilde{t}}(F)$. Soit $X\in \mathfrak{l}_{\tilde{t}}(F)$ assez proche de $0$. Via l'isomorphisme (8), $X$ se décompose en $\displaystyle X_0+\sum_{j=1}^u X_j$ où $X_0\in\mathfrak{g}_{0,\tilde{t}}(F)$ et $X_j\in\mathfrak{gl}_{d_j}(E)$, $j=1,\ldots,u$. Posons $\tilde{x}=\tilde{t}exp(X)$. On a alors $\tilde{x}_0=\tilde{t}exp(X_0)$ et ${}^t(\tilde{x}_{-j}^c)^{-1}\tilde{x}_j=exp(2X_j)$. Par définition de $\Gamma^{\tilde{L}}$, on a

$$\displaystyle \Gamma^{\tilde{L}}(\tilde{t}exp(X))=\Gamma_0(\tilde{t}exp(X_0))\prod_{j=1}^u \Gamma_j(exp(2X_j))$$ 

\noindent D'après les propriétés d'homogénéités des fonctions $\hat{j}(\mathcal{O}_j,.)$, on en déduit que

$$\displaystyle c_{\Gamma^{\tilde{L}},\mathcal{O}'}(\tilde{t})=|2|_F^bc_{\Gamma_0,\mathcal{O}_0}(\tilde{t})\prod_{j=1}^u m_{geom}(\Gamma_j)$$

\noindent où $b=-\sum_{j=1}^u d_j(d_j-1)$. On vérifie, toujours à partir de la description 3.2 [B], que l'application qui à $\mathcal{O}$ associe $\mathcal{O}_0$ est une bijection entre $Nil(\mathfrak{u}(V''))_{reg}$ et $Nil(\mathfrak{u}(V_0''))_{reg}$. Par conséquent, on a

$$c_\Gamma(\tilde{t})=|2|_F^bD^{\tilde{G}}(\tilde{t})^{-1/2}D^{\tilde{L}}(\tilde{t})^{1/2}c_{\Gamma_0}(\tilde{t})\prod_{j=1}^u m_{geom}(\Gamma_j)$$

On vérifie facilement que

$$D^{\tilde{L}}(\tilde{t})=|2|_F^c D^{\tilde{G}_0}(\tilde{t})$$

\noindent où $c=2(d_1^2+\ldots+d_u^2)$. Puisque $b+\frac{c}{2}=a$, cela établit l'égalité (6). \\

Passons maintenant au cas $d_0>d_W$. Adoptons pour cela une notation bien commode. Soient $(\Theta_j)_{j=1,\ldots,t}$ une famille de quasi-caractères de $GL_{d'_j}(E)$ et $\Theta_0$ un quasi-caractère de $\tilde{G}_{d'_0}(F)$ où $\tilde{G}_{d'_0}$ est le groupe tordu définit de la même façon que $\tilde{G}$ en remplaçant $V$ par un espace de dimension $d'_0$. Alors on peut appliquer la construction faite au début de ce paragraphe pour obtenir un quasi-caractère de $\tilde{G}_{d'}(F)$ où $d'=d'_0+\ldots+d'_t$. On note $\Theta_0\times \Theta_1\times\ldots\times \Theta_t$ ce quasi-caractère. Il ne dépend pas des différents choix effectués, ce qui retire toutes ambiguïtés à la notation. En particulier, on a

$$\Gamma=\Gamma_0\times \Gamma_1\times\ldots\times \Gamma_u$$

\noindent Soient $\Theta'$ un quasi-caractère de $GL_{r+1}(E)$ et $\Gamma'$ un quasi-caractère de $GL_1(E)$ tels que $m_{geom}(\Theta')=m_{geom}(\Gamma')=1$ (il en existe). Puisque $\epsilon_{geom,\nu}(\Gamma,\Theta)=\epsilon_{geom,-\nu}(\Theta,\Gamma)$, d'après le premier cas, on a l'égalité

$$\epsilon_{geom,\nu}(\Gamma,\Theta)=\epsilon_{geom,-\nu}(\Theta\times \Theta',\Gamma)$$

\noindent Puisque $\epsilon_{geom,-\nu}(\Theta\times \Theta',\Gamma)=\epsilon_{geom,\nu}(\Gamma,\Theta\times \Theta')$, d'après le premier cas on a aussi l'égalité

$$\epsilon_{geom,-\nu}(\Theta\times \Theta',\Gamma)=\epsilon_{geom,\nu}(\Gamma\times \Gamma',\Theta\times \Theta')$$

\noindent Puisque $\Gamma\times \Gamma'=\Gamma_0\times \ldots\times \Gamma_u\times \Gamma'$, on a, toujours d'après le premier cas,

\[\begin{aligned}
\displaystyle \epsilon_{geom,\nu}(\Gamma\times \Gamma',\Theta\times \Theta') & =\epsilon_{geom,\nu}(\Gamma_0,\Theta\times\Theta')\prod_{j=1}^u m_{geom}(\Gamma_j) \\
 & =\epsilon_{geom,\nu}(\Gamma_0,\Theta)\prod_{j=1}^u m_{geom}(\Gamma_j)
\end{aligned}\]

\noindent Ce qui permet de conclure $\blacksquare$

\subsection{Preuve du théorème 6.1.1}

On démontre le théorème 6.1.1 par récurrence sur la dimension de $V$. On distingue deux cas pour la représentation $\tilde{\pi}$: le cas où $\tilde{\pi}$ est une induite parabolique tordue d'une représentation tempérée et le cas où $\tilde{\pi}$ est une représentation elliptique.

\subsubsection{Le cas proprement induit}
Supposons qu'il existe un parabolique tordu propre $\tilde{Q}=\tilde{L}U_Q$, $\tilde{Q}\neq\tilde{G}$, et une représentation $\tilde{\tau}\in Temp(\tilde{L})$ tels que $\tilde{\pi}=i_Q^G(\tilde{\tau})$. On écrit $\tilde{L}$ et $\tilde{\tau}$ comme en 2.3. En particulier, le prolongement de $\pi$ à $\tilde{G}(F)$ détermine un prolongement de $\tau_0$ à $\tilde{G}_0(F)$. D'après 2.5(2), on a

$$\displaystyle \mbox{(1)}\;\;\; \epsilon_\nu(\tilde{\pi},\tilde{\sigma})=\epsilon_\nu(\tilde{\tau}_0,\tilde{\sigma})\prod_{j=1}^u \omega_{\tau_j}(-1)^m$$

\noindent D'autre part, on a $\Theta_{\tilde{\pi}}=Ind_{\tilde{L}}^{\tilde{G}}(\Theta_{\tilde{\tau}})$. D'après le lemme 2.3.1, on a $\Theta_{\tilde{\tau}}=\Gamma^{\tilde{L}}$ avec les notations de la section précédente où on a posé $\Gamma_0=\Theta_{\tilde{\tau}_0}$ et $\Gamma_j=\omega_{\tau_j}(-1)^{d+1}\Theta_{\tau_j}$ pour $j=1,\ldots,u$. On déduit de la proposition 6.2.1 l'égalité

$$\displaystyle \mbox{(2)} \;\;\; \epsilon_{geom,\nu}(\tilde{\pi},\tilde{\sigma})=\epsilon_{geom,\nu}(\tilde{\tau}_0,\tilde{\sigma})\prod_{j=1}^u \omega_{\tau_j}(-1)^{d+1} m_{geom}(\Theta_{\tau_j})$$

\noindent D'après un résultat de Rodier ([Ro]), le terme $m_{geom}(\Theta_{\tau_j})$ vaut $1$ si $\tau_j$ admet un modèle de Whittaker et $0$ sinon. Comme $\tau_j$ est une représentation tempérée, elle admet un modèle de Whittaker. Donc $m_{geom}(\Theta_{\tau_j})=1$ pour $j=1;,\ldots,u$. On a aussi $(-1)^{d+1}=(-1)^m$ puisque $d$ et $m$ sont de parités différentes. Par hypothèse de récurrence appliquée au couple $(\tilde{G}_0,\tilde{H})$, on a $\epsilon_\nu(\tilde{\tau}_0,\tilde{\sigma})=\epsilon_{geom,\nu}(\tilde{\tau}_0,\tilde{\sigma})$. Comparant (1) et (2), on obtient l'égalité $\epsilon_\nu(\tilde{\pi},\tilde{\sigma})=\epsilon_{geom,\nu}(\tilde{\pi},\tilde{\sigma})$.

\subsubsection{Le cas des représentations elliptiques}

Supposons maintenant que $\tilde{\pi}\in \Pi_{ell}(\tilde{G})$. Soit $\tilde{f}\in C_c^\infty(\tilde{G}(F))$ une fonction très cuspidale. D'après les théorèmes 3.5.1 et 5.1.1, on a l'égalité

$$\mbox{(1)}\;\;\; J_{geom}(\Theta_{\tilde{\sigma}},\tilde{f})=J_{spec}(\Theta_{\tilde{\sigma}},\tilde{f})$$

\noindent Par définition, on a $J_{geom}(\Theta_{\tilde{\sigma}},\tilde{f})=\epsilon_{geom,\nu}(\Theta_{\tilde{f}},\Theta_{\tilde{\sigma}})$. Reprenons les notations du paragraphe 1.7, on a la décomposition

$$\displaystyle \mbox{(2)}\;\;\; \Theta_{\tilde{f}}=\sum_{\tilde{L}\in\mathcal{L}^{\tilde{G}}} |W^L||W^G|^{-1} (-1)^{a_{\tilde{L}}} Ind_{\tilde{L}}^{\tilde{G}}(I\Theta_{\phi_{\tilde{L}}(\tilde{f})})$$

\noindent où la fonction $\phi_{\tilde{L}}(\tilde{f})$ doit être interprétée comme la somme localement finie

$$\displaystyle \sum_{\zeta\in \mathcal{A}_{\tilde{L},F}} \mathbf{1}_{H_{\tilde{L}}=\zeta}\phi_{\tilde{L}}(\tilde{f})$$

\noindent D'après la proposition 1.9.2, on a aussi

$$\displaystyle \mbox{(3)}\;\;\; I\Theta_{\phi_{\tilde{L}}(\tilde{f})}=\sum_{\zeta\in \mathcal{A}_{\tilde{L},F}}\sum_{\mathcal{O}\in\{\Pi_{ell}(\tilde{L})\}} c(\mathcal{O})mes(i\mathcal{A}_{\tilde{L},F}^*)\Theta_{\tilde{\rho}^\vee}(\mathbf{1}_{H_{\tilde{L}}=\zeta}\phi_{\tilde{L}}(\tilde{f})) \mathbf{1}_{H_{\tilde{L}}=\zeta}\Theta_{\tilde{\rho}}$$

\noindent où $\tilde{\rho}\in\mathcal{O}$ est un point base. Tout comme dans le paragraphe précédent, on peut appliquer la proposition 6.2.1 pour calculer $\epsilon_{geom,\nu}(Ind_{\tilde{L}}^{\tilde{G}}(\mathbf{1}_{H_{\tilde{L}}=\zeta}\Theta_{\tilde{\rho}}),\Theta_{\tilde{\sigma}})$. Pour $\zeta\neq 0$, au moins  un des quasi-caractères $\Gamma_j$, $j=1,\ldots,u$, est nul au voisinage de $1$ ce qui entraîne $m_{geom}(\Gamma_j)=0$ et par conséquent le terme que l'on cherche à évaluer est nul. Seul le terme pour $\zeta=0$ est possiblement non nul. Posons $\tilde{f}_{\tilde{L}}=\mathbf{1}_{H_{\tilde{L}}=0}\phi_{\tilde{L}}(\tilde{f})$ et $\tilde{\pi}(\tilde{\rho})=i^G_Q(\tilde{\rho})$ ($\tilde{Q}\in\mathcal{P}(\tilde{L})$). D'après (3) et ce qui précède, on a

$$\displaystyle \mbox{(4)}\;\;\; \epsilon_{geom,\nu}(I\Theta_{\phi_{\tilde{L}}(\tilde{f})},\tilde{\sigma})=\sum_{\mathcal{O}\in\{\Pi_{ell}(\tilde{L})\}} mes(i\mathcal{A}_{\tilde{L},F}^*)c(\mathcal{O})\Theta_{\tilde{\rho}^\vee}(\tilde{f}_{\tilde{L}})\epsilon_{geom,\nu}(\tilde{\pi}(\tilde{\rho}),\tilde{\sigma})$$

\noindent Pour $\tilde{L}\in\mathcal{L}^{\tilde{G}}$ et $\mathcal{O}\in\{\Pi_{ell}(\tilde{L})\}$, posons

$$X_{geom}(\tilde{L},\mathcal{O})=mes(i\mathcal{A}_{\tilde{L},F}^*)c(\mathcal{O})\Theta_{\tilde{\rho}^\vee}(\tilde{f}_{\tilde{L}})\epsilon_{geom,\nu}(\tilde{\pi}(\tilde{\rho}),\tilde{\sigma})$$

\noindent D'après (2) et (4), on a

$$J_{geom}(\Theta_{\tilde{\sigma}},\tilde{f})=\sum_{\tilde{L}\in\mathcal{L}^{\tilde{G}}} |W^L||W^G|^{-1} (-1)^{a_{\tilde{L}}} \sum_{\mathcal{O}\in\{\Pi_{ell}(\tilde{L})\}} X_{geom}(\tilde{L},\mathcal{O})$$

D'autre part, d'après la définition 5.1, on a

$$\displaystyle J_{spec}(\Theta_{\tilde{\sigma}},\tilde{f})=\sum_{\tilde{L}\in\mathcal{L}^{\tilde{G}}} |W^L||W^G|^{-1} (-1)^{a_{\tilde{L}}} \sum_{\mathcal{O}\in\{\Pi_{ell}(\tilde{L})\}} X_{spec}(\tilde{L},\mathcal{O})$$

\noindent où on a posé

$$\displaystyle X_{spec}(\tilde{L},\tilde{O})=c(\mathcal{O})\epsilon_\nu(\tilde{\rho},\tilde{\sigma}) \int_{i\mathcal{A}_{\tilde{L},F}^*} J^{\tilde{G}}_{\tilde{L}}(\tilde{\rho}^\vee_{\lambda},\tilde{f})d\lambda$$

\noindent D'après la proposition 1.7.1, on a l'égalité

$$\displaystyle \int_{i\mathcal{A}_{\tilde{L},F}^*} J^{\tilde{G}}_{\tilde{L}}(\tilde{\rho}^\vee_{\lambda},\tilde{f})d\lambda=mes(i\mathcal{A}_{\tilde{L},F}^*)\Theta_{\tilde{\rho}^\vee}(\tilde{f}_{\tilde{L}})$$

\noindent D'après 2.5(2) et le cas proprement induit, si $\tilde{L}\neq\tilde{G}$, on a donc

$$X_{geom}(\tilde{L},\mathcal{O})=X_{spec}(\tilde{L},\mathcal{O})$$

\noindent L'égalité (1) peut donc se réécrire

$$\mbox{(5)}\;\;\; \displaystyle \sum_{\tilde{\pi}'\in\Pi_{ell}(\tilde{G})} c(\tilde{\pi}') \Theta_{(\tilde{\pi}')^\vee}(\tilde{f}) (\epsilon_{geom,\nu}(\tilde{\pi}',\tilde{\sigma})-\epsilon_\nu(\tilde{\pi}', \tilde{\sigma}))=0$$

\noindent D'après [W4], la théorie des pseudo-coefficients est valable pour les groupes tordus. En particulier, il existe une fonction cuspidale $\tilde{f}\in C_c^\infty(\tilde{G}(F))$ (et que l'on peut même prendre très cuspidale d'après le lemme 1.9.1) telle que

\begin{itemize}
\item $\Theta_{\tilde{\pi}^\vee}(\tilde{f})=1$;
\item $\Theta_{\tilde{\pi}'}(\tilde{f})=0$ pour toute $\tilde{\pi}'\in \Pi_{ell}(\tilde{G})$ telle que $\tilde{\pi}'\not\simeq \tilde{\pi}^\vee$.
\end{itemize}

En appliquant l'égalité (5) à un tel pseudo-coefficient, on obtient l'égalité $\epsilon_\nu(\tilde{\pi}, \tilde{\sigma})=\epsilon_{geom,\nu}(\tilde{\pi}, \tilde{\sigma})$.

\bigskip

{\bf Bibliographie}

\bigskip

[AGRS] A. Aizenbud, D. Gourevitch, S. Rallis, G. Schiffmann: {\it Multiplicity one theorems}, Ann. of Math. (2) 172 (2010),  no. 2, 1407-1434.

[A1] J. Arthur: {\it The trace formula in invariant form}, Annals of Math. 114 (1981), p.1-74

[A2] -----------: {\it Intertwining operators and residues I. Weighted characters}, J. Funct. Analysis 84 (1989), p. 19-84

[B] R. Beuzart-Plessis: {\it La conjecture locale de Gross-Prasad pour les représentations tempérées de groupes unitaires}, prépublication 2012

[C] L.Clozel: {\it Characters of nonconnected, reductive p-adic groups}, Canad. J. Math. 39 (1987), no. 1, 149-167

[GGP]  W. T. Gan, B. Gross, D. Prasad:{\it Symplectic local root numbers, central critical $L$-values and restriction problems in the representation theory of classical groups}, Astérisque 346 (2012)

[JPSS] H. Jacquet, I.I. Piatetskii-Shapiro, J. Shalika: {\it Rankin-Selberg convolutions},  Amer. J. Math.  105  (1983),  no. 2, 367-464

[Ro] F. Rodier: {\it Modèle de Whittaker et caractères de représentations}, in Non commutative harmonic analysis, J. Carmona, J. Dixmier, M. Vergne éd. Springer LN 466 (1981), p.151-171

[Sh] F. Shahidi: {\it On certain L-functions}, Amer. J. Math. 103 (1981), no. 2, 297-355

[W1] J.-L. Waldspurger: {\it Une formule intégrale reliée à la conjecture de Gross-Prasad}, Compos. Math. 146 (2010), no. 5, p. 1180-1290

[W2] -----------------------: {\it Une formule intégrale reliée à la conjecture locale de Gross-Prasad, $2^{eme}$ partie: extension aux représentations tempérées}, Astérisque 346 (2012)

[W3] -----------------------: {\it Calcul d'une valeur d'un facteur $\epsilon$ par une formule intégrale}, Asterique 347 (2012)

[W4] -----------------------: {\it La formule des traces locale tordue}, prépublication 2012

[W5] -----------------------: {\it La formule de Plancherel pour les groupes $p$-adiques, d'après Harish-Chandra}, J. of the Inst. of Math. Jussieu 2 (2003), p.235-333
\bigskip

Institut de mathématiques de Jussieu 2 place Jussieu 75005 Paris \\
 e-mail: rbeuzart@math.jussieu.fr
\end{document}